\documentclass[12pt]{amsart}
\usepackage{amssymb}
\usepackage[usenames]{color}

\newtheorem{theorem}{Theorem}

\newtheorem{proposition}{Proposition}
\newtheorem{corollary}{Corollary}

\newtheorem{remark}{Remark}
\newtheorem{definition}{Definition}
\newtheorem{lemma}{Lemma}

\begin{document}

{

\title[Isotopy classification of degree four polynomials]{Isotopy classification of Morse polynomials of degree four on ${\mathbb R}^2$}

\author{V.A.~Vassiliev}
\address{Weizmann Institute of Science, Rehovot, Israel}
 \email{vavassiliev@gmail.com}

\begin{abstract}
We introduce a system of  invariants of isotopy classes of Morse polynomials ${\mathbb R}^2 \to {\mathbb R}^1$, prove its completeness for polynomials of degrees $\leq 4$, calculate all  71 possible values of these invariants  for the case of degree 4, and realize them by concrete Morse polynomials. Also we count all 45460 classes of {\em strictly} Morse polynomials of degree four with the maximal possible number (nine) of real critical points.


\end{abstract} 

\thanks{This work was supported by the Absorption Center in Science of the Ministry of Immigration and Absorption of the State of Israel}

\maketitle

\section{Introduction}

\subsection{}
In  local bifurcation theory, {\em discriminants} and {\em caustics} (i.e., sets of functions with singular zero level sets and with non-Morse critical points, respectively) are considered in parallel and in close relation to each other, see e.g.  \cite{Kluwer}, \cite{AGLV2}, \cite{AVGZ}. These studies include the (also related and parallel) enumeration of complementary domains of functions without such features. 
On the contrary, the global algebraic version of the enumeration of complements of discriminants, i.e. the classification of real polynomials with non-singular  zero level sets, is very popular, especially  since the 16th Hilbert problem (see e.g. the papers \cite{klein}, \cite{Pe}, \cite{Ro}, \cite{DK}, \cite{FK}, \cite{Kha}, \cite{schlaf}, \cite{viro1}, \cite{viro2}, \cite{zeuten}), while a similar study of spaces of Morse polynomials is less developed. For some early steps in the latter direction, see 
\cite{BE} where topological types of real Morse polynomials in one variable were enumerated,  and 
\cite{ANN}, \cite{AN} containing some  estimates of the number of possible topological types of {\em strictly} Morse polynomials of degree four with negative principal homogeneous parts in ${\mathbb R}^2$ and maximal possible number of real critical points. See also \cite{Nic} where the numbers of topological types of strictly Morse functions with a given number of critical points on $S^2$ are computed; this implies upper bounds on the numbers of topological types of Morse polynomials with positive or negative principal homogeneous parts of a given degree on ${\mathbb R}^2$. In the corresponding local theory, all topological types of arbitrarily small Morse perturbations of real simple singularities were enumerated, see \cite{gor}.

In the following, we study the  classification of all (not necessarily strictly) Morse polynomials of a fixed degree on ${\mathbb R}^2$ up to the more restrictive {\em isotopy} equivalence relation: two polynomials are isotopic if they belong to the same connected component of the space of Morse polynomials. Although the numbers of such classes are obviously smaller than these of strictly Morse polynomials, the classification of them is more complicated in the same way as the classification of knots is more complicated than that of knot diagrams.

We define  invariants of isotopy classes in terms of the topology of zero sets of these polynomials in ${\mathbb C}^2$. Computing all possible values of these invariants is a combinatorial problem that can be solved algorithmically. These values also provide clues to the construction of polynomials realizing them.  We solve these problems explicitly in the case of polynomials of degree $\leq 4$, listing  all possible values of these invariants, presenting the polynomials having these values and proving the completeness of our invariants in this case.  Many  obtained examples of Morse polynomials with many critical points seem to be unexpected. 

We also enumerate the isotopy classes of strictly Morse polynomials of degree four with the maximal possible number (equal to 9) of real critical points.

It turns out that the list of isotopy classes of Morse polynomials of degree four with nine real critical points  almost coincides with the classification of degree four polynomials having two real critical points with the sum of Milnor numbers equal to 9. We explicitly realize all existing classes of such polynomials (i.e. the pairs of singularity classes of such critical points) and prove the non-existence of the others. This result is a real version of the list from  \cite{Jaw2} of possible decompositions of complex $X_9$ singularities into pairs of  simple singularities. These examples provide many obstructions to the equisingularity of bifurcation sets along $\mu= \mbox{const}$ varieties of real function singularities of class $X_9$. 
\medskip

By a {\em function singularity} we mean a germ of a smooth function at its critical point.

\subsection{Stating the problem}
\label{statpr}
\begin{definition} \rm
A polynomial ${\mathbb R}^2 \to {\mathbb R}$ is {\em Morse} if all its real critical points are Morse. It is {\em strictly Morse} if, in addition, all its critical values {\em at real critical points} are different.
\end{definition}

Let $d$ be a natural number. 
Consider the space of all degree $d$ polynomials ${\mathbb R}^2 \to {\mathbb R}$. In general, all critical points of such a polynomial are Morse, and its principal (of degree $d$) homogeneous part is {\it non-dis\-cri\-minantal}, i.e. it vanishes on $d$ different lines in ${\mathbb C}^2$. The set of polynomials satisfying these genericity conditions  is disconnected: the variety of non-generic polynomials separates it into connected components. In this paper we count these components for $d \leq 4$. 

For $d=2$ there are three such components characterized by the Morse index of the single critical point of a polynomial.

The analogous problem for $d=3$ essentially coincides with the enumeration of components of complements of caustics of $D_4$ singularities, see e.g. \cite{AVGZ}. 
Namely, the non-discriminantal principal homogeneous parts of degree three polynomials form two classes $D_4^+$ and $D_4^-$ consisting respectively of polynomials vanishing either on one or on three real lines. By  linear transformations of ${\mathbb R}^2$, they can be reduced to the normal forms $x^2y + y^3$ and $x^2y - y^3$, respectively. The space $[D_4^+]$ or $[D_4^-]$
of all polynomials with principal parts of one of these two classes contains the 3-dimensional affine subspace consisting of polynomials of the form
\begin{equation}
f_0 + \alpha x + \beta y + \gamma y^2 ,
\label{d4}
\end{equation}
where $f_0 $ is the corresponding normal form $x^2y + y^3$ or $x^2y - y^3$.
All other polynomials in the space $[D_4^+]$ or $[D_4^-]$ can be reduced to polynomials in this subspace by orientation preserving affine transformations and adding the constants. Therefore,   enumerating the components of the set of Morse polynomials in this entire space reduces to  studying the intersection of this set with the subspace (\ref{d4}).

\unitlength 1.0mm
\linethickness{0.4pt}
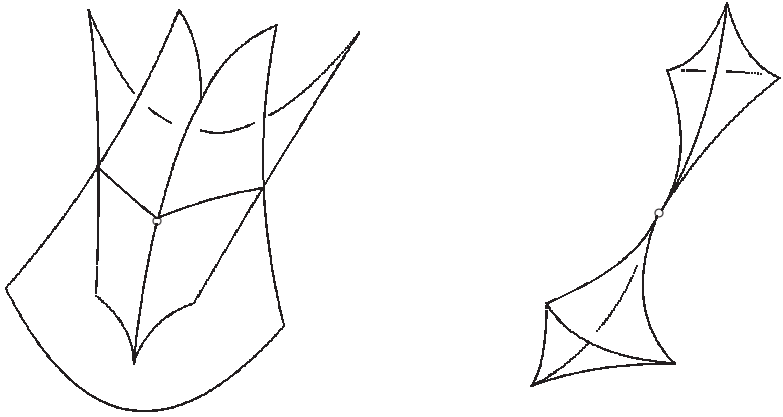
\begin{figure}
\begin{center}
\begin{picture}(47.00,54.00)
\bezier{400}(0.00,17.00)(15.00,-13.00)(37.00,12.00)
\bezier{88}(12.00,16.00)(17.00,12.00)(17.00,7.00)
\bezier{104}(17.00,7.00)(19.00,13.00)(25.00,15.00)
\bezier{328}(37.00,12.00)(32.00,32.00)(36.00,52.00)
\bezier{352}(0.00,17.00)(15.00,34.00)(23.00,54.00)
\bezier{304}(12.00,16.00)(13.00,41.00)(11.00,54.00)
\bezier{336}(25.00,15.00)(38.00,37.00)(47.00,51.00)
\bezier{104}(26.00,42.00)(26.00,50.00)(23.00,54.00)
\bezier{60}(26.00,42.00)(29.00,49.00)(36.00,52.00)
\put(20.10,26.00){\circle{1.00}}
\bezier{104}(19.70,26.20)(16.00,29.50)(12.20,33.20)
\bezier{128}(20.60,26.25)(27.10,29.15)(34.20,30.35)
\bezier{136}(20.30,26.40)(23.00,37.00)(26.00,42.00)
\bezier{152}(20.00,25.50)(18.00,17.00)(17.00,7.00)
\bezier{104}(17.00,43.00)(14.00,47.00)(11.00,54.00)
\bezier{128}(35.00,40.00)(40.00,42.50)(47.00,51.00)
\bezier{64}(33.00,39.00)(29.00,37.00)(26.00,38.00)
\bezier{32}(19.00,41.00)(20.00,40.00)(22.00,39.00)
\end{picture}
\qquad \qquad 
\begin{picture}(40.00,58.00)
\put(22.00,27.00){\circle{1.00}}
\bezier{168}(22.30,27.40)(27.00,35.00)(23.00,46.00)
\bezier{248}(22.30,27.40)(29.00,38.00)(31.00,55.00)
\bezier{192}(22.30,27.40)(29.00,38.00)(38.00,45.00)
\bezier{112}(23.00,46.00)(28.00,47.00)(31.00,55.00)
\bezier{112}(38.00,45.00)(33.00,47.00)(31.00,55.00)
\bezier{168}(21.70,26.60)(19.00,20.00)(7.00,15.00)
\bezier{168}(7.00,15.00)(11.00,8.00)(24.00,7.00)
\bezier{160}(24.00,7.00)(13.00,7.00)(5.00,4.00)
\bezier{96}(5.00,4.00)(7.00,7.00)(7.00,15.00)
\bezier{72}(5.00,4.00)(9.00,6.00)(12.50,9.60)
\bezier{192}(24.00,7.00)(17.00,14.00)(21.70,26.60)
\bezier{24}(25.00,46.00)(27.00,46.00)(28.00,45.90)
\bezier{32}(31.00,45.85)(32.00,45.84)(35.50,45.50)
\bezier{60}(13.7,11.3)(17,15)(19,20)
\end{picture}
\caption{Caustics for degree 3 polynomials: {\em purse} ($D_4^+$) and {\em pyramid} ($D_4^-$)}
\label{pyramid}
\end{center}
\end{figure}

The shapes of the caustics in the subspaces (\ref{d4}) are shown in two parts of Fig.~\ref{pyramid} (see \S I.3 in \cite{AVGZ} or \S II.2 in \cite{AGLV2}). In both cases the components of their complements are distinguished from each other by the {\em passport} invariant of the isotopy classes of the Morse functions, i.e. by the set \begin{equation}
\label{passp}
(m_-, m_\times, m_+)
\end{equation} of numbers of their 
minima, saddlepoints and maxima. For four components of the caustic complements of polynomials of type $D_4^+$ (see Fig.~\ref{pyramid} left), the passports are equal to $(1,  2, 1)$, $(1, 1, 0)$, $(0, 1, 1),$ and $(0, 0,0)$; for three components of such complements of type $D_4^-$, they are equal to $(1, 3, 0)$, $(0, 3, 1),$ and $(0, 2, 0)$. \medskip

Non-discriminantal principal homogeneous parts of real polynomials of degree four have four different types called  $X_9^+, X_9^-, X_9^1,$ and $X_9^2$:  positive $(X_9^+)$ or negative ($X_9^-$) quartic forms, and polynomials vanishing on exactly two $(X_9^1)$ or four $(X_9^2)$ different real lines. Each of these types is fibered into a 1-parametric family of orbits of the obvious  $\mbox{GL}(2, {\mathbb R})$-action. These orbits are characterized by the {\em $j$-invariant} (see e.g. \cite{Mukai}), which separates all $\mbox{GL}(2, {\mathbb C})$-orbits but sometimes takes equal values on different real orbits. The cases $X_9^+$ and $X_9^-$ can be reduced to each other by the multiplication by $-1$, therefore we will consider only the first of them. We will show 13 different isotopy classes of Morse functions with the principal part of type $X_9^+$, 29 classes of type $X_9^1$, and 16 classes of type $X_9^2$, and prove that these lists are complete. 

\subsection{Interpretation in the local theory of caustics}

The principal parts of these polynomials, i.e. non-discriminantal homogeneous polynomials of degree four, can also be considered as functions with isolated critical points at the origin. 
The space of all such polynomials with a given principal part $f_4$ is a deformation of the function singularity $f_4:$ the parameters of this deformation are the coefficients at lower monomials.  The property of being Morse is invariant under the addition of constants, so we can assume that the considered polynomials have zero constant term. If our polynomial $f$ splits into homogeneous summands as $f_4 + f_3 + f_2 + f_1$, then all functions from the one-parametric family $t^4f(t^{-1}x, t^{-1}y) \equiv f_4 + t f_3 + t^2 f_2 + t^3 f_1,$ $t >0$, have the same topological type, in particular they are simultaneously Morse or not. Therefore, we can assume that the considered polynomials $f$ are arbitrarily small perturbations of their principal parts $f_4$, and thus our study includes the enumeration of local components of complements of caustics of small deformations of function  singularities of class $X_9$. 

The similar problem for {\em simple} function singularities was solved in \cite{sed}, \cite{Vcau}. In particular, in these cases the passport invariant (\ref{passp})  is a complete invariant of the complements of the caustics. We will see below that this is not the case  for $X_9$ singularities. 
The analogous problem concerning Morse polynomials of degree three in ${\mathbb R}^3$ is considered in \cite{Vspace}. The analogous problem for the remaining class of parabolic singularities, $J_{10}$, is solved in \cite{VJ10}.

\subsection{Atypical adjacencies of multisingularities and equisingularity problem} \label{pr2} Our calculations show that the real caustics are not equisingular along the set of smooth functions with critical points of type $X_9$. We present explicit exceptional adjacencies of multisingularity strata which prevent such equisingularity. In the case of complex singularities, this effect has been considered in \cite{Jaw2}. 

The space of all polynomials of  degree four, whose principal parts are close to a given non-discriminantal homogeneous polynomial $f_4$, is  a {\em versal deformation} (see e.g. \cite{AVGZ}) of the corresponding function singularity. This implies that the local geometry of the caustic in this space determines the geometry of the caustic of the space of all smooth functions $M^2 \to {\mathbb R}$ in a neighborhood of any function having a singularity locally equivalent  to $f_4$ (i.e. reducible to $f_4$ by a choice of local coordinates and adding a constant function) and no other non-Morse critical points. For any particular polynomial $f_4$ of this class, the set of functions with a singularity locally equivalent to $f_4$ has codimension 8 in the space of all functions.  The set $\{X_9\}$ of functions with arbitrary $X_9$ singularities is swept out by the 1-parametric family of these sets (parametrized by the $j$-invariant) and thus has codimension 7.

\begin{table}
\caption{Real simple and $X_9$ singularities in two variables}
\label{t1}
\begin{center}
\begin{tabular}{|l|l|l|}
\hline
Notation & Normal form  & Restriction \\ 
\hline 
$A_{2k-1}$   & $\pm x^{2k} \pm y^2 $  &  $k \ge 1$ \\ 
$A_{2k}$  & $x^{2k+1} \pm y^2$  &  $k \ge 1$ \\
[3pt]
$D_{2k}^{\pm}$  & $x^2y \/\pm y^{2k-1} $  &  $k \ge 2$ \\ [4pt]
$D_{2k-1}$  & $\pm (x^2y \/+ y^{2k-2})$  & $k \ge 3$ \\ [3pt]
$E_6^{\pm}$ & $x^3 \pm y^4 $  & \cr
$E_7$ & $x^3 + x y^3 $ &  \cr
$E_8$ & $x^3 + y^5$ &  \cr
\hline
$X_9^{\pm}$ & $\pm(x^2 + 2A x^2 y^2 + y^4)$ & $|A|<1$ \\
$X_9^1$ & $x y (x^2 + 2A x y + y^2)$ & $|A|<1$ \\
$X_9^2$ & $\left\{\mbox{\begin{tabular}{l}
$x y (x^2 + 2A x y + y^2)$ \\
 $x^4 + 2 A x^2y^2 + y^4$ 
\end{tabular}}\right.$ 
& $\mbox{\begin{tabular}{l}
 $|A|>1$ \\
 $A<-1$
\end{tabular}
}
$
\\
\hline \end{tabular} \end{center}
\end{table}

Now, let $\Xi$ and $ \tilde \Xi$ be two {\em simple} singularity classes (see e.g. \cite{AVGZ}), the sum of whose Milnor numbers is 9. (The normal forms, to which simple and $X_9$ function singularities can be  reduced by a choice of local coordinates and adding constants, are listed in Table \ref{t1}; their Milnor numbers $\mu(\Xi)$ are the subscripts of the notations.)

The set $\{\Xi + \tilde \Xi\}$ of  smooth functions $M^2 \to {\mathbb R}$ with a critical point of class $\Xi$ and a critical point of class $\tilde \Xi$ (with arbitrary critical values) also has codimension 7 in the space of all smooth functions.  Therefore, most functions from the set $\{X_9\}$ are not approximated by functions of the set $\{ \Xi + \tilde \Xi\}$. However, some special functions of $\{X_9\}$ (corresponding to special values of the $j$-invariant) are approximated in this way. 

\begin{table}
\begin{center}
\caption{Approximation of $X_9$ singularities by bisingularities}
\label{tabadj}
\begin{tabular}{|c|l|l|l|}
\hline
Type & \qquad \qquad $X_9^+$ & \qquad \qquad $X_9^1$ & \qquad \qquad $ X_9^2$ \\
\hline
$E_7+A_2$ & No (index) & {\bf Yes}, \S~\ref{reale7a2}, $j=0$ & No (index) \\
$D_7 + A_2$ & No (Siersma) & No (Siersma) & No (Siersma) \\
$A_7 + A_2$ & {\bf Yes}, \S~\ref{reala7a2}, $j=\frac{5^37^3}{3^5}$ & No (D-graph) & No (index) \\
$E_6 + A_3$ & {\bf Yes,} \S \ref{+X9B}, $j=1$ & {\bf Yes,} \S \ref{adj38}, $j=1$ & No (index) \\
$A_6+A_3$ & No (D-graph) & {\bf Yes}, \S~\ref{reala6a3}, $j=-49$ & No (index) \\
$D_6^- + A_3$ & No (index) & No (D-graph) & {\bf Yes,} \S~\ref{reald6a3}, $ j=5$  \\
$D_6^+ + A_3$ & No $\left(\overline{\{D_6^+\}} \cap X_9=\emptyset \right)$ &  No $\left(\overline{\{D_6^+\}} \cap X_9=\emptyset \right)$ & No $\left(\overline{\{D_6^+\}} \cap X_9=\emptyset \right)$ \\
$D_5 + A_4$ & No (index) & {\bf Yes,} \S~\ref{reald5a4}, $j=-\frac{5}{3}$ & No (index) \\ 
$A_5+ A_4$ & {\bf Yes}, \S~\ref{reala5a4}, $j=\frac{245}{3}$ & {\bf Yes}, \S~\ref{reala5a4a}, $j=-96$ & No (index) \\
$D_5 + D_4^-$ & No (index, Bezout) & No (index, Bezout) & No (Bezout) \\
$D_5+ D_4^+$ & No (index, Bezout) & No (Bezout) & No (index, Bezout) \\
$A_5+D_4^-$ & No (index) & No (D-graph) & {\bf Yes}, \S~\ref{reala5d4m}, $j=\frac{25}{9}$ \\
$A_5+D_4^+$ & {\bf Yes}, \S~\ref{adj1}, $j=\frac{5^3}{3^3}$ & {\bf Yes}, \S~\ref{adj3}, $j=0$ & No (index) \\
\hline
\end{tabular}
\end{center}
\end{table}

\begin{definition}[cf. \cite{Jaw2}, \S 2] \rm
For any pair of simple singularity classes $\Xi$ and $\tilde \Xi$ with $\mu(\Xi) + \mu(\tilde \Xi) = 9$ and a class of functions $X_9^\ast$
$(\ast = +, 1$ or $2)$, the notation $\{\Xi + \tilde \Xi\}\rightsquigarrow \{X_9^*\}$ means that there exists a smooth function $f:{\mathbb R}^2 \to {\mathbb R}$ with an $X_9^*$ critical point, and its 1-parametric deformation $F: {\mathbb R}^2 \times [0, \varepsilon) \to {\mathbb R}$, $F(\cdot, 0) \equiv f$, such that for any $\tau \in (0, \varepsilon) $ the corresponding function $f_\tau \equiv F(\cdot, \tau)$ has a critical point of type $\Xi$ and a critical point of type $\tilde \Xi$ in such a way that these two critical points depend continuously on $\tau$ and tend to the $X_9^\ast$ singular point of the function $f$ as $\tau$ tends to 0.
\end{definition}

\begin{theorem}
\label{tabadjp}
For any pair of simple singularity classes $\Xi$ and $\tilde \Xi$ with $\mu(\Xi) + \mu(\tilde \Xi) = 9$ and any class of functions $X_9^\ast$
$(\ast = +, 1$ or $2)$, we have $\{\Xi + \tilde \Xi\}\rightsquigarrow X_9^*$
 if and only if ``Yes'' is written in Table~\ref{tabadj} at the intersection of the row $\Xi + \tilde \Xi$ and the  column $X_9^{\ast}$. The $j$-invariants of some functions of type $X_9^*$ which are approached by the strata  $\{\Xi + \tilde \Xi\}$ are indicated in the ``Yes'' cells.
\end{theorem}

\begin{corollary}
All existing types of splitting of complex $X_9$ singularities into pairs of simple singularities with $\mu(\Xi)+\mu(\tilde \Xi)=9$ listed in \cite{Jaw2} can be realized by splitting real $X_9^*$ singularities into pairs of real simple singularities.
\end{corollary}

All adjacencies mentioned in Table \ref{tabadj} are realized in the sections indicated in the table. The non-existence of the adjacencies in the cases when ``No'' stays in the table is proved in \S~\ref{adj2}.

It turns out that almost all examples of Morse polynomials of degree four with the maximal number of real critical points are closely related to these exceptional adjacencies of strata $\{\Xi+\tilde \Xi\}$. Curiously, the topological code (in the terms of which our invariants are formulated) of the unique example not related to such an adjacency  contains the {\em extended Coxeter-Dynkin graph} \ \ 
\unitlength 0.5mm
\begin{picture}(60,10)
\put(0,0){\line(1,0){60}}
\put(0,0){\circle*{1.5}}\put(10,0){\circle*{1.5}}\put(20,0){\circle*{1.5}}\put(30,0){\circle*{1.5}}\put(40,0){\circle*{1.5}}\put(50,0){\circle*{1.5}}\put(60,0){\circle*{1.5}}\put(30,10){\circle*{1.5}}
\put(30,0){\line(0,1){10}}
\end{picture} \ 
of $\tilde E_7$ type (which is one of the alternative names of the $X_9$ singularity, see e.g. \cite{bo} p. 199 or \cite{siersma} p. 86). Moreover,
the same situation (with the extended graphs of types $\tilde E_6$ and $\tilde E_8$) occurs for Morse polynomials of degree three in ${\mathbb R}^3$, which are closely related to the singularity class $P_8$,  and for $J_{10}$ singularities; see \cite{Vspace} and \cite{VJ10}.

\subsection{About the methods}

The main tool of this work is the so-called {\em Lyashko--Looijenga map} (see e.g. \cite{Lo0} which allows one to implement the deformations of critical values of functions by deformations of these functions themselves, and thus to construct (or prove the non-existence of) the functions with different restrictions on their critical points and critical values. Using this technique,  O.V.~Lyashko \cite{lyashko} enumerated all possible decompositions of simple complex singularities. In the application to real functions, it was then used in \cite{Izv} for the construction and enumeration of non-discriminantal perturbations of simple singularities with some prescribed topological properties motivated by the theory of hyperbolic PDE. An important part of this investigation is the construction of singularity decompositions with only two critical values (called in \cite{Izv} {\em characterizing functions}). Analogous calculations for some non-simple singularities (see \cite{umn}, \cite{AGLV2}) were done with the help of a computer program enumerating possible topological types of morsifications of real functions. The restriction of  the Lyashko--Looijenga map to the real algebraic functions was then used in \cite{chis} for the regular study of decompositions of simple singularities, in \cite{gor} and some related works. 

A large part of the computations in the present paper is also done by the above mentioned computer program; the assertions about the completeness of the obtained lists are largely based on the results of P.~Jaworski \cite{Jaw}, \cite{Jaw2}, who extended the nice properties of the Lyashko--Looijenga map (originally established only for simple singularities) to parabolic function singularities.  

\section{Invariants of isotopy classes of Morse polynomials}

\subsection{Trivial invariant} 
\label{itriv}

The simplest invariant of isotopy classes of Morse polynomials in ${\mathbb R}^2$ is their passport (\ref{passp}). 
In the case of degree three polynomials this is a complete invariant, see \S \ref{statpr}.

By the index considerations, for any degree four Morse polynomial with non-discriminantal principal part the sum of three numbers of the passport is an odd number not greater than 9, and the {\em Euler number} $m_- - m_{\times} + m_+$ is  1 for $ X_9^+$ and $X_9^-$ polynomials,  $-1$ for $X_9^1$, and  $-3$ for $X_9^2$.  Therefore, when studying polynomials of one of these three types we will characterize the triplets $(m_-, m_\times, m_+)$ only by pairs of numbers: the total number $M= m_- + m_\times + m_+$ of real critical points and the number $m_+$ of maxima.

\begin{table}
\caption{Numbers of classes for $X_9^+$ (left), $X_9^1$ (center) and $X_9^2$ (right)}
\label{TX+}
\begin{tabular}{|c| c | c | c | c | c |}
\hline
$m_+ \backslash M$ & 1 & 3 & 5 & 7 & 9 \\
\hline
0 & 1 & 1 & 1 & 2 & 4 \\
1 &   & 1 & 1 & 1 & 1 \\
\hline
\end{tabular}
 \ \
\begin{tabular}{|c| c | c | c | c | c |}
\hline
$m_+ \backslash M$ & 1 & 3 & 5 & 7 & 9 \\
\hline
0 & 1 & 1 & 1 & 1 & 8 \\
1 &   & 1 & 1 & 1 & 1 \\
2 &   &   & 1 & 1 & 1 \\
3 &   &   &   & 1 & 1 \\
4 &   &   &   &   & 8 \\
\hline 
\end{tabular} \ \ 
\begin{tabular}{|c|| c |c | c | c |}
\hline
$m_+ \backslash M$ & 3 & 5 & 7 & 9 \\
\hline
0 & 1 & 1 & 1 & 4 \\
1 &  & 1 & 1 & 1 \\
2 &   &   & 1 & 1 \\
3 &   &   &   & 4 \\
\hline 
\end{tabular}

\end{table}

In Table \ref{TX+},  
we list all pairs $(M, m_+)$ which (according to the results formulated in \S \ref{mthms} below) can be realized by Morse polynomials with principal parts of types $X_9^+$, $X_9^1$ and $X_9^2$. Namely, we give there the numbers of isotopy classes of these realizations with given pairs $(M, m_+)$. The analogous table for $X_9^-$ can be obtained from the table for $X_9^+$ by replacing the numbers 0 and 1 in the left column by 5 and 4, respectively. 

\subsection{Set-valued invariant and virtual Morse functions}

\label{svinv}

In this subsection, we describe a much stronger combinatorial invariant which in particular separates all isotopy classes of degree four Morse polynomials with non-discriminantal principal parts up to reflections in ${\mathbb R}^2$. 
\smallskip

First, with any {\em generic} Morse polynomial $f$  we associate a set of its discrete topological invariants, called a {\em virtual Morse function}. 

\begin{definition} \rm
A polynomial $f: ({\mathbb C}^2, {\mathbb R}^2) \to ({\mathbb C}, {\mathbb R})$  of degree $d$ is {\em generic} if it has only Morse critical points in ${\mathbb C}^2$, all of their critical values are different and not equal to 0, and the zero set of the principal homogeneous part of $f$ consists of $d$ different lines  in ${\mathbb C}^2$.
\end{definition}

\unitlength 1.00mm
\linethickness{0.4pt}
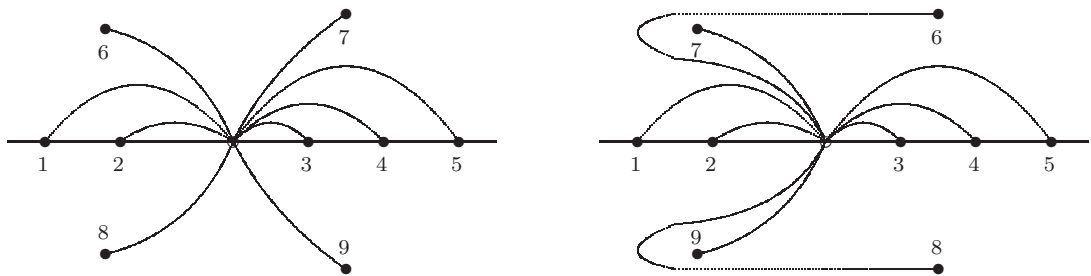
\begin{figure}
\begin{center}
\begin{picture}(65,40)
\put(0,20){\line(1,0){65}}
\put(5,20){\circle*{1.5}}
\put(15,20){\circle*{1.5}}
\put(40,20){\circle*{1.5}}
\put(50,20){\circle*{1.5}}
\put(60,20){\circle*{1.5}}
\bezier{100}(30,20)(17,35)(5,20)
\bezier{50}(30,20)(22,25)(15,20)
\bezier{50}(30,20)(35,25)(40,20)
\bezier{90}(30,20)(40,30)(50,20)
\bezier{140}(30,20)(45,40)(60,20)
\put(13,35){\circle*{1.5}}
\put(13,5){\circle*{1.5}}
\put(45,37){\circle*{1.5}}
\put(45,3){\circle*{1.5}}
\bezier{100}(30,20)(25,32)(13,35)
\bezier{100}(30,20)(25,8)(13,5)
\bezier{100}(30,20)(35,30)(45,37)
\bezier{100}(30,20)(35,10)(45,3)
\put(30,20){\circle{1.5}}
\put(4,16){{\tiny 1}}
\put(14,16){{\tiny 2}}
\put(39,16){{\tiny 3}}
\put(49,16){{\tiny 4}}
\put(59,16){{\tiny 5}}
\put(12,31){{\tiny 6}}
\put(44,33){{\tiny 7}}
\put(12,7){{\tiny 8}}
\put(44,5){{\tiny 9}}
\end{picture} \qquad \quad
\begin{picture}(65,40)
\put(0,20){\line(1,0){65}}
\put(5,20){\circle*{1.5}}
\put(15,20){\circle*{1.5}}
\put(40,20){\circle*{1.5}}
\put(50,20){\circle*{1.5}}
\put(60,20){\circle*{1.5}}
\bezier{100}(30,20)(17,35)(5,20)
\bezier{50}(30,20)(22,25)(15,20)
\bezier{50}(30,20)(35,25)(40,20)
\bezier{90}(30,20)(40,30)(50,20)
\bezier{140}(30,20)(45,40)(60,20)
\put(13,35){\circle*{1.5}}
\put(13,5){\circle*{1.5}}
\put(45,37){\circle*{1.5}}
\put(45,3){\circle*{1.5}}
\bezier{100}(30,20)(25,32)(13,35)
\bezier{100}(30,20)(25,8)(13,5)
\bezier{100}(30,20)(25,30)(10,31)
\bezier{50}(10,31)(0,35)(10,37)
\bezier{100}(10,37)(35,37)(45,37)
\bezier{100}(30,20)(25,10)(10,9)
\bezier{50}(10,9)(0,5)(10,3)
\bezier{100}(10,3)(35,3)(45,3)
\put(30,20){\circle{1.5}}
\put(4,16){{\tiny 1}}
\put(14,16){{\tiny 2}}
\put(39,16){{\tiny 3}}
\put(49,16){{\tiny 4}}
\put(59,16){{\tiny 5}}
\put(12,31.3){{\tiny 7}}
\put(44,33){{\tiny 6}}
\put(12,6.4){{\tiny 9}}
\put(44,5){{\tiny 8}}
\end{picture}
\end{center}
\caption{Standard systems of paths}
\label{standd}
\end{figure}

If $f$ is a generic polynomial of degree $d$, then the set $V_f \subset {\mathbb C}^3$ defined by the equation  $f(x, y) + z^2=0$ is a smooth complex surface homotopy equivalent to the wedge of $(d-1)^2$ two-dimensional  spheres (see e.g. \cite{M}, \cite{AVGZ}). The homology group $H_2(V_f)$ is generated by {\em vanishing cycles} (see e.g. \cite{AVGZ}, \cite{APLT}) defined by some system of non-intersecting paths in ${\mathbb C}^1$ connecting the non-critical value 0 with all critical values of $f$, see Fig.~\ref{standd}. We agree to choose these paths so that the paths going to complex conjugate non-real critical values are symmetric to each other about the real axis, and the paths going to real values lie in the upper half-plane in the area where the imaginary part is smaller than the absolute values of the imaginary parts of all non-real critical values. 

Let us somehow fix an orientation of ${\mathbb R}^3$. The orientations of all vanishing cycles can then be defined by such a system of paths  in a standard way, see \S V.1.6 of \cite{APLT}. In particular, it is required that the complex conjugation in ${\mathbb C}^3$ takes the oriented vanishing cycles defined by complex conjugate paths into each other. Also, a canonical order of these vanishing cycles can be defined, in particular the cycles vanishing at real critical points go first in the ascending order of the corresponding critical values.

\begin{definition}[see \cite{AGLV2}, \S V.3]\rm
\label{vm1}

A {\em virtual Morse function associated with} a generic degree $d$  Morse {\em polynomial $f:({\mathbb C}^2, {\mathbb R}^2) \to ({\mathbb C}, {\mathbb R})$} is a collection of its topological data consisting of

a) the $(d-1)^2  \times (d-1)^2$ matrix of intersection indices in  $V_f$
of canonically ordered and oriented vanishing cycles $\Delta_i \in H_2(V_f)$ corresponding to all critical values of $f$ and defined by a system of paths as above,

b) the string of $(d-1)^2$  intersection indices in $V_f$ of these vanishing cycles with the naturally oriented set of real points, 
$V_f \cap {\mathbb R}^3$;

c) the string of Morse indices of all real critical points of $f$, ordered by the increase of their critical values, and

d) the numbers of negative and positive real critical values of $f$.
\end{definition}

\begin{table}
\caption{Two virtual  Morse functions with 7 real critical values}
\label{virtmatr}
\begin{tabular}{|cccc|ccc|cc|}
\hline
-2  & 0  & 0  & 0  & 1  & 0  & 1  & 1  & 1 \\
0 & -2 & 0 & 0  & 0  & 1 &  0  & 0  & 0 \\
 0  & 0 & -2  & 0  & 1 &  0 &  0 &  0 &  0 \\
 0  & 0  & 0 & -2  & 0  & 1  & 1  & 1 &  1 \\
 1 &  0 &  1 &  0 & -2 & 0 & 0 & -1 & -1 \\
 0 & 1 & 0  & 1 & 0 & -2 & 0 & -1 & -1 \\
 1  & 0  & 0 & 1 & 0 & 0 & -2 & -1 & -1 \\
 1  & 0 &  0  & 1 & -1 & -1 & -1&  -2 & -2 \\
 1  & 0 &  0 & 1 & -1 & -1 & -1 & -2 & -2 \\
\hline
2 & 2 & 2 & 2 & -2 & -2 & -2 & -2 & -2 \\
\hline
0 & 0 & 0 & 0 & 1 & 1 & 1 & & \\
\hline
\end{tabular} \ \quad
\begin{tabular}{|cccc|ccc|cc|}
\hline
-2  & 0 &  0 & 0 & 1 & 1 & 1 & 1 & 1 \\
 0 & -2 & 0 &  0 &  0 & 1 & 0 & 0 &  0 \\
 0 & 0 & -2  & 0  & 1 & 0 & 0  & 0 & 0 \\
 0  & 0 & 0 & -2 & 0 & 0 &  1 & 1 & 1 \\
 1 & 0 & 1 & 0 &-2 & 0 & 0 & -1 & -1 \\
 1 & 1 & 0 & 0 & 0 & -2 & 0 & 0 & 0 \\
 1 & 0  & 0 &  1&  0 & 0 & -2 & -1 & -1 \\
 1 & 0 & 0  & 1 & -1 & 0 & -1 & -2 & -2 \\
 1 & 0 & 0  & 1 & -1 & 0 & -1 & -2 & -2 \\
\hline
2 & 2 & 2 & 2 & -2 & -2 & -2 & -2 & -2 \\
\hline
0 & 0 & 0 & 0 & 1 & 1 & 1 & & \\
\hline
\end{tabular}
\end{table}

\noindent
{\bf Example.} 
Two virtual Morse functions with $d=4$ and seven real critical points are shown in two parts of Table \ref{virtmatr}. Two vertical lines in each of these tables indicate the last element of the corresponding virtual morsification: in both cases, the numbers of negative, positive and non-real critical values are 4, 3, and 2, respectively.

\begin{definition} \rm 
A {\em critical point of a virtual Morse function} is any column of its data set as in Table \ref{virtmatr}, i.e. a column of the intersection matrix, the intersection index with the real set, and a Morse index or the information that the point is non-real.
\end{definition}

\begin{remark} \rm
If the number of pairs of non-real critical values of $f$ is greater than one then there can be more than one virtual Morse function associated with the same real polynomial, because the choice of a system of paths is not homotopically unique: see Fig.~\ref{standd}.
\end{remark}

\begin{definition} \rm
\label{elsur}
{\em Elementary virtual surgeries}\index{elementary surgery} of virtual Morse functions include six transformations of their data, modeling the standard local topological surgeries of the corresponding real Morse polynomials, namely
\begin{itemize}
\item[$s1, s2$:] \ \ \ collision of two neighboring real critical values at a non-zero value, after which the corresponding two critical points either ($s1$) meet and go into the complex domain, or (s2) change the order in ${\mathbb R}^1$ of their critical values; 

\item[$s3, s4$:] \ \ \ collision of two complex conjugate critical values at a point of the line ${\mathbb R}^1 \setminus \{0\}$, after which the corresponding critical points either ($s3$) meet at a real point and come to the real space, or ($s4$) miss each other in the complex domain, while the imaginary parts of their critical values change their signs; 

\item[$s5, s6$:] \ \ \ jumps of real critical values up ($s5$) or down ($s6$) through 0; \\ \hspace*{4cm} and additionally

\item[$s7$:] \ \ \ specifically virtual transformations within the classes of virtual Morse functions associated with the same real Morse polynomials, caused by flips of standard systems of paths going from 0 to non-real critical values (see Fig.~\ref{standd} and \cite{AVGZ}, Figs. 19--21).
\end{itemize}
\end{definition}

The results of all these virtual surgeries are determined by the data of the original virtual Morse functions. For a detailed description of these standard flips of data, see  \S V.8 of \cite{APLT}, the explicit formulas for them are given in the comments to our program, see the reference in \S~\ref{onproofs}.
. In particular, an attempt to perform the surgery $s1$ or $s2$ over real critical values $v_i$ and $v_{i+1}$ begins by examining the intersection index $\langle \Delta_i, \Delta_{i+1}\rangle$ of the corresponding vanishing cycles. If this 
index is $0$ then the surgery $s2$ happens, if it is 1 then the surgery $s1$, in all other cases the surgery fails. Similarly, a collision of two complex conjugate critical values at a real point not separated from 0 by other critical values follows scenario $s4$ if the intersection index is 0, scenario $s3$ if it is $1$ or $-1$, and fails in all other cases; in the second  case the sign of this intersection index allows us to predict  the Morse indices of the newborn real critical points.

We will denote by $s1, \dots, s6$ both the real surgeries of real Morse functions and the corresponding elementary virtual surgeries.
\medskip

Let $f: {\mathbb R}^2 \to {\mathbb R}$ be a generic Morse polynomial of degree $d$.

\begin{definition} \rm
An (abstract) {\em virtual Morse function} of type $f$  is any collection of data as in Definition \ref{vm1} (i.e. a matrix, two strings and two numbers) obtained from an arbitrary virtual Morse function associated with $f$ by an arbitrary chain of elementary virtual surgeries.  

The {\em formal graph} of type $f$ is the graph, whose vertices are all virtual  Morse functions of type $f$, and two such vertices are connected by an edge if and only if they are obtained from each other by an elementary virtual surgery.

A {\em virtual component} $S(f)$ of the formal graph of type $f$  is its subgraph, whose vertices are only those virtual Morse functions of type $f$ which can be obtained from virtual Morse functions associated with $f$ by arbitrary chains of virtual surgeries $s2, s4, s5, s6$ and $s7$ from Definition \ref{elsur} (i.e. all surgeries which do not model the collision of critical points).
\end{definition}

\begin{proposition}
\label{protrivi}
If the principal homogeneous parts of generic Morse polynomials $f, \tilde f$ of the same degree are topologically equivalent $($i.e. they vanish on the same non-zero number of real lines, or are simultaneously negative or positive definite$)$,
then the formal graphs of types $f$ and $\tilde f$ are the same.

If two generic Morse polynomials $f$ and $\tilde f$ belong to the same connected component of the set of Morse  polynomials with non-discriminantal principal parts, then their virtual components $S(f)$, $S(\tilde f)$ of the formal graph are the same.
\end{proposition}

\noindent
{\it Proof.} The first statement follows from the connectedness of the space of all polynomials with principal parts of the same topological class: any two real Morse polynomials from this space can be connected by a path in it, which intersects the variety of non-Morse polynomials at finitely many points corresponding to standard surgeries. The second statement follows immediately from the definitions. \hfill $\Box$
\medskip

\begin{figure}
\begin{picture}(60,55)
\put(10,5){\circle*{2}}\put(10,15){\circle*{2}}\put(10,25){\circle*{2}}\put(10,35){\circle*{2}}\put(10,45){\circle*{2}}
\put(30,15){\circle*{2}}\put(30,25){\circle*{2}}\put(30,35){\circle*{2}}\put(50,5){\circle*{2}}\put(50,15){\circle*{2}}\put(50,25){\circle*{2}}\put(50,35){\circle*{2}}\put(50,45){\circle*{2}}
\put(10,5){\line(0,1){40}}
\put(30,15){\line(0,1){20}}
\put(50,5){\line(0,1){40}}
\put(10,5){\line(2,1){20}}
\put(10,25){\line(2,-1){20}}
\put(10,35){\line(2,-1){20}}
\put(10,45){\line(2,-1){20}}
\put(50,5){\line(-2,1){20}}
\put(50,15){\line(-2,1){20}}
\put(50,25){\line(-2,1){20}}
\put(50,45){\line(-2,-1){20}}
\bezier{150}(10,5)(18,7)(18,5)
\bezier{150}(10,5)(18,3)(18,5)

\bezier{150}(10,15)(18,17)(18,15)
\bezier{150}(10,15)(18,13)(18,15)

\bezier{150}(10,25)(18,27)(18,25)
\bezier{150}(10,25)(18,23)(18,25)

\bezier{150}(10,45)(18,47)(18,45)
\bezier{150}(10,45)(18,43)(18,45)

\bezier{150}(30,15)(38,17)(38,15)
\bezier{150}(30,15)(38,13)(38,15)

\bezier{150}(50,5)(58,7)(58,5)
\bezier{150}(50,5)(58,3)(58,5)

\bezier{150}(50,15)(58,17)(58,15)
\bezier{150}(50,15)(58,13)(58,15)

\bezier{150}(50,35)(58,37)(58,35)
\bezier{150}(50,35)(58,33)(58,35)

\bezier{150}(50,45)(58,47)(58,45)
\bezier{150}(50,45)(58,43)(58,45)

\bezier{150}(10,5)(2,7)(2,5)
\bezier{150}(10,5)(2,3)(2,5)

\bezier{150}(10,15)(2,17)(2,15)
\bezier{150}(10,15)(2,13)(2,15)

\bezier{150}(10,35)(2,37)(2,35)
\bezier{150}(10,35)(2,33)(2,35)

\bezier{150}(10,45)(2,47)(2,45)
\bezier{150}(10,45)(2,43)(2,45)

\bezier{150}(30,35)(22,37)(22,35)
\bezier{150}(30,35)(22,33)(22,35)

\bezier{150}(50,5)(42,7)(42,5)
\bezier{150}(50,5)(42,3)(42,5)

\bezier{150}(50,25)(42,27)(42,25)
\bezier{150}(50,25)(42,23)(42,25)

\bezier{150}(50,35)(42,37)(42,35)
\bezier{150}(50,35)(42,33)(42,35)

\bezier{150}(50,45)(42,47)(42,45)
\bezier{150}(50,45)(42,43)(42,45)

\bezier{150}(30,15)(32,7)(30,7)
\bezier{150}(30,15)(28,7)(30,7)

\bezier{150}(30,35)(32,42)(30,42)
\bezier{150}(30,35)(28,42)(30,42)

\bezier{150}(30,25)(32,32)(30.5,32)
\bezier{150}(30,25)(28,32)(29.5,32)
\end{picture}

\caption{Formal graph of $D_4^-$ singularity}
\label{graph}
\end{figure}
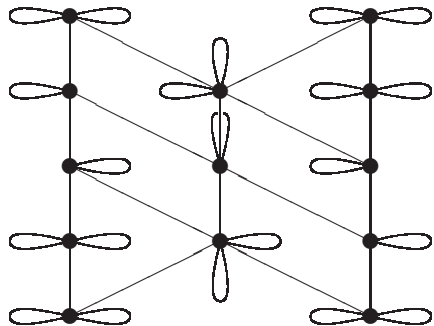

\noindent
{\bf Example.}
The formal graph of the $D_4^-$ singularity is shown in Fig.~\ref{graph}:
vertical segments (respectively, oblique segments, ``horizontal'' loops, and ``vertical'' loops)
in it denote surgeries of type $s5/s6$ (respectively,  $s1/s3$,  $s2$, and $s4$). (In the case of more complicated singularities the loop edges are very rare.) Removing the oblique segments splits this graph into three virtual components corresponding to three domains in Fig.~\ref{pyramid}.
\medskip

\begin{remark} \rm
In the parallel study of complements of discriminants, we define virtual components by removing the edges corresponding to surgeries $s5$ and $s6$ represented by the vertical segments in Fig.~\ref{graph}. This splits the formal graph $D_4^-$ into five components. One of these components (containing the central node of our graph) corresponds to three different components of the complement of the discriminant, so there are seven of them.
\end{remark}

Proposition \ref{protrivi} justifies the following definitions.

\begin{definition} \rm
The {\em formal graph} of a (not necessarily Morse) degree $d$ polynomial ${\mathbb R}^2 \to {\mathbb R}$ with non-discriminantal homogeneous principal part is the formal graph of an arbitrary generic Morse polynomial with the principal part of the same topological type. 

The {\it virtual component} of a (not necessarily strictly) Morse polynomial with non-discriminantal principal part is the virtual component of an arbitrary generic polynomial from the same connected component of the space of Morse polynomials with non-discriminantal principal parts.
\end{definition}

In the case $d=4$, Proposition \ref{protrivi} has a much stronger version.

\begin{proposition}
\label{propmain}
Let  $f$ be a generic polynomial of degree four, then

a$)$ every virtual Morse function of type $f$ is associated with a generic real polynomial $\tilde f$ whose principal homogeneous part is topologically equivalent to that of $f$;

b$)$ every virtual Morse function from the virtual component of $f$ is associated with some generic real polynomial from the same connected component as $f$ of the space of Morse polynomials with non-discriminantal principal parts.
\end{proposition}

\noindent
{\it Proof}. We need to show that any virtual surgery of a virtual Morse function associated with $f$ can be realized by a real surgery of $f$ occurring along a path in the space of polynomials with non-discriminantal principal parts.
The realization of the surgeries $s5$ and $s6$ is achieved by adding constants to $f$, and surgery $s7$ is completely virtual and does not require any change of $f$ at all.
Therefore, we only need to realize any surgery of types $s1$---$s4$ related with a collision of critical values. 

By a choice of affine coordinates in ${\mathbb R}^2$, we can assume that $f$ belongs to the standard miniversal deformation of functions of the corresponding class $X_9^*$, i.e. to the family of polynomials
\begin{equation} f_A(x) + \lambda_1 + \lambda_2 x + \lambda_3 y + \lambda_4 x^2 + \lambda_5 x y + \lambda_6 y^2 + \lambda_7 x^2 y + \lambda_8 x y^2 ,
\label{vers0}
\end{equation}
where $f_A$ is the normal form from Table \ref{t1}
(in the case $X_9^2$ we take the second normal form); here $A$ and $\lambda_1, \dots \lambda_8$ are the {\em complex} parameters of the deformation, $A$ takes all values except for $\pm 1$ which correspond to discriminantal homogeneous polynomials. Denote by $\Theta$ the parameter space $\left({\mathbb C}^1 \setminus \{1, -1\} \right) \times {\mathbb C}^8$ of this family.

Our construction is based on the results of  \cite{Jaw2} and exploits the {\it Lyashko--Looijenga map} $\Lambda: \Theta \to \mbox{Sym}^9({\mathbb C}^1) \equiv {\mathbb C}^9 / S(9)$ which sends any polynomial (\ref{vers0}) with a non-dis\-cri\-mi\-nantal principal part to the unordered collection of its critical values. This map is a submersion over the space $B({\mathbb C^1}, 9)$ of all collections of nine distinct points in ${\mathbb C}^1$, i.e. in the restriction to the space of all strictly Morse complex polynomials (\ref{vers0}) (cf. \cite{Jaw}). Indeed,  in a neighborhood of the line $\{\lambda_1 = \dots = \lambda_8 =0\} \cap \Theta $ this follows from the fact that some neighborhood of any point $(A, 0, \dots, 0)$ of this line is the parameter space of a versal deformation of the corresponding singularity of the polynomial $f_{A,0}$ of type $X_9$; the submersivity in it follows directly from the versality property. The one-parametric group of quasihomogeneous dilations   \begin{equation}
\label{dil}
f_{A,\lambda}(x, y) \mapsto t^{-4}f_{A,\lambda} (t x, t y)
\end{equation}
 extends this property to the whole $\Theta$.

Consider first the case of polynomials with principal parts of types $X_9^+$ and $X_9^2$.
Let  $I:[0,1] \to \mbox{Sym}^9({\mathbb C}^1)$ be a piecewise algebraic path such that 
\begin{itemize}
\item
$I(0)$ is the collection of critical values of $f$, 
\item
$I(\tau) \in B({\mathbb C}, 9)$ for $\tau < 1$, 
\item some seven points of the configuration $I(0)$ remain fixed for all $\tau \in [0,1]$, and two points of  $I(\tau)$ collide as in the prescribed surgery when $\tau$ tends to $1$; 
\item
in the case of surgeries $s3$ and $s4$ these two points should be complex conjugate along the path, and in the case of surgeries $s1$ and $s2$ they should be real all the time. 
\end{itemize}

Let us try to lift the map $I$ to a map $\tilde I: [0,1] \to \Theta$  such that $\tilde I (0) = f$ and $I \equiv \Lambda \circ \tilde I$.
By the covering homotopy property of the submersion $\Lambda$, this map is unique (if exists). Let $T \in (0,1]$ be the upper bound of values $\tau$ for which this covering homotopy defines the map $\tilde I$. By the algebraicity of the path $I$, the limit 
\begin{equation}\lim_{\tau \to T-0}\tilde I(\tau)
\label{limlim}\end{equation}
 in the compactification ${\mathbb C}P^1 \times {\mathbb C}P^8$ of the space $\Theta$ is well defined. By  Proposition 2 of \cite{Jaw2}, this limit point belongs to the  space $\left({\mathbb C}^1 \setminus \{1, -1\}\right) \times {\mathbb C}P^8$: indeed, by the construction of surgeries $s1$--$s4$ the intersection index of the vanishing cycles corresponding to two meeting critical values is $0, 1$ or $-1$ and so the arising multisingularity is elliptic. Therefore,  all values $\tilde I(\tau)$ for $\tau <T $ sufficiently close to $T$ are polynomials, whose principal homogeneous parts have the normal form from Table \ref{t1}
with the module $A$ lying in some compact disk $\bar U \subset {\mathbb C}^1 \setminus \{-1,1\}$. Let $\bar B$ be the unit ball centered at the origin in ${\mathbb C}^8$. 
It follows from the connectedness of the Dynkin diagram (see \cite{Gab}) that for any $A \in \bar U$ the function $f_A \equiv f_{A,0}$ is the unique function $f_{A,\lambda}$ with the same principal part $f_A$ such that all its critical values are equal to 0. Therefore, there exists a number $\delta>0$ such that all functions $f_{A, \lambda}$, $A \in \bar U$, all of whose critical values lie in the $\delta$-neighborhood of zero, lie in the compact set $\bar U \times \bar B$. In particular, the point (\ref{limlim}) belongs to $\Theta$. So, the point $T\in (0,1]$ cannot be less than $1$, because otherwise $\tilde I(T)$ is a strictly Morse polynomial,  the map $\Lambda$ is a submersion at it, and the path $\tilde I$ can be continued through $T$. If $T=1$, then the limit point (\ref{limlim}) is a polynomial with eight different critical values: it has one critical point of type $A_2$ if our surgery is of type $s1$ or $s3$, or  two Morse critical points with the same real critical value in the case of surgery $s2$ or $s4$. By the construction of our maps $I$ and $\Lambda$, all the polynomials $\tilde I(t),$ $t \in [0,1)$ are real, so also this limit double critical value is real. The desired Morse perturbation of the corresponding real critical point of type $A_2$ (in  the case of surgery $s1$ or $s3$) or the shift of the critical values of two real (in the case of $s2$) or complex conjugate (in the case of $s4$) Morse critical points can be performed within the real part of the family (\ref{vers0}) because of the versality property. \smallskip

In the case of the class $X_9^1$, the normal form $x y(x^2 +2 A x y + y^2)$ of principal parts of our polynomials does not coincide with the canonical form $x^4 + 2A x^2y^2 + y^4$ used in Proposition 2 of \cite{Jaw2}, in particular the corresponding family (\ref{vers0}) does not coincide with the versal family used there. Therefore, we consider the path $\tilde I$ together with its image  in the parameter space of the latter versal family under the map of parameter spaces defined by the substitution 
\begin{equation}
x = \frac{1}{\Omega(A)}(\tilde x + \alpha(A) \tilde y), \ y= \frac{1}{\Omega(A)}(\tilde y +\alpha(A) \tilde x), 
\end{equation}
taking $x y(x^2 + 2A x y + y^2)$ to a polynomial of the form $\tilde x^4 + 2\tilde A(A) \tilde x^2 \tilde y^2 +\tilde y^4$
(no matter that it gives non-real polynomials); here $$\alpha(A) = -A+\sqrt{A^2-1}+ \sqrt{2 A^2 - 2 A \sqrt{A^2-1}-2},$$  $\Omega(A) = (\alpha(A) +2A \alpha^2(A) + \alpha^3(A))^{-1/4}.$
This map is regular for $A \neq \pm 1$ and preserves the $j$-invariant of the homogeneous principal parts, therefore the image of the point $\tilde I(T)$ from the previous argument again belongs to the space $\Theta$, and so the same is true for the point $\tilde I$. The rest of the proof is the same as before. \hfill $\Box$

\begin{proposition}
\label{refprop}
If  a virtual Morse function is associated with two generic polynomials $f$, $\tilde f$ of degree four, then either $f$ and $\tilde f$ belong to the same connected component of the space of generic polynomials of that degree, or their connected components are mapped to each other by the involution taking each polynomial $f(x,y) $ to $f(x, -y)$. \end{proposition}

Proofs of this proposition for polynomials with the principal parts of type $X_9^{\pm}, X_9^1$ and $X_9^{2}$ will be given separately in subsections \ref{proofref0}, \ref{proofref1} and \ref{proofref2}. 

\begin{definition} \rm
\label{her}
A generic polynomial $f$ of degree four is called {\em achiral} (respectively, {\em chiral}) if this involution takes its connected component of the set of Morse polynomials to itself (respectively, to a different component).
\end{definition}

\begin{definition} \rm
\label{svi}
The {\em set-valued  invariant} of a real Morse polynomial $f$ with a non-discriminantal principal part is the set of virtual Morse functions corresponding to all
vertices of the virtual component  of  $f$. 

The {\em  $\mbox{Card}$ invariant} of such a Morse polynomial $f$ is the cardinality of this set of vertices.
\end{definition}

The set-valued invariant determines the ``passport'' invariant (which can be read from the bottom line of any virtual Morse function representing it).
\smallskip

 It turns out (see Theorems \ref{mthmdX90}--\ref{mthmdX92} below) 
that in the case of degree four polynomials the passport and $\mbox{Card}$ invariants separate all virtual components $S(f)$. The values of $\mbox{Card}$ for all isotopy classes mentioned in Table \ref{TX+} are given in Tables \ref{TXpm}--\ref{TX2}. 
The numbers $\mbox{Card}$ for the virtual components associated with chiral polynomials are underlined in these tables. For example, the cell $(M, m_+) = (9,0)$ of Table \ref{TXpm} contains three numbers, one of which is underlined and therefore should be counted twice to get the number 4 of components with these values of $M$ and $m_+$, see the corresponding cell of Table \ref{TX+}.
\medskip

\begin{table}
\caption{Values of the Card invariant for isotopy classes of type  $X_9^+$ with different passports}
\label{TXpm}
\begin{tabular}{|c| c | c | c | c | c |}
\hline
$m_+ \backslash M$ & 1 & 3 & 5 & \ \ \ 7 & \ \ 9 \\
\hline
0 & 1136 & 768 & 1584 & 2528 + 2912 & \underline{7320} + 2460 + 6220 \\
1 & 0 & 256 & 192 & 384 & 1360 \\
\hline
$\Sigma$ & 1136 & 1024 & 1776 & 5824 & 17360 \\
\hline
\end{tabular}
\end{table}

\begin{table}
\caption{The Card invariants for isotopy classes of type  $X_9^1$ \qquad \qquad \qquad \qquad}
\label{TX1}
\begin{tabular}{|c| c | c | c | c | c |}
\hline
$m_+ \backslash M$ & 1 & 3 & 5 & 7 & 9 \\
\hline
0 & 1396 & 1608 & 3420 & 11920 & \underline{5220} + \underline{7400} + \underline{13000} + \underline{9600} \\
1 & 0 & 1608 & 1356 & 2384 & 9060 \\
2 & 0 & 0 & 3420 & 2384 & 2660 \\
3 & 0 & 0 & 0 & 11920 & 9060 \\
4 & 0 & 0 & 0 & 0 & \underline{5220} + \underline{7400} + \underline{13000} + \underline{9600} \\
\hline 
$\Sigma$ & 1396 & 3216 & 8196 & 28608 & 91220 \\
\hline
\end{tabular}
\end{table}

\begin{table}
\caption{The Card invariants for isotopy classes of type  $X_9^2$ \qquad \qquad \qquad \qquad}
\label{TX2}
\begin{tabular}{|c|| c |c | c | c |}
\hline
$m_+ \backslash M$ & 3 & 5 & 7 & 9 \\
\hline
0 & 1320 & 2988 & 9872 & 2880 + \underline{11360} + 16180 \\
1 & 0 & 2988 & 5648 & 20260 \\
2 & 0 & 0 & 9872 & 20260 \\
3 & 0 & 0 & 0 & 2880 + \underline{11360} + 16180 \\
\hline 
$\Sigma$ & 1320 & 5976 & 25392 & 101360 \\
\hline
\end{tabular}
\end{table}

Finally, our system of invariants of  Morse polynomials consists of

a) the topological type of the principal homogeneous part,

b) the set-valued invariant of Definition \ref{svi}, and 

c) the reflection class in the case of chiral polynomials.

By Propositions \ref{propmain} and \ref{refprop}, this system of invariants separates all isotopy classes of Morse polynomials of degree four.

\subsection{D-graph invariant}
\label{dinv}

In this subsection we consider only the real polynomials $({\mathbb C}^2,{\mathbb R}^2) \to ({\mathbb C}, {\mathbb R}),$ all of whose critical points are real. In this case  the set-valued invariant of \S \ref{svinv} has the following transparent interpretation.

Let $f$ be a generic polynomial ${\mathbb R}^2 \to {\mathbb R}$ of degree $d$ with only real critical points, in particular all its $(d-1)^2$ critical values are real and different, and   $0$ is a non-critical  value of $f$. The matrix of intersection indices of vanishing cycles $\Delta_i \in H_2(V_f)$ (numbered in the ascending order of the corresponding critical values)
 can be depicted by its  {\em Coxeter-Dynkin graph}  (see e.g. \cite{AVGZ}) with $(d-1)^2$ numbered vertices corresponding to all critical values of $f$. Namely, if  the intersection index $\langle \Delta_i, \Delta_j \rangle$ is positive, then the corresponding vertices  $v_i$ and $v_j$ are connected by  $ \langle \Delta_i, \Delta_j \rangle$  solid segments; if $\langle \Delta_i, \Delta_j \rangle$ is negative then they are connected  by $-\langle \Delta_i, \Delta_j \rangle$ dashed segments. 

\begin{definition} \rm
\label{dfdinv}The {\it D-graph}  of a generic real Morse polynomial $f$ with only real critical points is (the isomorphism class of) the 
oriented graph with vertices labeled with indices 0, 1 and 2, which is obtained from the Coxeter-Dynkin graph of $f$ by

1) orienting each edge of this graph from the vertex corresponding to the critical point with the lower critical value to that with the higher critical value;

2) the indication of the Morse index of each critical point of $f$ at the corresponding vertex of the graph, and

3) forgetting the numbering of the vertices.
\end{definition} 

\begin{remark} \rm
The notion of a $D$-graph is very close to that of an $R$-diagram, which is used in \cite{chis} 
for classifying decompositions of simple function singularities. Two differences are that $D$-graphs are oriented and show the integer values of the intersection indices of vanishing cycles, and not only the information about their non-triviality. 
For the simple singularities considered in \cite{chis}, the orientation of the graph, i.e.  the comparison of the critical values corresponding to the endpoints of each edge, follows directly from the Morse indices of their critical points, while for general Morse functions the {\em tunnel} edges (see \S \ref{norma} below) with counterintuitive orientation provide the main part of the intrigue.
\end{remark}

\noindent{\bf Notation.}
In the pictures of the D-graphs (see Figs.~\ref{miss90a}--\ref{X902}), instead of numerical indices, 
we will mark the vertices  corresponding to minima, saddlepoints and maxima by white circles, black circles and white squares, respectively.

\begin{remark} \rm
The D-graph  of a real Morse polynomial is determined by an arbitrary virtual Morse function associated with this polynomial. Indeed, the intersection matrix and the Morse indices are the elements of the virtual Morse function, and the orientation of the edges follows from the order of the rows and columns of this matrix determined by the order of the corresponding critical values.  In this way also D-graphs of arbitrary (abstract) virtual Morse functions with only real critical points are well defined.
\end{remark}

\begin{theorem}
In the restriction to generic Morse polynomials with only real critical points, the D-graphs are invariants of isotopy classes of Morse functions, and this
 invariant is equivalent to the set-valued invariant from \S \ref{svinv}.
\end{theorem}

\noindent 
{\it Proof} of this theorem consists of the following two lemmas.

\begin{lemma}
\label{p11}
Virtual surgeries that model surgeries of Morse polynomials with only real critical points, without changing the Morse isotopy class, do not change the corresponding D-graphs.
In particular, the D-graph of a virtual Morse function is determined by its set-valued invariant.
\end{lemma}

\begin{figure}
\begin{picture}(40,15)
\put(0,0){\line(1,0){40}}
\bezier{150}(5,0)(10,10)(15,0)
\bezier{100}(20,0)(18,3))(15,0)
\bezier{250}(15,0)(25,15)(35,0)
\put(5,0){\circle*{1}}
\put(20,0){\circle*{1}}
\put(35,0){\circle*{1}}
\put(15,0){\circle{1}}
\end{picture} $\Longrightarrow$
\begin{picture}(40,15)
\put(0,0){\line(1,0){40}}
\bezier{200}(5,0)(15,10)(20,5)
\put(20,0){\line(0,1){5}}
\bezier{200}(20,5)(25,10)(35,0)
\put(5,0){\circle*{1}}
\put(20,0){\circle*{1}}
\put(35,0){\circle*{1}}
\put(20,5){\circle{1}}
\end{picture} $\Longrightarrow$
\begin{picture}(40,15)
\put(0,0){\line(1,0){40}}
\bezier{150}(5,0)(15,15)(25,0)
\bezier{100}(20,0)(22,3))(25,0)
\bezier{250}(25,0)(30,10)(35,0)
\put(5,0){\circle*{1}}
\put(20,0){\circle*{1}}
\put(35,0){\circle*{1}}
\put(25,0){\circle{1}}
\end{picture}
\caption{Surgery $s6$}
\label{ss6}
\end{figure}
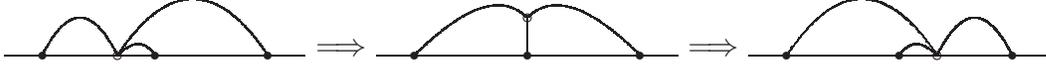

\noindent
{\it Proof}. These are only the surgeries $s2$, $s5$ and $s6$ from Definition \ref{elsur}. Surgeries of type $s2$ preserve all the ingredients of D-graphs
because they cannot permute the $i$th and $(i+1)$th critical points when $\langle \Delta_i, \Delta_{i+1} \rangle  \neq 0$. The surgery s5 or s6 can be realized as the composition of a) moving the non-critical value $0$ along an arc over a neighboring critical value (during which all paths to this value and the corresponding vanishing cycles in the fibers over these values are deformed continuously, and the  intersection indices do not change), see Fig.~\ref{ss6},  and b) adding a real constant  to $f$ moving this non-critical value to its original position $0$. These operations also preserve D-graphs.
 \hfill $\Box$
\medskip

The orientations of the edges of the D-graph of any abstract virtual Morse function with only real critical points define a partial order on the set of its vertices. Indeed, in the case of virtual Morse functions associated with real Morse functions these orientations  obviously cannot form oriented cycles, and this property is preserved by the elementary virtual surgeries.  

\begin{lemma}
\label{le13}
The set of all abstract virtual Morse functions defining the same D-graph is in a natural one-to-one correspondence with the set of pairs, consisting of 

a$)$ an extension of the partial order of vertices defined by this graph to a total order $($considered up to isomorphisms of graphs with numbered vertices, i.e. the extensions \
\unitlength 1mm
\begin{picture}(17,8)
\put(0.4,0){\circle{1}}
\put(15.6,0){\circle{1}}
\put(8,4){\circle*{1}}
\put(1,0.5){\vector(2,1){7}}
\put(15,0.5){\vector(-2,1){7}}
\put(0,2){\scriptsize 1}
\put(15,2){\scriptsize 2}
\put(7,1){\scriptsize 3}
\end{picture}
and \
\begin{picture}(17,8)
\put(0.4,0){\circle{1}}
\put(15.6,0){\circle{1}}
\put(8,4){\circle*{1}}
\put(1,0.5){\vector(2,1){7}}
\put(15,0.5){\vector(-2,1){7}}
\put(0,2){\scriptsize 2}
\put(15,2){\scriptsize 1}
\put(7,1){\scriptsize 3}
\end{picture}
of the D-graph \
\begin{picture}(17,8)
\put(0.4,0){\circle{1}}
\put(15.6,0){\circle{1}}
\put(8,4){\circle*{1}}
\put(1,0.5){\vector(2,1){7}}
\put(15,0.5){\vector(-2,1){7}}
\end{picture}
are the same$)$, and 

b$)$ an additional choice of an integer number in the range from $0$ to $(d-1)^2$ $($describing the last element of the virtual Morse function, i.e. the number of negative critical values$)$.

All abstract virtual Morse functions from such a set can be obtained from each other by elementary surgeries $s2$, $s5$ and $s6$.
\end{lemma}

\noindent
{\it Proof.} In the case of functions with only real critical points, the second element of a virtual Morse function (i.e., the string of intersection indices of vanishing cycles with the set of real points)  can be reconstructed from the Morse indices of the corresponding critical points and the intersection indices of these vanishing cycles with each other, see Theorem 3.1 in \S V.3  of \cite{APLT}. Therefore, this element does not imply any restrictions, and any virtual Morse function defining a D-graph is determined by this D-graph, a total order of its vertices and the number of negative critical values.
The latter number in our case does not participate in any way in the calculations and can be chosen arbitrarily.

Let us somehow number the vertices of the D-graph. The arithmetic space ${\mathbb R}^{(d-1)^2}$ is divided by hyperplanes into $(d-1)^2!$ open {\em Weyl chambers}, each of which consists of strings $(x_1, \dots, x_{(d-1)^2})$ with a fixed order of values of these numbers $x_i$ (which are all pairwise different in these chambers). 
Any total order of vertices of a D-graph with $(d-1)^2$ vertices defines one of these  Weyl chambers, namely the chamber consisting of vectors $(x_1, \dots, x_{(d-1)^2})$ such that $x_i< x_j$ every time the $i$-th vertex is subordinate to the $j$-th vertex in the sense of this order. The orientation of any edge $(i, j)$ of the D-graph  defines the half-space $\{x_i<x_j\}$ in this space, and the chambers corresponding to the extensions of its partial order defined by these edges to a total order
are exactly the ones lying in the intersection of all these half-spaces. This intersection is convex, so some points of any two chambers in it can be connected by a segment crossing only the walls defined by equations of the form $x_i = x_j$, such that the $i$-th and $j$-th vertices of the D-graph are not connected by its edges. All such passages correspond to 
 virtual Morse surgeries of type $s2$. Therefore, any total order of vertices extending the partial order encoded in the D-graph can be obtained from any other one by a sequence of such surgeries, in particular it also corresponds to an abstract virtual Morse function from the same set $S(f)$. \hfill $\Box$ $\Box$

\subsection{Main theorems}
\label{mthms}

\begin{figure}
\unitlength=0.6mm
\begin{picture}(120,40)
\put(0,5){\circle{2}}
\put(30,5){\circle{2}}
\put(60,5){\circle{2}}
\put(90,5){\circle{2}} 
\put(120,5){\circle{2}}
\put(89,5){\vector(-1,0){28}}
\put(0,6){\vector(0,1){28}} 
\put(29.3,5.7){\vector(-1,1){28.2}} 
\put(119.4,5.9){\vector(-1,2){13.7}} 
\put(90.5,6){\line(1,2){3}}
\put(95,15){\line(1,2){4}}
\put(100,25){\vector(1,2){4}}
\put(31,5){\vector(1,0){28}} 
\put(30.7,5.7){\vector(1,1){28.3}} 
\put(59.25,5.75){\line(-1,1){13.8}}
\put(44.5,20.5){\vector(-1,1){13.8}}
\put(89,6){\vector(-1,1){28}} 
\put(59.5,35){\line(0,-1){7}} 
\put(60.5,35){\line(0,-1){7}} 
\put(59.5,25){\line(0,-1){7}}
\put(60.5,25){\line(0,-1){7}} 
\put(59.5,15){\vector(0,-1){7}} 
\put(60.5,15){\vector(0,-1){7}}
\put(60,35){\vector(1,0){43.8}}
\put(59,35){\vector(-1,0){27.8}}
\put(60.8,5.6){\line(3,2){16.5}}
\put(78.5,17.2){\vector(3,2){26}}
\put(0,35){\circle*{2}}
\put(60,35){\circle*{2}}
\put(30,35){\circle*{2}}
\put(105,35){\circle*{2}}
\end{picture}
\caption{ D-graph for $X_9^+$ (no local maxima, case A)}
\label{miss90a}
\end{figure}
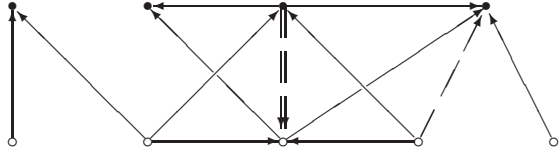

\begin{figure}
\unitlength=0.6mm
\begin{picture}(120,55)
\put(0,5){\circle{2}}
\put(30,5){\circle{2}}
\put(60,5){\circle{2}}
\put(90,5){\circle{2}} 
\put(120,5){\circle{2}}
\put(15,35){\circle*{2}}
\put(60,35){\circle*{2}}
\put(75,50){\circle*{2}}
\put(105,35){\circle*{2}}
\put(59.1,5.6){\line(-3,2){16.4}}
\put(60.9,5.6){\line(3,2){16.5}}
\put(41.5,17.2){\vector(-3,2){25.5}}
\put(78.5,17.2){\vector(3,2){25.5}}
\put(89,5){\vector(-1,0){28}}
\put(0.4,5.8){\vector(1,2){14}} 
\put(29.4,6.2){\line(-1,2){3}} 
\put(25,15){\line(-1,2){4}} 
\put(20,25){\vector(-1,2){4}} 
\put(119.6,5.8){\vector(-1,2){13.8}} 
\put(90.4,5.8){\line(1,2){3.5}}
\put(95,15){\line(1,2){4}}
\put(100,25){\vector(1,2){4}}
\put(31,5){\vector(1,0){28}} 
\put(30.7,5.7){\vector(1,1){28}} 
\put(60.3,5.9){\line(1,3){6.9}}
\bezier{100}(67.8,28.4)(68,29)(69.7,34.1)
\put(70.2,35.6){\vector(1,3){4.4}}
\put(89.3,5.7){\vector(-1,1){28.2}} 
\put(59.5,35){\line(0,-1){7}} 
\put(60.5,35){\line(0,-1){7}} 
\put(59.5,25){\line(0,-1){7}}
\put(60.5,25){\line(0,-1){7}} 
\put(59.5,15){\vector(0,-1){7}} 
\put(60.5,15){\vector(0,-1){7}}
\put(59,35){\vector(-1,0){43}}
\put(61,35){\vector(1,0){43}}
\put(60,35){\vector(1,1){14.3}}
\end{picture}
\caption{D-graph $X_9^+$ (no local maxima, case B)}
\label{miss90b}
\end{figure}
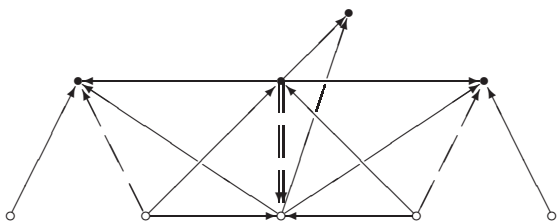

\begin{figure}
\unitlength=0.6mm
\begin{picture}(120,55)
\put(0,5){\circle{2}}
\put(30,5){\circle{2}}
\put(60,5){\circle{2}}
\put(90,5){\circle{2}} 
\put(120,5){\circle{2}}
\put(89,5){\vector(-1,0){28}}
\put(0.4,5.8){\vector(1,2){14}} 
\put(29.6,5.8){\vector(-1,2){14}} 
\put(119.6,5.8){\vector(-1,2){14}} 
\put(90.4,5.8){\vector(1,2){14}}
\put(31,5){\vector(1,0){28}} 
\put(30.7,5.7){\vector(1,1){28.7}} 
\put(60,35){\vector(1,1){14.5}} 
\put(89.3,5.7){\vector(-1,1){28.6}} 
\put(59.5,35){\line(0,-1){7}} 
\put(60.5,35){\line(0,-1){7}} 
\put(59.5,25){\line(0,-1){7}}
\put(60.5,25){\line(0,-1){7}} 
\put(59.5,15){\vector(0,-1){7}} 
\put(60.5,15){\vector(0,-1){7}}
\put(60.3,5.9){\line(1,3){6.9}}
\put(68,29){\vector(1,3){6.5}}
\put(15,35){\circle*{2}}
\put(60,35){\circle*{2}}
\put(75,50){\circle*{2}}
\put(105,35){\circle*{2}}
\end{picture}
\caption{D-graph $X_9^+$ (no local maxima, case C)}
\label{miss90c}
\end{figure}
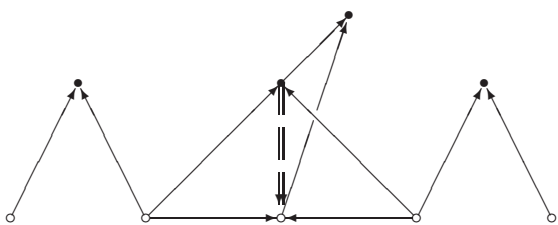

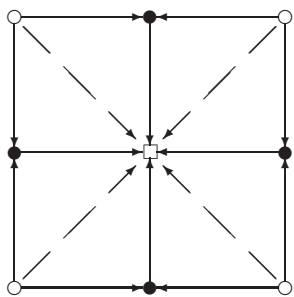
\begin{figure}
\unitlength=0.6mm
\begin{picture}(70,70)
\put(0,5){\circle{3}}
\put(30,5){\circle*{3}}
\put(60,5){\circle{3}}
\put(0,35){\circle*{3}}
\put(30,35){\makebox(0,0)[cc]{\tiny $\square$}}
\put(60,35){\circle*{3}}
\put(0,65){\circle{3}}
\put(30,65){\circle*{3}}
\put(60,65){\circle{3}}
\put(1.5,5){\vector(1,0){27}}
\put(58.5,5){\vector(-1,0){27}}
\put(1,35){\vector(1,0){27.5}}
\put(58.5,35){\vector(-1,0){27}}
\put(1.5,65){\vector(1,0){27}}
\put(58.5,65){\vector(-1,0){27}}
\put(0,6.5){\vector(0,1){27}}
\put(0,63.5){\vector(0,-1){27}}
\put(30,6.5){\vector(0,1){27}}
\put(30,63.5){\vector(0,-1){27}}
\put(60,6.5){\vector(0,1){27}}
\put(60,63.5){\vector(0,-1){27}}
\put(21,26){\vector(1,1){6}}
\put(18,23){\line(-1,-1){7}}
\put(8,13){\line(-1,-1){6}}
\put(39,26){\vector(-1,1){6}}
\put(42,23){\line(1,-1){7}}
\put(52,13){\line(1,-1){6}}
\put(21,44){\vector(1,-1){6}}
\put(18,47){\line(-1,1){7}}
\put(8,57){\line(-1,1){6}}
\put(39,44){\vector(-1,-1){6}}
\put(42,47){\line(1,1){7}}
\put(52,57){\line(1,1){6}}
\end{picture}
\caption{D-graph for $X_9^+$, one local maximum}
\label{X901mx}
\end{figure}

\begin{theorem}[$X_9^+$]
\label{mthmdX90}
The ``passport''  invariant $(M, m_+)$ $($see \S \ref{itriv}$)$ of a
Morse polynomial ${\mathbb R}^2 \to {\mathbb R}$ of degree four with a positive non-discriminantal principal homogeneous part can be equal to only one of the pairs 
\begin{equation}
\label{z1}
(1, 0), \ (3, 0), \ (5, 0), \ (7, 0), \ (9, 0), \ (3, 1), \ (5, 1), \ (7, 1), \ (9, 1),
\end{equation}
 in particular such a polynomial cannot have more than one local maximum. 

 The  D-graph invariant takes three different values on polynomials with passport $(M, m_+) = (9, 0) $: they are shown in Figs.~\ref{miss90a}, \ref{miss90b} and \ref{miss90c}.  The values of the $\mbox{Card}$ invariant of the polynomials with these D-graphs are 7320, 2460, and 6220, respectively. The set of Morse polynomials with the D-graph shown in Fig.~\ref{miss90a} consists of two different connected components that can be moved to each other by a reflection in any line in ${\mathbb R}^2$. Each of the sets of Morse polynomials with the D-graphs from Figs.~\ref{miss90b} and \ref{miss90c} consists of a single connected component.

There is exactly one connected component of the space of Morse polynomials with $(M,m_+)=(9,1) $; its D-graph is shown in Fig.~\ref{X901mx}, and  its Card invariant is  1360.

The polynomials with $(M, m_+) = (7, 0)$ belong to exactly two connected components,  the $\mbox{Card}$ invariant takes on them values    2528 and 2912.

For each remaining value of the passport invariant $(M, m_+)$ from the list $($\ref{z1}$)$, all Morse functions with this value form a single non-empty isotopy class. The values of the  $\mbox{Card}$ invariant  for them are given in Table \ref{TXpm} on page \pageref{TXpm}.
\end{theorem}

\unitlength=0.8mm
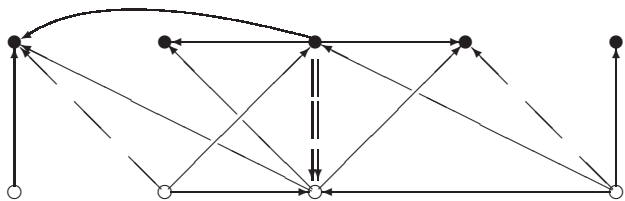
\begin{figure}
\begin{picture}(100,42)
\put(0,5){\circle{2}}
\put(0,30){\circle*{2}}
\put(25,5){\circle{2}}
\put(50,30){\circle*{2}}
\put(50,5){\circle{2}}
\put(75,30){\circle*{2}}
\put(100,5){\circle{2}}
\put(100,30){\circle*{2}}
\put(25,30){\circle*{2}}
\put(0,6){\vector(0,1){23}}
\put(7,23){\vector(-1,1){6}}
\put(16, 14){\line(-1,1){6}}
\put(24.4,5.6){\line(-1,1){5.4}}
\put(49.4,5.3){\line(-2,1){15.5}}
\put(32.6,13.7){\vector(-2,1){31.6}}
\put(26,5){\vector(1,0){23}}
\put(25.6,5.6){\vector(1,1){23.5}}
\put(50.6,5.6){\vector(1,1){23.7}}
\put(50,30){\vector(1,0){24.2}}
\put(99.4,5.3){\line(-2,1){32.1}}
\put(66,22){\vector(-2,1){15}}
\put(82,23){\vector(-1,1){6}}
\put(91,14){\line(-1,1){6}}
\put(99.5,5.5){\line(-1,1){5.5}}
\put(100,6){\vector(0,1){23}}
\put(49.5,5.5){\line(-1,1){11.5}}
\put(37,18){\vector(-1,1){11.5}}
\put(49.5,12){\vector(0,-1){5}}
\put(50.5,12){\vector(0,-1){5}}
\put(49.5,14){\line(0,1){6}}
\put(50.5,14){\line(0,1){6}}
\put(49.5,21){\line(0,1){6}}
\put(50.5,21){\line(0,1){6}}
\put(49,30){\vector(-1,0){23}}
\bezier{500}(50,30.5)(15,40)(2,31)
\put(2,31){\vector(-2,-1){1}}
\put(99,5){\vector(-1,0){48}}
\end{picture}
\caption{ D-graph for $X_9^1$ (no local maxima, case A)}
\label{DGX921A}
\end{figure}

\unitlength=0.8mm
\begin{figure}
\begin{picture}(100,35)
\put(0,5){\circle{2}}
\put(0,30){\circle*{2}}
\put(25,5){\circle{2}}
\put(50,30){\circle*{2}}
\put(50,5){\circle{2}}
\put(75,30){\circle*{2}}
\put(100,5){\circle{2}}
\put(100,30){\circle*{2}}
\put(25,30){\circle*{2}}
\put(0,6){\vector(0,1){23}}
\put(24.4,5.6){\vector(-1,1){23.7}}
\put(26,5){\vector(1,0){23}}
\put(25.6,5.6){\vector(1,1){23.7}}
\put(50.6,5.6){\vector(1,1){23.7}}
\put(50,30){\vector(1,0){24.2}}
\put(99.4,5.6){\line(-1,1){5.4}}
\put(91,14){\line(-1,1){6}}
\put(82,23){\vector(-1,1){6}}
\put(100,6){\vector(0,1){23}}
\put(50,30){\vector(-1,0){24}}
\put(99.2,5.4){\line(-2,1){31.8}}
\put(66,22){\vector(-2,1){15.5}}
\put(49.4,5.6){\line(-1,1){11.3}}
\put(37,18){\vector(-1,1){11}}
\put(49.5,12){\vector(0,-1){5}}
\put(50.5,12){\vector(0,-1){5}}
\put(49.5,14){\line(0,1){7}}
\put(50.5,14){\line(0,1){7}}
\put(49.5,23){\line(0,1){7}}
\put(50.5,23){\line(0,1){7}}
\put(99,5){\vector(-1,0){48}}
\end{picture}
\caption{ D-graph for $X_9^1$ (no local maxima, case B)} \label{DGX921B} 
\end{figure}
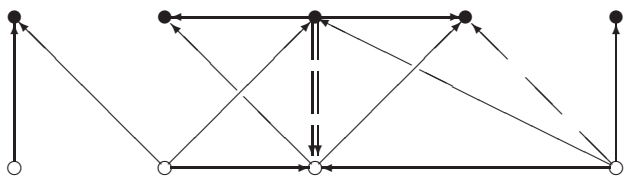

\unitlength=0.8mm
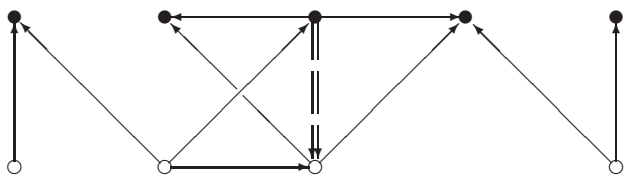
\begin{figure}
\begin{picture}(100,35)
\put(0,5){\circle{2}}
\put(0,30){\circle*{2}}
\put(25,5){\circle{2}}
\put(50,30){\circle*{2}}
\put(50,5){\circle{2}}
\put(75,30){\circle*{2}}
\put(100,5){\circle{2}}
\put(100,30){\circle*{2}}
\put(25,30){\circle*{2}}
\put(0,6){\vector(0,1){23.1}}
\put(24.3,5.7){\vector(-1,1){23.3}}
\put(25.7,5.7){\vector(1,1){23.2}}
\put(26,5){\vector(1,0){23}}
\put(50.7,5.7){\vector(1,1){23.2}}
\put(50,30){\vector(1,0){24.2}}
\put(99.3,5.7){\vector(-1,1){23}}
\put(100,6){\vector(0,1){23}}
\put(50,30){\vector(-1,0){24}}
\put(49.3,5.7){\line(-1,1){11.3}}
\put(37,18){\vector(-1,1){11}}
\put(49.5,12){\vector(0,-1){5.5}}
\put(50.5,12){\vector(0,-1){5.5}}
\put(49.5,14){\line(0,1){7}}
\put(50.5,14){\line(0,1){7}}
\put(49.5,23){\line(0,1){7}}
\put(50.5,23){\line(0,1){7}}
\end{picture}
\caption{ D-graph $X_9^1$ (no local maxima, case C)}
\label{DGX921C}
\end{figure}

\unitlength 0.8 mm
\begin{figure}
\begin{picture}(100,40)
\put(0,5){\circle{2}}
\put(0,30){\circle*{2}}
\put(25,5){\circle{2}}
\put(50,30){\circle*{2}}
\put(50,5){\circle{2}}
\put(75,30){\circle*{2}}
\put(100,5){\circle{2}}
\put(100,30){\circle*{2}}
\put(25,30){\circle*{2}}
\put(0,6){\vector(0,1){23}}
\put(7,23){\vector(-1,1){6}}
\put(16, 14){\line(-1,1){6}}
\put(24.3,5.7){\line(-1,1){5.3}}
\put(49.4,5.3){\line(-2,1){15.4}}
\put(32.6,13.7){\vector(-2,1){31.5}}
\put(26,5){\vector(1,0){23}}
\put(25.7,5.7){\vector(1,1){23.3}}
\put(50.7,5.7){\vector(1,1){23.3}}
\put(50,30){\vector(1,0){24.2}}
\put(99.3,5.7){\vector(-1,1){23.3}}
\put(100,6){\vector(0,1){23}}
\put(50,30){\vector(-1,0){23.5}}
\put(49.3,5.7){\line(-1,1){11}}
\put(37,18){\vector(-1,1){11}}
\put(49.5,12){\vector(0,-1){5}}
\put(50.5,12){\vector(0,-1){5}}
\put(49.5,14){\line(0,1){7}}
\put(50.5,14){\line(0,1){7}}
\put(49.5,23){\line(0,1){7}}
\put(50.5,23){\line(0,1){7}}
\bezier{300}(50,30)(15,40)(2,31)
\put(2,31){\vector(-2,-1){1}}
\end{picture}
\caption{ D-graph $X_9^1$ (no local maxima, case D)}
\label{DGX921D}
\end{figure}
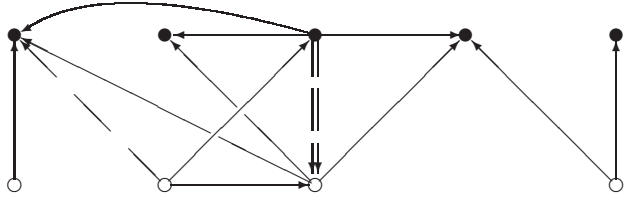

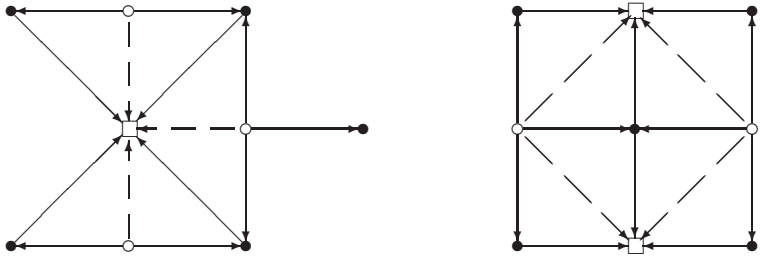
\begin{figure}
\unitlength=0.52mm
\begin{picture}(95,70)
\put(0,65){\circle*{3}}
\put(90,35){\circle*{3}}
\put(60,65){\circle*{3}}
\put(0,5){\circle*{3}}
\put(30,35){\makebox(0,0)[cc]{\tiny $\square$}}
\put(60,5){\circle*{3}} 
\put(30,65){\circle{3}}
\put(30,5){\circle{3}}
\put(60,35){\circle{3}}
\put(28.5,65){\vector(-1,0){27.5}} 
\put(31.5,65){\vector(1,0){27.5}}
\put(60,33.5){\vector(0,-1){27.5}}
\put(60,36.5){\vector(0,1){27.5}}
\put(28.5,5){\vector(-1,0){27.5}}
\put(31.5,5){\vector(1,0){27.5}}
\put(61.5,35){\vector(1,0){27.5}}
\put(37,35){\vector(-1,0){5}}
\put(30,62){\line(0,-1){6}}
\put(30,52){\line(0,-1){6}}
\put(30,42){\vector(0,-1){5}}
\put(41,35){\line(1,0){6}}
\put(51,35){\line(1,0){6}}
\put(1,6){\vector(1,1){27.3}}
\put(1,64){\vector(1,-1){27.3}}
\put(59,6){\vector(-1,1){27.3}}
\put(59,64){\vector(-1,-1){27.3}}
\put(30,27){\vector(0,1){5}}
\put(30,23){\line(0,-1){6}}
\put(30,13){\line(0,-1){6}}
\end{picture} \qquad \qquad \quad
\begin{picture}(70,70)
\put(0,5){\circle*{3}}
\put(30,5){\makebox(0,0)[cc]{\tiny $\square$}}
\put(60,5){\circle*{3}}
\put(0,35){\circle{3}}
\put(30,35){\circle*{3}}
\put(60,35){\circle{3}}
\put(0,65){\circle*{3}}
\put(30,65){\makebox(0,0)[cc]{\tiny $\square$}}
\put(60,65){\circle*{3}}
\put(0,33.5){\vector(0,-1){27.5}}
\put(0,36.5){\vector(0,1){27.5}}
\put(60,33.5){\vector(0,-1){27.5}}
\put(60,36.5){\vector(0,1){27.5}}
\put(30,34){\vector(0,-1){27}}
\put(30,36){\vector(0,1){27}}
\put(1,5){\vector(1,0){27}}
\put(59,5){\vector(-1,0){27}}
\put(1.5,35){\vector(1,0){27.5}}
\put(58.5,35){\vector(-1,0){27.5}}
\put(1,65){\vector(1,0){27}}
\put(59,65){\vector(-1,0){27}}
\put(2,37){\line(1,1){7}}
\put(12,47){\line(1,1){7}}
\put(21,56){\vector(1,1){7}}
\put(2,33){\line(1,-1){7}}
\put(12,23){\line(1,-1){7}}
\put(21,14){\vector(1,-1){7}}
\put(58,37){\line(-1,1){7}}
\put(48,47){\line(-1,1){7}}
\put(39,56){\vector(-1,1){7}}
\put(58,33){\line(-1,-1){7}}
\put(48,23){\line(-1,-1){7}}
\put(39,14){\vector(-1,-1){7}}
\end{picture}
\caption{$X_9^1$, one maximum (left) and two maxima (right)}
\label{X9011}
\end{figure}

\begin{theorem}[$X_9^1$]
\label{mthmdX91}
Morse polynomials of degree four with non-discriminantal principal homogeneous part vanishing on exactly two real lines can have arbitrary values of the ``passport'' invariant $(M, m_+)$ such that $M$ is equal to 1, 3, 5, 7 or 9, and $0 \leq m_+ < M/2$. 

The D-graph invariant takes exactly four values on polynomials of this type with $(M, m_+) = (9, 0)$; these values are 
shown in Figs.~\ref{DGX921A}, \ref{DGX921B}, \ref{DGX921C} and  \ref{DGX921D}. The values of the $\mbox{Card}$  invariant on the corresponding polynomials are 5220, 7400, 13000, and  9600, respectively.
The set of Morse polynomials with either of these four values of the D-graph invariant contains exactly two connected components  which can be transformed into each other by a reflection in any line in ${\mathbb R}^2$.

All statements of the previous paragraph remain  true for the polynomials with $(M, m_+) = (9, 4)$ if we reverse all arrows in the D-graphs and replace minima by maxima. 

Morse polynomials with a principal part of type $X_9^1$, nine real critical points and exactly one $($respectively, two, respectively, three$)$ local maxima form a single connected component;
their D-graphs are shown in Fig.~\ref{X9011} left
$($respec\-ti\-vely, Fig.~\ref{X9011} right, respectively, the graph of Fig.~\ref{X9011} left with reversed arrows and minima replaced by maxima and vice versa$)$. Their Card invariants are 9060, 2660, and 9060, respectively.

For any other admissible value of the passport invariant $(M, m_+)$,  all polynomials with that value form a single  isotopy class. The values of the $\mbox{Card}$ invariant for them are given in Table \ref{TX1}.
\end{theorem}

\unitlength=0.8mm
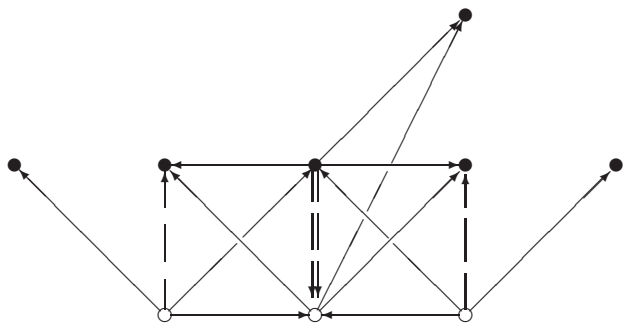
\begin{figure}
\begin{picture}(100,60)
\put(0,30){\circle*{2}} 
\put(25,30){\circle*{2}} 
\put(50,30){\circle*{2}} 
\put(75,55){\circle*{2}} 
\put(75,30){\circle*{2}} 
\put(100,30){\circle*{2}}
\put(25,5){\circle{2}}
 \put(50,5){\circle{2}} 
\put(75,5){\circle{2}}
\put(50,30){\vector(-1,0){24}}
\put(24.3,5.7){\vector(-1,1){23.5}} 
\put(25,6){\line(0,1){5}} 
\put(25,14){\line(0,1){6}} 
\put(25,23){\vector(0,1){6}}
\put(26,5){\vector(1,0){23}}
\put(25.7,5.7){\line(1,1){11.3}}
\put(38,18){\vector(1,1){11.5}}
\put(62.2,17.9){\vector(-1,1){11.5}}
\put(74.3,5.7){\line(-1,1){11.3}}
\put(74,5){\vector(-1,0){23}}
\put(49.3,5.7){\vector(-1,1){23.4}} 
\put(50.4,5.8){\line(1,2){7.6}} 
\put(58.7,22.4){\line(1,2){3.3}}
\put(62.7,30.4){\vector(1,2){11.8}}
\put(75.7,5.7){\vector(1,1){23.3}}
 \put(75,6){\line(0,1){6}} 
 \put(75,14){\line(0,1){7}} 
\put(75,23){\vector(0,1){5.5}} 
\put(50.7,5.7){\vector(1,1){23.3}}
 \put(50,30){\vector(1,1){24}} 
\put(50,30){\vector(1,0){24}} 
\put(49.5,14){\vector(0,-1){6}} 
\put(50.5,14){\vector(0,-1){6}} 
\put(49.5,22){\line(0,-1){6}} 
\put(50.5,22){\line(0,-1){6}} 
\put(49.5,30){\line(0,-1){6}} 
\put(50.5,30){\line(0,-1){6}} 
\end{picture}
\caption{ D-graph of type $X_9^2$ (no local maxima, case A)}
\label{DGX92A}
\end{figure}

\unitlength=0.8mm
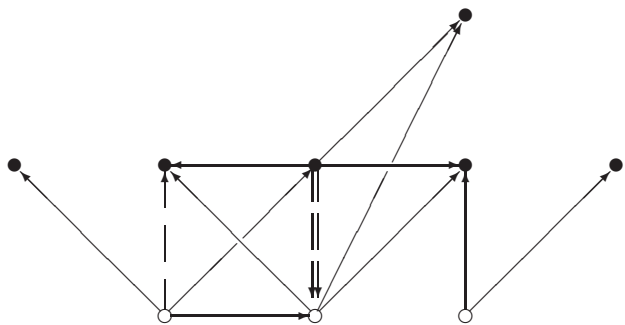
\begin{figure}
\begin{picture}(100,60)
\put(0,30){\circle*{2}} 
\put(25,30){\circle*{2}} 
\put(50,30){\circle*{2}} 
\put(75,55){\circle*{2}} 
\put(75,30){\circle*{2}} 
\put(100,30){\circle*{2}}
\put(25,5){\circle{2}}
 \put(50,5){\circle{2}} 
\put(75,5){\circle{2}}
\put(50,30){\vector(-1,0){24}}
\put(24.3,5.7){\vector(-1,1){23.3}} 
\put(25,6){\line(0,1){5}} 
\put(25,14){\line(0,1){6}} 
\put(25,23){\vector(0,1){6}}
\put(26,5){\vector(1,0){23}}
\put(25.7,5.7){\line(1,1){11.3}}
\put(38,18){\vector(1,1){11.5}}
\put(49.3,5.7){\vector(-1,1){23.3}} 
\put(50.4,5.8){\line(1,2){11.6}} 
\put(62.8,30.6){\vector(1,2){11.5}}
\put(75.7,5.7){\vector(1,1){23.3}}
 \put(75,6){\vector(0,1){23}} 
 \put(50.7,5.7){\vector(1,1){23.3}}
 \put(50,30){\vector(1,1){24}} 
\put(50,30){\vector(1,0){24}} 
\put(49.5,14){\vector(0,-1){6}} 
\put(50.5,14){\vector(0,-1){6}} 
\put(49.5,22){\line(0,-1){6}} 
\put(50.5,22){\line(0,-1){6}} 
\put(49.5,30){\line(0,-1){6}} 
 \put(50.5,30){\line(0,-1){6}} 
\end{picture}
\caption{ D-graph $X_9^2$ (no local maxima, case B)}
\label{DGX92B}
\end{figure} 

\unitlength 0.8mm
\begin{figure}
\begin{picture}(100,60)
\put(0,30){\circle*{2}} 
\put(25,30){\circle*{2}} 
\put(50,30){\circle*{2}} 
\put(75,55){\circle*{2}} 
\put(75,30){\circle*{2}} 
\put(100,30){\circle*{2}}
\put(25,5){\circle{2}}
 \put(50,5){\circle{2}} 
\put(75,5){\circle{2}}
\put(50,30){\vector(-1,0){24}}
\put(24.3,5.7){\vector(-1,1){23.3}} 
\put(25,6){\vector(0,1){23}}
\put(49.3,5.7){\vector(-1,1){23.3}} 
\put(50.4,5.8){\line(1,2){11.6}} 
\put(62.9,30.8){\vector(1,2){11.8}}
\put(75.7,5.7){\vector(1,1){23.3}}
 \put(75,6){\vector(0,1){23}} 
\put(50.7,5.7){\vector(1,1){23.3}}
 \put(50,30){\vector(1,1){24}} 
\put(50,30){\vector(1,0){24}} 
\put(49.5,14){\vector(0,-1){6}} 
\put(50.5,14){\vector(0,-1){6}} 
\put(49.5,22){\line(0,-1){6}} 
\put(50.5,22){\line(0,-1){6}} 
\put(49.5,30){\line(0,-1){6}} 
 \put(50.5,30){\line(0,-1){6}} 
\end{picture}
\caption{ D-graph $X_9^2$ (no local maxima, case C)}
\label{DGX92C}
\end{figure}
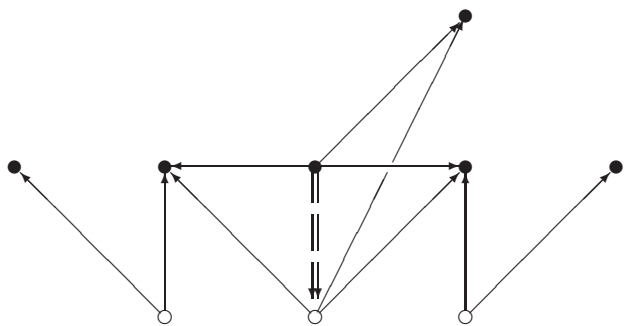

\unitlength 0.9mm
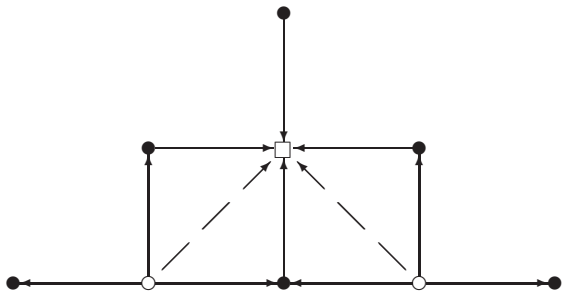
\begin{figure}
\begin{picture}(60,48)
\put(0,5){\circle*{2}}
\put(80,5){\circle*{2}}
\put(20,5){\circle{2}}
\put(40,5){\circle*{2}}
\put(60,5){\circle{2}}
\put(20,25){\circle*{2}}
\put(60,25){\circle*{2}}
\put(40,25){\makebox(0,0)[cc]{\tiny $\square$}}
\put(40,45){\circle*{2}}
\put(19,5){\vector(-1,0){18}}
\put(61,5){\vector(1,0){18}}
\put(20,6){\vector(0,1){18.2}}
\put(60,6){\vector(0,1){18.2}}
\put(21,5){\vector(1,0){18}}
\put(59,5){\vector(-1,0){18}}
\put(21,25){\vector(1,0){18}}
\put(59,25){\vector(-1,0){18}}
\put(40,6){\vector(0,1){18}}
\put(40,44){\vector(0,-1){18}}
\put(34,19){\vector(1,1){4.3}}
\put(32,17){\line(-1,-1){4}}
\put(26,11){\line(-1,-1){4}}
\put(46,19){\vector(-1,1){4.3}}
\put(48,17){\line(1,-1){4}}
\put(54,11){\line(1,-1){4}}
\end{picture}
\caption{$X_9^2$, one maximum}
\label{X902}
\end{figure}

\begin{theorem}[$X_9^2$]
\label{mthmdX92}
Morse polynomials of degree four with the principal homogeneous  parts vanishing on four real lines can have arbitrary values of the invariant $(M, m_+)$ such that $M$ is 3, 5, 7 or 9, and $0 \leq m_+ < M/2 - 1$. 

All three possible D-graphs of polynomials of this type with $(M, m_+) = (9, 0)$ are shown in Figs.~\ref{DGX92A}, \ref{DGX92B} and \ref{DGX92C}. The corresponding values of the $\mbox{Card}$ invariant are equal respectively to 2880, 11360, and 16180. The set of polynomials with the D-graph of Fig.~\ref{DGX92B} consists of two connected components which can be sent one to the other by a reflection in any line in ${\mathbb R}^2$.
The sets of polynomials with the D-graphs shown in Figs.~\ref{DGX92A} and \ref{DGX92C} each consist of a single connected component.

All statements of the previous paragraph remain true for the polynomials with $(M, m_+) = (9, 3)$ if we reverse the orientation of all arrows and replace local minima by maxima and vice versa. 

Polynomials with  a principal part of type $X_9^2$
and $(M, m_+)= (9,1)$ $($respec\-ti\-vely, $(9,2))$ form a single connected component, their D-graph is shown in Fig.~\ref{X902} $($respectively, is obtained from this one by reversing the orientations and replacing minima by maxima and vice versa$)$;  the value of the Card invariant for these two cases is  20260.

For any other admissible value of the passport invariant $(M, m_+)$, all Morse polynomials with this value form a single  isotopy class. The values of the $\mbox{Card}$ invariant for all of them are given in Table \ref{TX2}.
\end{theorem}

\subsection{On the proofs of main theorems}
\label{onproofs}
All possible values of the $\mbox{Card}$ and D-graph invariants mentioned in Theorems   \ref{mthmdX90} -- \ref{mthmdX92}
were found using the computer program described in \cite{AGLV2}, \cite{Vcau}, \cite{VS}.
 Namely, for any particular type $X_9^\ast$ we use the Gusein-Zade--A'Campo method \cite{AC}, \cite{GZ} to compute the intersection matrix of vanishing cycles of a Morse polynomial of this type having only real critical points. Using theorem 1.4 of \S 5.1 of \cite{APLT}, we then calculate the intersection indices of these vanishing cycles with the set of real points, and thus obtain a particular virtual Morse function associated with this polynomial.
Starting from this initial data,  the combinatorial program

{\scriptsize\begin{verbatim}https://drive.google.com/file/d/185T0euaR7Gs_AXqt6_otm3vQJmnMEM7g/view?usp=sharing
\end{verbatim}}
\noindent
runs through the entire formal graph of this type, counting in particular  
the number of all virtual Morse functions of this type, as well as the numbers of all these virtual functions with any value of the passport invariant $(M, m_+)$. For each passport value, for which this number is not zero, it then  (on special request) finds a virtual Morse function with this value. Starting from this virtual Morse function, a slightly modified version 

{\scriptsize \begin{verbatim}https://drive.google.com/file/d/1TqXsP1msOjkppCRczMJW1pHpzOKhFAvt/view?usp=sharing
\end{verbatim}}
\noindent
of the same program (with virtual surgeries of types $s1$ and $s3$ disabled)
counts the number of virtual functions in its virtual component, i.e. the Card invariant of this virtual function. If this number is smaller than the number of all virtual functions with this passport, then we find another virtual function with the same passport but not in the same virtual component, and calculate its $\mbox{Card}$ value again. We continue this search until the sum of the different values of the Card invariant of virtual functions with any given passport becomes equal to the total number of virtual functions with that passport. This computation proves all statements of Theorems \ref{mthmdX90}--\ref{mthmdX92} concerning the lists of possible virtual components and the values of the $\mbox{Card}$ invariant for them, and also gives us all values of the D-graph invariants of polynomials with only real critical points. 

By Proposition \ref{propmain} all these virtual components can be realized by some Morse polynomials. 
In the following sections  \ref{reaax0}--\ref{reaax2} we present such realizations and investigate their chirality. Based on their topological properties, we then prove Proposition \ref{refprop} for all types $X_9^*$ of degree four polynomials with non-discriminantal principal parts, in particular justify the numbers of real components of the space of Morse functions for each virtual component. \medskip

The sets of fourth degree polynomials whose homogeneous principal parts vanish on two or four real lines are invariant under the multiplication by $-1$. This multiplication takes the Morse polynomials with any value  $(M, m_+)$ of the ``passport'' invariant to polynomials with the value $(M, [M/2]-m_+)$ in the first case and $(M, [M/2]-1 -m_+) $ in the second.  Therefore, it is sufficient to study only such polynomials with $m_+ \leq [M/4]$ (respectively, $m_+ \leq [M/4]-1$).

\begin{remark} \rm
In the similar problems for other classes of polynomials and singularities, it can happen that several virtual components have the same values of the Card invariant; this makes it even more difficult to analyze the results of the program. This does not happen for the degree 4 polynomials.
\end{remark}

\begin{remark} \rm 
In the case of higher degree polynomials, some of the subgraphs $S(f)$ may be infinite (as was pointed out by V.I.~Arnold in a discussion of the prospects of our program, see e.g. problem 1984-9 in \cite{AP}). Fortunately,  this is not the case for $X_9$ and other parabolic singularities. If this is the case, we can define the restricted versions of the set-valued invariant, in whose computation only the virtual Morse functions described by numbers satisfying some upper bounds are counted. 
In any case, the direct computations by the previous method for $d \geq 5$ are very long, so that additional simplifying tricks will be needed to make progress in these computations.
\end{remark}

\subsection{On the classification of strictly Morse polynomials}
\label{stric}

In \cite{AN}, a topological classification of strictly Morse polynomials ${\mathbb R}^2 \to {\mathbb R}$ of degree four with negative definite principal part and nine real critical points was considered, and some estimates on the number of classes were proved. This problem is analogous to the one considered in the present paper for the case $X_9^-$ (which is completely isomorphic to $X_9^+$), but  concerns classes of {\em strictly} Morse functions (not allowing equal critical values at different critical points). In particular, its classes are generally not preserved by elementary surgeries of type $s2$.

\begin{figure}
\unitlength=0.9mm
\begin{picture}(45,54)
\put(0.7,0.7){\line(1,1){39}}
\put(10,10){\line(1,-1){4.3}}
\put(20,20){\line(1,-1){4.3}}
\put(30,30){\line(1,-1){4.3}}
\put(40,40){\line(1,-1){4.3}}
\put(40,40){\line(0,1){4}}
\put(40,48){\line(0,1){4}}
\put(0,0){\circle{2}}
\put(10,10){\circle*{2}}
\put(20,20){\circle*{2}}
\put(30,30){\circle*{2}}
\put(40,40){\circle*{2}}
\put(15,5){\circle{2}}
\put(25,15){\circle{2}}
\put(35,25){\circle{2}}
\put(45,35){\circle{2}}
\end{picture}
\caption{Simplest Reeb graph for polynomials of class $X_9^+$}
\label{reeb}
\end{figure}
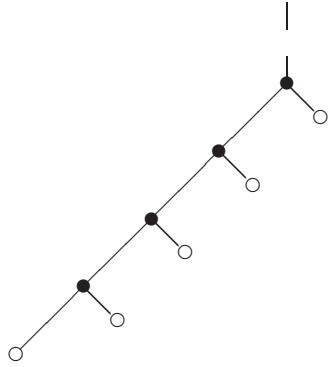

On the other hand, the {\em topological equivalence} considered there is weaker than the isotopy equivalence.   For example, all three classes of polynomials without local maxima, which have the D-graphs of Figs.~
\ref{miss90a}, \ref{miss90b} and \ref{miss90c} and are thus non-isotopic to each other, contain Morse polynomials, whose Kronrod--Reeb graphs are as shown in Fig.~\ref{reeb}, and hence these polynomials are equivalent from the point of view of \cite{AN}. In this subsection we consider the more rigorous isotopy classification of strictly Morse polynomials.

\begin{theorem}
\label{mthst}
For any chiral isotopy class 
$($respectively, any achiral isotopy class$)$ of degree four Morse polynomials with nine real critical points, the number of isotopy classes of {\em strictly} Morse polynomials within this class is equal to $\frac{1}{10}$  $($respectively, $\frac{1}{5})$ of the value of the Card invariant  of that class. 
\end{theorem}

\noindent
{\it Proof.} For any generic polynomial $f$, denote by $f^\uparrow$ an arbitrary polynomial of the form $f + c$, where the constant $c$ is such that all real critical values of $f+c$ are positive. 

\begin{lemma}
\label{lemel}
1. The generic polynomials $f$ and $\tilde f$ are isotopic in the class of strictly Morse polynomials if and only if $f^\uparrow$ and $\tilde f^\uparrow$ are isotopic in the class of strictly Morse polynomials, all of whose critical values at real critical points are positive.

2. If two generic  polynomials of degree $d$, whose critical points are all real and  whose critical values are all positive,  are isotopic in the class of strictly Morse polynomials, then the virtual Morse functions associated with them are the same.

3. A generic polynomial $f(x,y)$ of degree four with only real critical points cannot be isotopic to its mirror image $f(x,-y)$ in the class of strictly Morse polynomials.
\end{lemma}

\noindent
{\it Proof.} 1. Every strictly Morse polynomial $f$ is isotopic to all polynomials $f+c$ (in particular, to $f^\uparrow$)  in the class of strictly Morse polynomials. If $f$ is isotopic to a polynomial $\tilde f$, let $c_-$ be the lowest critical value at real critical points of all intermediate polynomials of this isotopy. The polynomials $f- c_- +1$ and $\tilde f- c_- +1$ are then isotopic in the class of strictly Morse polynomials with only positive critical values.

2. By statement 1, we can assume that our isotopy is in the class of polynomials with only positive critical values. All surgeries allowed in the isotopy of strictly Morse polynomials are s4, s5, s6 and s7. However, surgeries s4 and s7 cannot be applied to polynomials with only real critical points, and surgeries s5 and s6 cannot appear in the isotopy with only positive critical values.

3. First, suppose that the generic polynomial $f$ has at least three local minima or at least  three local maxima. Let us order these three points of ${\mathbb R}^2$ in ascending order of their critical values. These points cannot lie on the same affine line, because the restriction of $f$ to any line is a polynomial of degree at most four. Therefore, their order defines an orientation of ${\mathbb R}^2$. This orientation remains throughout any isotopy of strictly Morse functions. However, the polynomial $f$ and its mirror image define different orientations. 

All generic polynomials of classes $X_9^+$ and $X_9^-$ with nine critical points satisfy  the above condition 
\begin{equation}
\label{condi4}
m_- \geq 3 \ \mbox{ or } \ m_+ \geq 3.
\end{equation} 

There is only one isotopy class of (not necessarily strictly) Morse polynomials of class $X_9^1$ with nine critical points that  does not satisfy this condition, see Fig.~\ref{X9011} (right). Each strictly Morse polynomial of this class is related to two lines in ${\mathbb  R}^2$. The first line passes through two local minima, while the second passes through two local maxima. Both lines are oriented from the critical point with the lower critical value towards the point with the higher critical value. These two lines always cross. Indeed, restricting our polynomial to the first line gives a polynomial of degree at most four with two Morse minima. Therefore, its degree is exactly four, and the line asymptotically connects two opposite domains at infinity in ${\mathbb R}^2$ in which the polynomial $f$ grows to $+\infty$. Similarly, the second line connects two domains in which $f$ falls to $-\infty$. Thus, these two oriented lines define an orientation of the plane. This orientation remains constant along the connected components of the set of strictly Morse polynomials, but differs for any two mirror-image polynomials.

There are exactly two  isotopy classes of Morse polynomials of class $X_9^2$ with nine critical points that  do not satisfy the condition (\ref{condi4}). Their passport invariants are $(2, 6, 1)$ and $(1, 6, 2)$. These two classes are obtained from each other by the multiplication by $-1$, therefore it is sufficient to consider only the first of them. 

\begin{proposition}
\label{au}
The polynomials of class $X_9^2$ cannot have two local minima and one local maximum, 
all three of which lie on the same affine line in ${\mathbb R}^2$. 
\end{proposition}

\noindent
{\it Proof.} Suppose that a polynomial $f$ of class $X_9^2$ has two local minima and one maximum on the same line. By choosing appropriate affine coordinates, we can assume that this line is given by the condition $y=0$, and the local maximum point is at the origin. By adding a constant function we can assume that $f(0)=0$. Let $ q< p$ be the $x$-coordinates of the two local minima of $f$. Clearly, $q <0 <p$. Consider the 
decomposition  of the polynomial $f$ by the powers of $y$:
\begin{equation}
f = f_4(x) + y f_3(x) + y^2f_2(x) + y^3 f_1(x) + f_0 y^4
\end{equation}
 where each $f_i$ is a polynomial of degree at most $i$ in one variable. Both degree three polynomials $f'_4$ and $f_3$ are divisible by $x(x-p)(x-q)$. Multiplying $f$ by a constant we can assume that $f'_4 = x(x-p)(x-q)$, and  
\begin{equation}f_4(x) = \frac{1}{4}x^4 - \frac{1}{3}(p+q)x^3 + \frac{1}{2} p q x^2.
\end{equation}

First, suppose that $f_3 \not \equiv 0$. Then, by dilating the $y$ coordinate, we can assume that $f_3= x(x-p)(x-q)$. Let $a, b,$ and $c$ be the coefficients of the polynomial
\begin{equation}
\label{equ2}
f_2(x) \equiv a x^2 + b x +c . 
\end{equation} 
The Hessian determinants of $f$ at the points $(0,0)$, $(p,0)$ and $(q,0)$ are positive, which gives us the system of inequalities 
\begin{eqnarray}
c - \frac{1}{2} p q  \leq 0\\
a p^2+b p+c- \frac{1}{2}p(p-q) & > & 0  \\
a q^2 + b q +c - \frac{1}{2}q (q-p) & >& 0.
\end{eqnarray}
These inequalities easily imply that \begin{equation}
\label{cont}
a > \frac{3}{2}.\end{equation}
The principal (degree four) homogeneous part of the polynomial $f$ has the form $\frac{1}{4}x^4 + x^3y + a x^2 y^2 + \varphi x y^3 + f_0 y^4$ with some $\varphi$. Since $f$ is of class $X_9^2$, the corresponding polynomial 
$\frac{1}{4}t^4 + t^3 + a t^2 + \varphi t + f_0$ in one variable $t = \frac{x}{y}$ has four real roots, and hence its second derivative $3 t^2 + 6 t + 2 a$ has two real roots. This contradicts inequality (\ref{cont}). 

If $f_3 \equiv 0$ then the degree two polynomial (\ref{equ2}) is positive at the points $p$ and $q$ and is negative at the origin, therefore $a>0$. The principal homogeneous part of $f$ has then the form $\frac{1}{4}x^4 + a x^2 y^2 + \varphi x y^3 + f_4 y^4$ with $a>0$. Such a polynomial also cannot have four lines of real solutions.
\hfill $\Box$ \medskip

\noindent
{\it End of the proof of Lemma \ref{lemel}.}
Consider two vectors in ${\mathbb R}^2$ which are directed from the minimum point of the polynomial $f$ to the maximum points and are ordered in the correspondence with the order of the critical values of these maximum points. According to Proposition \ref{au}, these ordered vectors define an orientation of ${\mathbb R}^2$. This orientation remains constant along the connected components of the space of strictly Morse polynomials of class $X_9^2$. Any reflection in a line changes this orientation. \hfill $\Box$ 
\medskip

\noindent
{\it End of the proof of Theorem \ref{mthst}.}
Consider an arbitrary isotopy class $\Omega$ of Morse polynomials
 with nine real critical points. 
According to statement 1 of Lemma \ref{lemel}, the number of isotopy classes of strictly Morse polynomials within $\Omega$ is equal to the number of isotopy classes of such polynomials with only positive critical values. 

The number of virtual Morse functions in  the virtual component
associated with the class $\Omega$ that have only positive critical values, is equal to $\frac{1}{10} \mbox{Card}(\Omega).$ By Proposition \ref{propmain} all of these virtual functions are represented by some polynomials of class $\Omega$. 
 By statement 2 of Lemma \ref{lemel}, each component of the set of strictly Morse polynomials with only positive critical values is associated with only one such virtual Morse function. Conversely, Proposition \ref{refprop} states that each such virtual Morse function is associated with one or two such components. If the class $\Omega$ is chiral, then this virtual Morse function is associated with two isotopy classes of Morse polynomials, only one of which belongs to the class $\Omega$.
If the class $\Omega$ is achiral, then according to statement 3 of Lemma \ref{lemel} there are two different components of the set of strictly Morse polynomials with this virtual Morse function, both of which belong to $\Omega$.
 \hfill $\Box$

\begin{corollary}
There are exactly $3200 = (732 + 246 + 622) \times 2$ 
$($respec\-ti\-vely, $272=136 \times 2)$
isotopy classes of strictly Morse polynomials of degree four with positive principal homogeneous parts, nine real critical points and five $($res\-pec\-tively, four$)$ local minima.

There are exactly 7044 $($respectively, 1812, 532, 1812, 7044$)$ isotopy classes of strictly Morse polynomials of degree four with principal homogeneous parts vanishing on two real lines, nine real critical points, and four $($respectively, three, two, one, none$)$ local minima.  

There are exactly 6084 isotopy classes 
of strictly Morse polynomials of degree four with principal homogeneous parts vanishing on four real lines and three points of local minima; the same holds true if minima are replaced by maxima. 
There are exactly 4052 classes of strictly Morse polynomials with such principal parts, nine real critical points and exactly two local minima; the same is true for the case of polynomials with one local minimum. \hfill $\Box$
\end{corollary}

\subsection{Normalization of $D$-graphs}
\label{norma}

\begin{definition} \rm
An edge of a D-graph is called {\em normal} if

a) it is oriented from a vertex with a smaller Morse index to a vertex with a larger index, and

b) it is solid if the parities of Morse indices of the critical points corresponding to its ends are different; it is dashed if these parities are the same. 

 Otherwise, this edge is called a {\em tunnel} edge. The {\em normalization} of a $D$-graph consists in removing all its tunnel edges.
\end{definition}

Analyzing Figs.~\ref{miss90a}--\ref{X902}, we see that

a) all $D$-graphs of polynomials with $m_+ >0$ shown in these figures are already normal,

b) nine out of ten $D$-graphs of polynomials with $m_+ = 0$ are split by the normalization into two standard Coxeter--Dynkin graphs of some simple singularities;

c) the tenth $D$-graph (see Fig.~\ref{DGX92C}) splits into an isolated vertex and the extended Coxeter--Dynkin graph of type $\tilde E_7$.
\medskip

A very similar situation occurs for $D$-graphs of degree three polynomials ${\mathbb R}^3 \to {\mathbb R}$, see \cite{Vspace}, and also for $D$-graphs of Morse perturbations of $J_{10}$ singularities, see \cite{VJ10}.

\section{Realization of polynomials with positive principal parts  and prescribed values of invariants} 

\label{reaax0}

\subsection{Realization of the D-graph of Fig.~\ref{miss90a}} 
\label{reala5a4}

\begin{proposition}
The principal $($of degree four$)$ homogeneous part of the polynomial 
\begin{equation}
F = x^2 + 2 x y^2 + y^4 + \frac{16 \sqrt{3}}{15} x^2 y + \frac{16\sqrt{3}}{15} x y^3 - \frac{76}{45} x^3 - \frac{2}{5} x^2 y^2 - \frac{16\sqrt{3}}{15} x^3 y + \frac{23}{30} x^4
\label{mfq0}
\end{equation}
is non-degenerate in the complex domain and positive definite in ${\mathbb R}^2$. The $j$-invariant of this principal part is  $245/3$.
\end{proposition}

\noindent 
{\it Proof.} The first statement is obviously equivalent to the positivity and non-dis\-cri\-mi\-nant\-ness of the polynomial
$$ \Phi(t) \equiv t^4 + \frac{16\sqrt{3}}{15} t^3 - \frac{2}{5} t^2 - \frac{16\sqrt{3}}{15} t + \frac{23}{30}.$$

By the Descartes' rule, $\Phi'$ has only one positive root, and hence $\Phi$ has only one local minimum in ${\mathbb R}_+$. The Taylor decomposition of $\Phi$ at the point $t=\frac{1}{2}$ in the coordinate $\tau= t-\frac{1}{2}$ has the form
\begin{equation}\alpha + \beta \tau + \gamma \tau^2 + \delta \tau^3 + \varepsilon \tau^4, \label{ll98}
\end{equation} where $\alpha > 0.0363$, $- 0.362 < \beta <0$, $\gamma > 3.87$, and $\delta $ and $\varepsilon $ are positive. In particular, the first three summands of (\ref{ll98}) are an everywhere positive quadratic polynomial.
 Since $\beta<0$, the positive local minimum point of $\Phi$ lies in the domain where $\tau>0$, but in this domain both remaining summands $\delta \tau^3$ and $\varepsilon \tau^4$ are positive. 

$\Phi'$ also has two negative roots in the intervals $(-2, -1)$ and $(-1, 0)$. Only the first of these roots is a local minimum of $\Phi$ and is of interest to us. By  substituting $s = - t^{-1}$ we reduce our problem to the positivity of the polynomial $\frac{23}{30} s^4 + \frac{16\sqrt{3}}{15} s^3 - \frac{2}{5} s^2 - \frac{16\sqrt{3}}{15} s + 1$
in the interval $(\frac{1}{2}, 1)$. Its Taylor decomposition at the point $s= \frac{1}{2}$ (in the coordinate $\sigma = s -\frac{1}{2}$) has the form 
$a + b \sigma + c \sigma^2 + d \sigma^3 + e \sigma^4$, where $a> 0.25$, $-0.48 < b <0$, $c> 3.5$, and $d$ and $e$ are positive. Therefore, this polynomial is positive for all positive $\sigma$. 

Finally, a positive polynomial of degree four in one variable can have multiple complex roots only if it is the square of a quadratic polynomial; it is easy to check that $\Phi$ is not such a square.

The value of the $j$-invariant follows from an explicit formula for it, see e.g. pages 26--29 in \cite{Mukai}; in the sequel we will not discuss the proofs of similar statements.
\hfill $\Box$ \medskip

The polynomial (\ref{mfq0}) has a critical point of type $A_5$ at the origin and  has  a local minimum there.
Indeed, in the coordinates $\check x = x + y^2$ and $y$ the sum of its monomials $\lambda \check x^u y^v$, $\lambda \neq 0,$ satisfying the inequality $ 3u + v \leq 6$ is $$\check x^2 - \frac{16 \sqrt{3}}{15} \check x y^3 + \frac{76}{45} y^6 \equiv \left(\check x - \frac{8\sqrt{3}}{15} y^3\right)^2 + \frac{188}{225} y^6 .$$ This  is a quasihomogeneous function of type $A_5$, so all the higher monomials can be killed by the choice of local coordinates, see \cite{AVGZ}, vol. I, \S 12.

The Taylor expansion of (\ref{mfq0}) at the point $(1, 0)$ in the coordinates $\tilde x \equiv x - 1$ and $y$ is 
$$ \frac{7}{90} + \frac{8}{15} \left(\tilde x - \sqrt{3} y\right)^2 + 
\frac{62}{45} x^3 - \frac{32 \sqrt{3}}{15} x^2 y + \frac{6}{5} x y^2 + \frac{16 \sqrt{3}}{15} y^3 + 
$$
$$
+\frac{23}{30} \tilde x^4 - \frac{16\sqrt{3}}{15} \tilde x^3 y - \frac{2}{5} \tilde x^2 y^2 + \frac{16\sqrt{3}}{15} \tilde x y^3 + y^4 .
$$
In the coordinates $X \equiv \tilde x - \sqrt{3}y$ and $y$ it looks as 
$$ \frac{7}{90} + \frac{8}{15} X^2 + 
\frac{1}{45} X 
\left( 36 y^2 + X ( 90 \sqrt{3} y +62 X ) \right) + 
$$
$$
+ \frac{3}{10} y^4 + \frac{1}{30} X \left( -4\sqrt{3}y^3 + X \left(114 y^2 + X(60y + 23 X)\right) \right) .
$$
The lowest quasihomogeneous part of this polynomial (modulo the constant term $\frac{7}{90}$) is  
\begin{equation}
\label{polo}
\frac{1}{30}(4 X + 3 y^2)^2, 
\end{equation} hence this polynomial has
a critical point of class $A_r$ with $r \geq 4$ at $(1, 0)$. 

By the index considerations, the sum of the Milnor numbers of all critical points of the function (\ref{mfq0})  in ${\mathbb C}^2$ is  9, so the Milnor number of the latter singularity is equal exactly to four, and so it is of type $A_4$. 

The space of all polynomials of degree four can be considered as a versal deformation of the function singularity defined by the principal homogeneous part of (\ref{mfq0}). Hence it is also a versal deformation of any multisingularity obtained by its sufficiently small perturbation contained in it: in particular, of our bisingularity consisting of critical points of type $A_5$ and $A_4$ of the function (\ref{mfq0}). Therefore, we can slightly perturb  the obtained polynomial in this space so that the first of these 
  critical points splits into two saddlepoints with equal critical values and three local Morse minima,  the second critical point splits into two saddlepoints (also with equal critical values) and two local minima, and all critical values at points obtained by splitting the $A_5$ singularity are lower than the critical values of points obtained from the $A_4$ singularity. 

\unitlength 0.5mm
\begin{figure}
\begin{picture}(120,65)
\bezier{800}(10,0)(20,0)(80,40)
\bezier{200}(10,0)(0,0)(0,10)
\bezier{200}(0,10)(0,15)(5,23)
\bezier{200}(5,23)(10,30)(5,37)
\bezier{200}(5,37)(0,45)(0,50)
\bezier{200}(0,50)(0,60)(10,60)
\bezier{400}(10,60)(20,60)(40,40)
\bezier{400}(40,40)(55,30)(80,40)
\bezier{100}(80,40)(85,45)(90,45)
\bezier{100}(80,40)(85,40)(90,45)
\bezier{100}(90,45)(95,50)(98,50)
\bezier{100}(90,45)(95,45)(99,47)
\bezier{50}(98,50)(100,49)(99,47)
\bezier{100}(20,20)(40,30)(20,40)
\bezier{200}(20,14)(35,30)(20,46)
\bezier{200}(20,20)(10,10)(20,14)
\bezier{200}(20,40)(10,50)(20,46)
\end{picture}
\caption{Level sets of a polynomial realizing Fig.~\ref{miss90a}}
\label{miss90breal}
\end{figure}
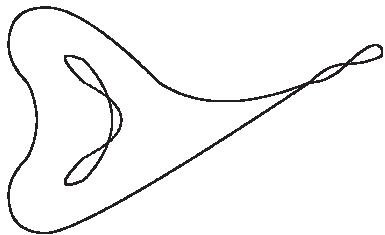

\begin{proposition}
The isotopy type of the pair of level sets of the obtained polynomial, corresponding to the critical values at its saddlepoints, is  shown in Fig.~\ref{miss90breal}.
\end{proposition}

\noindent 
{\it Proof.} The function (\ref{mfq0}) has only two critical points, therefore its critical level set $\{F=\frac{7}{90}\}$ is a compact one-component curve with a single cusp point of type $A_4$, and the interior part of this curve contains all points where $F< \frac{7}{90}$, including the other critical point of $F$. 
Such a curve in ${\mathbb R}^2$ can have only two isotopy types, see Fig.~\ref{csp}.
\begin{figure}
\unitlength 0.3 mm
\begin{picture}(60,60)
\bezier{300}(15,12)(-5,30)(15,48)
\bezier{300}(15,48)(30,60)(45,42)
\bezier{300}(15,12)(30,0)(45,18)
\bezier{300}(45,18)(55,30)(65,30)
\bezier{300}(45,42)(55,30)(65,30)
\end{picture} \qquad \qquad
\begin{picture}(60,60)
\bezier{300}(15,12)(-5,30)(15,48)
\bezier{300}(15,48)(30,60)(47,42)
\bezier{300}(15,12)(30,0)(47,18)
\bezier{300}(47,18)(60,30)(45,30)
\bezier{300}(47,42)(60,30)(45,30)
\end{picture}
\caption{Cuspidal ovals}
\label{csp}
\end{figure}
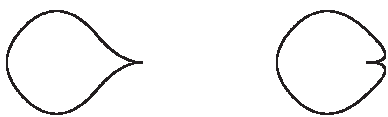

The quadratic part (\ref{polo}) of the Taylor expansion of our function at its cusp point is positive, so this function grows in the locally larger
component of the complement of this curve at the cusp point, and so we have the left-hand topological shape of Fig.~\ref{csp}. Applying the standard morsifications of two critical points gives us a shape isotopic to Fig.~\ref{miss90breal}.
\hfill $\Box$

\begin{proposition}
The D-graph of the obtained polynomial is shown in Fig.~\ref{miss90a}.
\label{proDg}
\end{proposition}

\noindent
{\it Proof.}
By the proved part of Theorem \ref{mthmdX90}, it can only be one of the three D-graphs shown in Figs. \ref{miss90a}, \ref{miss90b} and \ref{miss90c}. By the construction of our Morse polynomial, the set of nine vertices of this graph can be split into subsets of cardinalities 5 and 4 in such a way that, by removing all edges connecting points of different subsets we reduce our graph to two standard Coxeter-Dynkin graphs of types $A_5$ and $A_4$ with strictly alternating white and black vertices. Only one of the three allowed D-graphs admits such a splitting. 
\hfill $\Box$

\begin{proposition}
\label{achir}
The polynomial just constructed and the polynomial obtained from it by the
reflection of the plane with respect to any line are not isotopic in the class of Morse polynomials. 
\end{proposition}

\noindent
{\it Proof.} The D-graph of Fig.~\ref{miss90a} has no symmetries, therefore any order of its white vertices uniquely defines an order of local minima of any Morse function with this  D-graph. No non-discriminantal polynomial of degree four can have three Morse minima on the same affine line. Therefore,  if the number of its local minima is at least three (as in our case) then any order of them uniquely specifies an orientation of the plane, and isotopic Morse functions define the same orientation. Such orientations defined by any polynomial and its mirror image are opposite to each other. \hfill $\Box$

 \subsection{Realization of the D-graph of Fig.~\ref{miss90b}}

\label{+X9B}

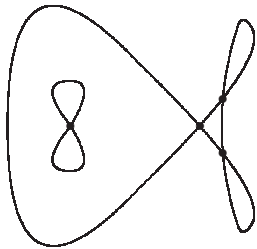
\begin{figure}
\unitlength 0.6mm
\begin{picture}(60,60)
\bezier{100}(3,20)(-0.5,30)(3,40)
\bezier{300}(30,30)(16,-15)(3,20)
\bezier{300}(3,40)(16,75)(30,30)
\bezier{100}(15,30)(7,20)(15,20)
\bezier{100}(15,30)(21,20)(15,20)
\bezier{100}(15,30)(7,40)(15,40)
\bezier{100}(15,30)(21,40)(15,40)
\put(15,30){\circle*{1.5}}
\put(30,30){\circle*{1.5}}
\end{picture} \qquad \qquad 
\begin{picture}(60,60)
\bezier{300}(30,44)(60,15)(55,8)
\bezier{70}(55,8)(53,5)(52,8)
\bezier{250}(52,8)(45,30)(52,52)
\bezier{70}(52,52)(53,55)(55,52)
\bezier{300}(55,52)(60,45)(30,16)
\bezier{300}(30,16)(1,-15)(1,30)
\bezier{300}(1,30)(1,75)(30,44)
\bezier{100}(15,30)(7,20)(15,20)
\bezier{100}(15,30)(21,20)(15,20)
\bezier{100}(15,30)(7,40)(15,40)
\bezier{100}(15,30)(21,40)(15,40)
\put(15,30){\circle*{1.5}}
\put(43.7,30){\circle*{1.5}}
\put(48.7,36){\circle*{1.5}}
\put(48.7,24){\circle*{1.5}}
\end{picture}
\caption{Level sets for a realization of Fig.~\ref{miss90b}}
\label{X901}
\end{figure}

For any $A \in (-1, \infty)$ and positive $\varepsilon$, the polynomial 
\begin{equation}
 x^4 +2A x^2 y^2 + y^4 +\varepsilon x^3 
\label{x9b}
\end{equation}
 has a positive principal homogeneous part (non-discriminantal if $A \neq 1$) and a critical point of type $E^+_6$ at the origin. The zero level set of this function is an oval with a single singular point at the origin, see the outer curve in Fig.~\ref{X901} left. If $A \in (-1,0)$ then the function (\ref{x9b}) has additionally two minima and one saddlepoint inside the domain of its negative values bounded by this oval. The level set of this saddlepoint looks topologically as shown by the inner curves in both parts of  Fig.~\ref{X901}.

The small perturbation of the $E_6$ singularity described on p. 16 of \cite{AC} or in \cite{GZ} splits it into three Morse saddlepoints with zero critical values and three Morse minima with lower values. The global zero level set of the obtained function is topologically situated as the outer curve in Fig.~\ref{X901} right. Thus, two level sets of the obtained polynomial, containing all its saddlepoints, look topologically as shown in this figure.

\begin{proposition}
The D-graph of the obtained Morse function in the case $A<0$ is as shown in Fig.~\ref{miss90b}.
\end{proposition}

\noindent
{\it Proof} is analogous to the proof of Proposition \ref{proDg}: only one of the three possible D-graphs  shown in  figures \ref{miss90a}, \ref{miss90b} or \ref{miss90c} can be split into the standard Coxeter--Dynkin diagrams of 
types $A_3$ and $E_6$. \hfill $\Box$ \medskip

If $A>0$, then the polynomial (\ref{x9b}) has only one local minimum in the domain of its negative values. 
If $A=0$, then it has a point of type $A_3$ there, which can be further perturbed in any of the ways arising in the cases $A<0$ and $A>0$. In particular, the polynomial 
$$x^4 +2A x^2 y^2 + y^4 $$ 
is approached by the stratum $\{A_3  +E_6\}$ of the function space if $A=0$, and is not for the neighboring values of $A$; see Theorem \ref{tabadjp} on page \pageref{tabadjp}.

\smallskip
The $j$-invariant of the principal part of the polynomial (\ref{x9b}) with $A=0$ is  1.

\subsection{Realization of the D-graph of Fig.~\ref{miss90c}}
\label{reala7a2}

For any $A \in (-1, 1),$ consider the polynomial 
\begin{equation}
x^4 +2A x^2 y^2 + y^4 + x^2 - 2 x y^2 - 2 A x^3 . 
\label{x9p3r}
\end{equation}
In the coordinates $z \equiv x - y^2$ and $y$, its principal quasihomogeneous part is  $z^2 - 2A z y^4 + y^8. $ The additional substitution $\tilde z \equiv z - A y^4 $ makes it a polynomial with the principal quasihomogeneous part $\tilde z^2 + (1- A^2) y^8$, i.e. with a local minimum point of class $A_7$ at the origin. 

Also, the polynomial (\ref{x9p3r}) has two critical points in ${\mathbb C}^2$ whose coordinates satisfy the conditions $2 x^2 - 3 A x +1=0$ and $y=0$. 
If $A^2 > \frac{8}{9}$ then these two points are real: if $A<0$ then they are a minimum and a saddlepoint, and if $A>0$ then they are a maximum and a saddlepoint.

For $A^2 = \frac{8}{9}$ these points collide and form a real critical point of type $A_2$.

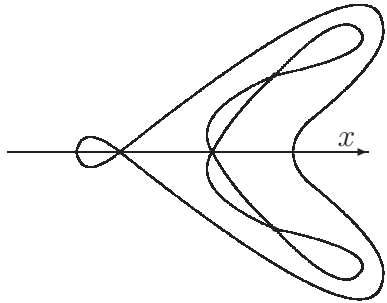
\begin{figure}
\unitlength=0.4mm
\begin{picture}(125,120)
\put(0,60){\vector(1,0){120}}
\put(110,62){$x$}
\put(37.5,60){\circle*{1.5}}
\put(68.15,60){\circle*{1.5}}
\put(89,34.4){\circle*{1.5}}
\put(89,85.5){\circle*{1.5}}
\bezier{200}(68.15,60)(60,73)(89,85.5)
\bezier{200}(68.15,60)(60,47)(89,34.4)
\bezier{200}(68.15,60)(75,73)(89,85.5)
\bezier{200}(68.15,60)(75,47)(89,34.4)
\bezier{200}(89,85.5)(112,110)(118,99)
\bezier{200}(89,34.4)(112,10)(118,21)
\bezier{200}(89,85.5)(120,92)(118,99)
\bezier{200}(89,34.4)(120,28)(118,21)
\bezier{800}(37.5,60)(125,130)(125,100)
\bezier{800}(37.5,60)(125,-10)(125,20)
\bezier{300}(125,100)(125,90)(100,70)
\bezier{300}(125,20)(125,30)(100,50)
\bezier{250}(100,70)(90,60)(100,50)
\bezier{100}(37.5,60)(25,70)(23,60)
\bezier{100}(37.5,60)(25,50)(23,60)
\end{picture}
\caption{Level sets of a Morse function with  D-graph of Fig.~\ref{miss90c}}
\label{miss90level}
\end{figure}

By slightly perturbing the critical point  of type $A_7$ at the origin, we can split it into seven Morse points (four minima and three saddlepoints) in such a way that all three saddlepoints have the same critical value 0. The topological disposition of some two level sets of such a perturbation in the case $A \in \left(-1, -\frac{\sqrt{8}}{3}\right)$ is shown in Fig.~\ref{miss90level}. The value $0$ of the polynomial on the outer level curve is higher than the value on the inner one. Indeed, the self-intersection points of these two curves lying on the horizontal axis $\{y=0\}$ are respectively the local maximum and a local minimum of the restriction of the function (\ref{x9p3r}) to this axis.

By the construction, the D-graph of the obtained Morse function should contain the standard Coxeter-Dynkin graph of type $A_7$ as a subgraph, therefore it can be only the D-graph of Fig.~\ref{miss90c}.

\begin{remark} \rm
 We see that both polynomials $$x^4 + 2 A x^2 y^2 + y^4 $$ of type $X_9^+$ with $A^2 = \frac{8}{9}$ are approached by the stratum $\{A_7+ A_2\}$ of the function space, cf. Theorem \ref{tabadjp} on page \pageref{tabadjp}.  
The $j$-invariants of these two polynomials are $\frac{5^3 7^3}{3^5}$.
\end{remark}

\subsection{Realization of degree four polynomials with positive principal homogeneous parts and nine real critical points, one of which is a local maximum}
\label{withmax}

The function 
\begin{equation} 
\label{197}
x^4 + y^4 - \varepsilon x^2 - \varepsilon y^2
\end{equation} provides such a realization. Moreover, every positive definite homogeneous polynomial of degree four in ${\mathbb R}^2$ is the product of two positive quadratic functions, which are non-proportional if our polynomial is non-discriminant. We can perturb these quadratic functions by subtracting constant terms in such a way that the corresponding zero sets in ${\mathbb R}^2$ will be two ellipses with four common points. The intersection of their inner parts then contains a local maximum point. Therefore, the desired perturbations exist for every positive  principal part. The D-graph of the obtained polynomial can be computed by the Gusein-Zade--A'Campo method (\cite{AC}, \cite{GZ}) and is shown in Fig.~\ref{X901mx}.

\subsection{Proof of Proposition \ref{refprop} for polynomials of class $X_9^+$}
\label{proofref0}

\begin{lemma}
\label{lem77}
The stratum $\{E_6 + A_3\}$  in the space of polynomials with principal parts of type $X_9^+$  consists of a single connected component. Every polynomial in this stratum has the form 
\begin{equation}x^4 + y^4 + \varkappa \cdot x^3, \quad \varkappa > 0, \label{stf} \end{equation} 
in appropriate affine coordinates. The set of real polynomials $($\ref{vers0}$)$ with the principal homogeneous parts $\tilde x^4 + 2 A \tilde x^2 \tilde y^2 + \tilde y^4,$ $|A|<1$, having the form $($\ref{stf}$)$ in appropriate affine coordinates 
and critical value $-\frac{27}{256}$,
consists of four polynomials $x^4  + y^4 \pm x^3$ and $x^4 + y^4 \pm y^3$, where, respectively, $(x, y) = (\tilde x - \frac{1}{4}, \tilde y)$   and $(x, y) = (\tilde x, \tilde y - \frac{1}{4})$.
\end{lemma}

\noindent{\it Proof.} The 3-jet of an $E_6$ singularity in some affine coordinates is $x^3$. In the same coordinates, the coefficient at $y^4$ of an $X_9^+$-polynomial should be positive; let us make it equal to $1$ by dilating  the coordinate. According to Bezout's theorem, the other critical point (of type $A_3$) cannot lie in the line $\{x=0\}$.
Using an appropriate substitution $\tilde x = x+ \theta y$, let us place it on the axis $\{y=0\}$; let $\alpha$ be its $x$-coordinate.
Thus, in these coordinates, our function $\Phi$ has the form
$$x^3 + a x^4 + b x^3 y + c x^2 y^2 + d x y^3 + y^4, \quad a>0 .$$

 By a characteristic property of $A_k$ singularities (see \cite{LL}, Lemma 10.4.5), there exists a smooth parametrized curve through the point $(\alpha, 0)$, along which the value $|\mbox{grad }\Phi|$ decreases as the third degree of the parameter. According to Bezout's theorem, this curve cannot be tangent to the line $\{y=0\}$. Let us choose the coordinate $y$ as its parameter, then the curve has the equation $$x(y) = \alpha + \beta y + \gamma y^2+ ... $$ for some $\beta$ and $\gamma$. The first three terms of the power expansions of the restrictions of $\Phi'_x$ and $\Phi'_y$ to this curve should vanish. This gives us six equations on seven coefficients $a,$ $ b,$ $ c, $ $d, $ $\alpha, $ $\beta,$ $ \gamma$. It is easy to compute that this system of equations is equivalent to the conditions $b=c=d=0=\beta=\gamma$, $3+4 a \alpha =0$. Any positive value of $a$ is satisfactory, so $f$ has the form $x^3 + y^4 + a x^4,$ which can be reduced to the form (\ref{stf}) by dilating the $x$ coordinate. 

The set of all polynomials of this form is connected. Reducing the original polynomial to this form, we used only the dilation of coordinates with positive coefficients and the reduction of the 3-jet of an $E_6$ singularity to the form $x^3$, which can be done by an orientation-preserving affine transformation of the plane. Therefore, the whole stratum $\{E_6+A_3\}$ is also path-connected in the space of polynomials with the principal parts of type $X_9^+$.

The last statement of Lemma \ref{lem77}
follows from the fact that $x^4+y^4$ is the unique 
polynomial of the form $x^4 + 2A x^2 y^2 +y^4$ with $|A|<1$, whose
$j$-invariant is  1. \hfill $\Box$ \medskip
 
Applying to the polynomials $f$ and $\tilde f$ from Proposition \ref{refprop} appropriate orientation-preserving affine transformations (not changing the connected component of the space of generic functions), we can assume  that they both belong to the parameter space 
$\Theta \equiv ({\mathbb C}^1 \setminus \{1, -1\}) \times {\mathbb C}^8$ of the deformation (\ref{vers0}) of $X_9^+$ singularities, see the proof of Proposition \ref{propmain}.

Using the Lyashko--Looijenga covering $\Lambda$ over $B({\mathbb C}^1,9)$, we can assume that $f$ and $\tilde f$ have  the same collection of critical values and the same system of paths defining their vanishing cycles. 
Consider an arbitrary path $I:[0, 1] \to \Theta$ consisting of real polynomials with the principal parts of type $X_9^+$, such that
\begin{enumerate}
\item $I(0) = f$;
\item $I(1) $ is the polynomial  \begin{equation}x^4 + y^4 + x^3  \label{nfa}
\end{equation}
of type $\{E_6+A_3\}$ with critical values 0 and $-\frac{27}{256}$, 
\item the path $I$ has only finitely many transversal intersections with the variety of non-generic polynomials, 
\item for $t \in (0,1)$ sufficiently close to $1$, the points $I(t)$ are Morse perturbations of the polynomial (\ref{nfa})  having only real critical points and realizing the $D$-graph of Fig.~\ref{miss90b}.
\end{enumerate}

 Let $\Lambda \circ I$ be the projection of the path $I$ to the space $\mbox{Sym}^9({\mathbb C}^1)$ by the Lyashko--Looijenga map, and $\tilde I$ be the lifting of this path $\Lambda \circ I$ to the space $\Theta$ with the starting point $\tilde f$, along which the polynomials undergo all the same elementary surgeries as the polynomials from the path $I$. By the construction, this path $\tilde I$ lies in the real part ${\mathbb R}^9 \cap \Theta$ of this space, and its final point is also a function with critical points of type $E_6$ and $A_3$ and corresponding critical values 0 and $-\frac{27}{256}$. By Lemma \ref{lem77}, both endpoints $I(1)$ and $\tilde I(1)$ of the paths $I$ and $\tilde I$ are polynomials of the form $x^4 +y^4 \pm x^3$ or $x^4 + y^4 \pm y^3$.
All these four polynomials belong to the same connected orbit of the action of the group of affine transformation of ${\mathbb R}^2$. Applying an appropriate element of this group to these endpoints and to the entire paths $I$ and $\tilde I$ including their starting points $f$ and $\tilde f$, we can assume that both these endpoints coincide with the polynomial $x^4 + y^4 + x^3$. A neighborhood of this point in $\Theta$ is the space of a miniversal deformation of the multisingularity of type $\{E_6 + A_3\}$. It has the structure of the direct product of the parameter spaces of miniversal deformations of the singularities $E_6$ and $A_3$: the projection maps of this product are the maps inducing this neighborhood (considered separately as a deformation of the $E_6$ singularity only and the $A_3$ singularity only) from these miniversal deformations. 

Let $f_-$ and $\tilde f_-$ be the near endpoints of the paths $I$ and $\tilde I$, which have the same sets of critical values. By the construction, the virtual Morse functions associated with them are also the same, in particular they have equal Coxeter--Dynkin subgraphs, whose vertices are the sets of 6 or 3 points obtained by the decompositions of two critical points of (\ref{nfa}). Then by \cite{Liv} the projections of these two points to the parameter space of the miniversal deformation of an $E_6$ singularity belong to the same orbit of the group of automorphisms of the Lyashko--Looijenga covering associated with this singularity. This group contains only two elements preserving the set of real functions: the identical map and the involution defined by the reflection $y \leftrightarrow -y$ of the normal form of the $E_6$ singularity from Table \ref{t1}.  This involution acts on the space of vanishing cycles related to this singularity according to the only non-trivial involution of the canonical Coxeter-Dynkin diagram $E_6$, i.e. permuting two pairs of vanishing cycles and leaving the remaining two vanishing cycles unmoved. 
Analogously, the projections of our two points to the miniversal deformation of an $A_3$ singularity either coincide or are obtained from each other by an involution of the corresponding parameter space inducing the non-trivial automorphism of the $A_3$-graph. 

It cannot happen that the projections of our two points $f_-$ and $\tilde f_-$ to two spaces of miniversal deformations of  $E_6$ and $A_3$ singularities coincide in one case and do not  in the other. Indeed, the permutation of vertices of the D-graph of Fig.~\ref{miss90b}, which keeps one of these subgraphs fixed and acts non-trivially on the other,  damages the entire D-graph. If in both cases these projections coincide then $f_- =\tilde f_-$, and by continuity also $f = \tilde f$. There is only one alternative to this coincidence, when projections are different in both cases (and mapped into each other by the only real non-trivial automorphism of the Lyashko--Looijenga coverings). 
But such an alternative situation is known to us: it occurs when $f_-(x,y) \equiv \tilde f_-(x,-y)$. In this case, the same identity 
$f(x,y) \equiv \tilde f(x, -y)$ holds by the continuity for the initial points of our  paths. \hfill $\Box$

\begin{corollary}
The sets of degree four polynomials with the D-graphs shown in Figs.~\ref{miss90b}, \ref{miss90c} and \ref{X901mx} consist of single connected components.
\end{corollary}

\noindent
{\it Proof.} The polynomials constructed in \S \ref{+X9B} (see Fig.~\ref{X901}), \S~\ref{reala7a2} (see Fig.~\ref{miss90level}) and \S~\ref{withmax} (formula (\ref{197})), which represent these   sets,
are invariant under some reflections. \hfill $\Box$

\subsection{Morse polynomials with positive principal part, seven real critical points and no maxima}
\label{x+7}

The algorithm of \S \ref{onproofs} proves that there are 
 only two virtual components associated with such polynomials, and they are represented by two virtual Morse functions given in Table~\ref{virtmatr} on page \pageref{virtmatr}.

The corresponding isotopy classes of Morse polynomials are realized as follows.

A. The polynomial (\ref{x9p3r}) with $A^2 < \frac{8}{9}$ has a critical point of type $A_7$  at the origin and two non-real critical points.
Take the standard small perturbation of its critical point at the origin, splitting it into four local minima and three saddlepoints (so that the corresponding level set looks like the inner curve in Fig.~\ref{miss90level}).

B. Repeat the construction from \S \ref{+X9B}, where the parameter $A$ is positive, that is, the zero level set of the resulting function looks like the outer curve in Fig.~\ref{X901} right, but its left-side loop surrounds a single (global) minimum of our function. Thus, this zero level set bounds four different domains of negative values, each of which contains exactly one point of local minimum.

\begin{proposition}
Two polynomials constructed in two previous paragraphs 
belong to different connected components of 
the space of Morse polynomials of degree four with strictly positive principal homogeneous parts. 
\end{proposition}

\noindent
{\it Proof.} All critical points of the function constructed in paragraph A lie in the zero set of the coordinate $\tilde z \equiv x - y^2 - A y^4 $, see \S \ref{reala7a2}. This set is a convex smooth curve near the origin. In particular, the convex hull of four local minima in ${\mathbb R}^2$ is a quadrilateral. In the case B one of the local minima of the constructed function belongs to the triangle in ${\mathbb R}^2$ spanned by other three local minima. 
Therefore, every path connecting our two polynomials in the space of Morse functions contains a polynomial with three local minima on the same affine line in ${\mathbb R}^2$.
A polynomial function of degree four with a positive principal part cannot have such three minima on a line because its restriction to any line is a polynomial of degree exactly four in one variable. \hfill $\Box$ 
\medskip

These two polynomials are self-symmetric in ${\mathbb R}^2$, so there are only these two connected components with the passport $(4,3,0)$.

\subsection{Realization of polynomials with positive principal parts and at most five real critical points, or with seven real critical points and one local maximum}

\begin{figure}
\unitlength 0.39mm
\begin{picture}(60,60)
\bezier{300}(30,44)(60,15)(55,8)
\bezier{70}(55,8)(53,5)(49,8)
\bezier{250}(49,8)(20,30)(49,52)
\bezier{70}(49,52)(53,55)(55,52)
\bezier{300}(55,52)(60,45)(30,16)
\bezier{300}(30,16)(1,-15)(1,30)
\bezier{300}(1,30)(1,75)(30,44)
\put(43.7,30){\circle*{1.5}}
\put(36.5,37.5){\circle*{1.5}}
\put(36.5,22){\circle*{1.5}}
\end{picture} \quad
\begin{picture}(60,60)
\bezier{70}(57,12)(53,5)(49,8)
\bezier{250}(49,8)(20,30)(49,52)
\bezier{70}(49,52)(53,55)(57,48)
\bezier{300}(30,11)(62,35)(57,12)
\bezier{300}(30,49)(62,25)(57,48)
\bezier{300}(30,11)(1,-15)(1,30)
\bezier{300}(1,30)(1,75)(30,49)
\put(39.2,42.3){\circle*{1.5}}
\put(39.2,17.7){\circle*{1.5}}
\end{picture} \quad 
\begin{picture}(60,60)
\put(20,30){\circle{40}}
\put(40,30){\circle{40}}
\put(30,16){\circle*{1.5}}
\put(30,44){\circle*{1.5}}
\end{picture} \quad
\begin{picture}(60,30)
\put(15,30){\oval(30,30)[l]}
\put(45,30){\oval(30,30)[r]}
\bezier{300}(15,15)(20,15)(30,30)
\bezier{300}(15,45)(20,45)(30,30)
\bezier{300}(45,15)(40,15)(30,30)
\bezier{300}(45,45)(40,45)(30,30)
\put(30,30){\circle*{1.5}}
\end{picture} \quad \
\begin{picture}(60,60)
\bezier{400}(35,2)(70,20)(45,40)
\bezier{400}(35,58)(70,40)(45,20)
\bezier{300}(35,2)(1,-15)(1,30)
\bezier{300}(1,30)(1,75)(35,58)
\bezier{250}(45,40)(25,50)(23,30)
\bezier{250}(45,20)(25,10)(23,30)
\put(53.5,30){\circle*{1.5}}
\end{picture}
\caption{Perturbations of $X_9^+$ singularities with few critical points}
\label{X901a}
\end{figure}
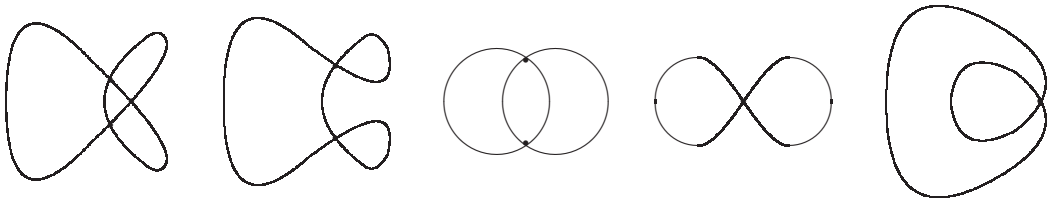

Consider the construction from \S \ref{+X9B} with $A>0$, but perturb the critical point of type $E_6$ arising on its first step in a different standard way, so that the zero level set of the obtained polynomial is as shown in Fig.~\ref{X901a} left.

This polynomial has seven real critical points: three minima, three saddlepoints and one maximum. 

Further, this critical point of type $E_6$ can be perturbed in other ways, so that the numbers $(m_-, m_{\times}, m_+)$ of minima, saddlepoints and maxima arising from it take any values $(2, 2, 0),$ $(1, 2, 1),$ $ (1, 1, 0),$ $(0, 1, 1) $ and $(0, 0, 0)$, see \S 2.5 in \cite{Vcau}. Adding to them the left-side minimum point, we realize all non-zero cells of the left part of Table~\ref{TX+} with $M \leq 5$. 
The shapes of the zero sets of some representatives of these classes are shown in Fig.~\ref{X901a}. All of them are invariant under some reflections and thus represent single connected components of the space of Morse polynomials.

\section{Realization of polynomials with two real asymptotes ($X_9^1$)} 

All four possible D-graphs for this case shown in Figs.~\ref{DGX921A}--\ref{DGX921D} have no symmetries, so by the arguments of Proposition \ref{achir} in each of the next four subsections we prove the existence of two different  isotopy classes of Morse polynomials with $(M, m_{+}) = (9, 0)$: one that we actually construct, and the other that is mirror symmetric to it.

 \subsection{Realization of Fig.~\ref{DGX921A} }
\label{reald5a4}

\begin{proposition}
The fourth degree homogeneous part of the polynomial 
\begin{equation}
-\frac{4}{3} x^2 y + y^4 + \frac{32}{5}x y^3 + 24 x^2 y^2 - 32 x^3 y -2 x^4
\label{ex91a}
\end{equation}
vanishes on exactly two lines. 
The $j$-invariant of this homogeneous part is  $-\frac{5}{3}$.
\end{proposition} 

\noindent
{\it Proof.} By the Descartes' rule the corresponding polynomial $$t^4 + \frac{32}{5}t^3 + 24 t^2 -32 t -2$$
in one variable has exactly one positive root. All its negative roots belong to the interval $(-1,0)$, since both its summands $t^4 + \frac{32}{5}t^3 + 24 t^2$ and $-32t-2$ are strictly positive on the domain $(-\infty, -1]$. A trivial estimate shows that its derivative is strictly negative on $(-1,0)$, so it also has only one negative root.
\hfill $\Box$ \medskip

The principal (lower) quasihomogeneous part of the function (\ref{ex91a}) at the origin is  $-\frac{4}{3}x^2y + y^4$, so it has a critical point of type $D_5$ there. It also has a critical point with the coordinates $(x_0, y_0) = (-\frac{5}{3^5}, \ \frac{5}{3^5})$. Its Taylor expansion at this point in the coordinates $\tilde x = x-x_0, \tilde y=y-y_0$ is equal (assuming the notation $\tau \equiv \frac{5}{3^5}$) to 
$$
-\frac{1}{3} \tau^3 + \frac{2\tau}{9}\left( 4 \tilde x^2 - 4 \tilde x \tilde y + \tilde y^2\right) - \frac{12 \tau}{5} \left(10 \tilde x^3 - 33 \tilde x^2 \tilde y + 12 \tilde x \tilde y^2 + \tilde y^3 \right) + 
$$
$$
+ \tilde y^4 + \frac{32}{5} \tilde x \tilde y^3 + 24 \tilde x^2 \tilde y^2 - 32 \tilde x^3 \tilde y -2 \tilde x^4.$$
In coordinates $\tilde x$ and $Y = \tilde y - 2 \tilde x$  this expression equals $$- \frac{\tau^3}{3} + \frac{2\tau}{9} Y^2 - \frac{12\tau}{5} Y (27 \tilde x^2 + 18 \tilde x Y + Y^2) + $$
$$+ \frac{486}{5} \tilde x^4 + Y \left(\frac{864}{5} x^3 + Y \left(\frac{432}{5} \tilde x^2 + Y \left( Y + \frac{72}{5} \tilde x \right) \right) \right).
$$
The principal (lower) quasihomogeneous part of its non-constant part in these coordinates with the weights $(2, 1)$ is $$
\frac{2 \tau}{9} Y^2 - \frac{12 \cdot 27 \tau}{5} Y \tilde x^2 + \frac{486}{5} \tilde x^4 \equiv \frac{2}{5 \cdot 3^7}(5 Y - 3^6 \tilde x^2)^2. $$ In particular, it is degenerate. Therefore, the Milnor number of this critical point is at least 4. It cannot be greater than 4 because the sum of the Milnor numbers of all critical points of the polynomial (\ref{ex91a}) in ${\mathbb C}^2$ is 9. This critical point has a non-zero quadratic part (with one positive square), hence it is of type $A_4$. If we  independently split two critical points  of types $D_5$ and $A_4$ of our function as shown (up to a sign)  on pp. 13--14 of \cite{AC}, we obtain a function with four local minima and five saddlepoints.

\begin{figure}
\unitlength 0.5mm
\begin{picture}(120,100)
\bezier{300}(0,65)(50,45)(55,70)
\bezier{350}(55,70)(65,120)(80,40)
\bezier{400}(80,40)(86,85)(120,100)
\bezier{600}(20,0)(55,60)(120,30)
\bezier{300}(0,55)(50,40)(60,65)
\bezier{200}(60,65)(55,40)(40,30)
\bezier{150}(40,30)(15,10)(10,0)
\bezier{400}(120,85)(80,50)(120,36)
\end{picture} \qquad \qquad 
\begin{picture}(120,100)
\bezier{300}(0,65)(50,45)(55,70)
\bezier{350}(55,70)(65,120)(80,35)
\bezier{400}(80,35)(86,85)(120,100)
\bezier{50}(80,35)(83,27)(80,27)
\bezier{50}(80,35)(77,27)(80,27)
\bezier{600}(20,0)(55,60)(120,30)
\bezier{300}(0,55)(50,40)(60,65)
\bezier{200}(60,65)(60,40)(40,30)
\bezier{150}(40,30)(15,10)(10,0)
\bezier{40}(60,65)(62,68)(60,71)
\bezier{40}(60,65)(58,68)(60,71)
\bezier{40}(60,71)(63,74)(60,74)
\bezier{40}(60,71)(58,74)(60,74)
\bezier{400}(120,85)(80,50)(120,36)
\end{picture}
\caption{Critical level sets of the polynomial (\ref{ex91a}) and of its perturbation  realizing Fig.~\ref{DGX921A}}
\label{??A}
\end{figure}

\begin{proposition}
The  D-graph of the obtained polynomial is shown in Fig.~\ref{DGX921A}.
\end{proposition}

\noindent
{\it Proof} is analogous to the proof of Proposition \ref{proDg}: only one of  the D-graphs shown in Figs.~\ref{DGX921A}, \ref{DGX921B}, \ref{DGX921C}, and \ref{DGX921D} admits a splitting into the standard Coxeter--Dynkin diagrams of types $D_5$ and $A_4$.  \hfill $\Box$

\subsection{Realization of Fig.~\ref{DGX921B}}
\label{reala6a3}

\begin{proposition}
The principal $($of degree four$)$ homogeneous part of the polynomial 
\begin{equation}
F = (3 x + y^2)^2 + 4 x^2 y + \frac{4}{3} x y^3 - \frac{40}{3} x^3 - 4 x^2 y^2 -4 x^3 y + \frac{11}{2} x^4
\label{mfq}
\end{equation}
 is non-discriminantal in ${\mathbb C}^2$ and vanishes on exactly two real lines. 
The $j$-invariant of this homogeneous part is $-49$.
\end{proposition}

\noindent 
{\it Proof.} The first  statement is equivalent to the assertion that the polynomial $$\varphi = t^4 + \frac{4}{3} t^3 - 4 t^2 - 4 t + \frac{11}{2} $$ in one variable has only two real roots. By the Descartes' rule, it has no more than two positive roots, and indeed has them, as $\varphi(0)>0$ and $\varphi(1)<0$. It remains to prove that its minimal value in the domain $\{t < 0\}$ is positive. It is easy to check that its derivative has two negative roots, one in the interval $(-2, -\frac{3}{2})$ and the other in $(-1,0)$. Only the first one corresponds to a local minimum of $\varphi$ and is interesting for us. Expressing $\varphi $ in the coordinate $\tau \equiv - \frac{3}{2} - t$ (after which our suspicious interval consists of values $\tau \in (0, \frac{1}{2})$) we get the polynomial $$ \tau^4 + \frac{14}{3} \tau^3 + \frac{7}{2} \tau^2 - \frac{7}{2} \tau + \frac{49}{16}.$$ Its last three summands form an everywhere positive polynomial, and the remaining two monomials are positive for all positive $\tau$. \hfill $\Box$

\begin{proposition}
The polynomial $($\ref{mfq}$)$ has a critical point of type $A_6$ at the origin and a critical point of type $A_3$ with the coordinates $(1, 0).$
\end{proposition}

\noindent
{\it Proof.} The substitution $\check x = x + \frac{1}{3}y^2$ turns it into a polynomial, whose lower quasihomogeneous (with the weights $(3, 1)$) part is $\left(3\check x - \frac{2}{9}y^3\right)^2$, in particular it is of type $A_k$ with $k \geq 6$.

The Taylor expansion of the polynomial  $($\ref{mfq}$)$ at the point $(1,0)$ (in the coordinates $\tilde x=x - 1$ and $y$) is
\begin{equation}
\label{76}
 \frac{7}{6} + 2 (\tilde x-y)^2 + \frac{1}{3}(\tilde x-y) (26 \tilde x^2 + 2 \tilde x y - 4 y^2) + 
\end{equation}
$$
+ \frac{4}{3} \tilde x y^3 - 4 \tilde x^2 y^2 - 4 \tilde x^3 y + \frac{11}{2} \tilde x^4.
$$
In the local coordinates $X \equiv \tilde x - y $ and $y$, its non-constant part contains only monomials $X^a y^b$ with $2a + b \geq 4$, and its quadratic part is non-trivial, so it has a singularity of type $A_r$ with $r \geq 3$. The sum of the Milnor numbers of all critical points of the polynomial $($\ref{mfq}$)$ in ${\mathbb C}^2$ is equal to 9, so the Milnor numbers of these two critical points are indeed 6 and 3. \hfill $\Box$ \medskip

By the index arguments, this function at the second critical point (of type $A_3$) is equal to $\frac{7}{6} + \xi^2 - \eta^4$ in some local coordinates. Indeed, the quadratic part of (\ref{76}) is non-negative, and in the case of the sign $+$ in front of $\eta^4$ we would have a local minimum point with the local index of the gradient vector field equal to 1. Together with the index 0 of the $A_6$ singularity,
this contradicts the index of the polynomial  $($\ref{mfq}$)$ at the infinity (which is  $-1$).

By perturbing independently two critical points of this polynomial, we can obtain a Morse function with four local minima and five saddlepoints. 

\begin{figure}
\unitlength 0.5mm
\begin{picture}(120,120)
\bezier{500}(60,60)(50,35)(0,10)
\bezier{500}(60,60)(55,35)(20,0)
\bezier{600}(120,110)(70,50)(110,120)
\bezier{800}(120,100)(65,60)(30,0)
\bezier{600}(100,120)(77,44)(64,60)
\bezier{600}(64,60)(25,140)(45,75)
\bezier{600}(45,75)(55,50)(0,15)
\end{picture} \qquad \qquad \qquad
\begin{picture}(120,120)
\bezier{500}(55,53)(45,37)(0,10)
\bezier{500}(53,53)(45,30)(20,0)
\bezier{50}(53,53)(54,55)(55,55)
\bezier{50}(55,53)(55,54)(55,55)
\bezier{50}(55,55)(55,59)(58,59)
\bezier{50}(55,55)(58,56)(58,59)
\bezier{50}(58,59)(58,62)(60,62)
\bezier{50}(58,59)(60,59)(60,62)
\bezier{600}(120,110)(70,50)(110,120)
\bezier{800}(120,100)(65,60)(30,0)
\bezier{600}(100,120)(80,40)(64,60)
\bezier{600}(64,60)(25,140)(45,75)
\bezier{600}(45,75)(55,50)(0,15)
\end{picture} 
\caption{Critical level sets of the polynomial (\ref{mfq}) and of its perturbation  realizing Fig.~\ref{DGX921B}}
\label{??B}
\end{figure}

\begin{proposition}
The  D-graph of the polynomial constructed in the previous paragraph is shown in Fig.~\ref{DGX921B}.
\end{proposition}

\noindent
{\it Proof} is analogous to that of Proposition \ref{proDg}: only one of the four possible $D$-graphs of Figs.~\ref{DGX921A}, \ref{DGX921B}, \ref{DGX921C} and \ref{DGX921D} can be  split  into the standard Coxeter-Dynkin graphs of types $A_3$ and $A_6$. \hfill $\Box$

\subsection{ Realization of Fig.~\ref{DGX921C}} 
\label{reala5a4a}

\begin{proposition}
The top $($of degree four$)$ homogeneous part 
of the polynomial 
\begin{equation}
f = \frac{1}{2} x^2 + x y ^3 -\frac{12\sqrt{3}}{5} x^2 y + 9 x^2 y^2 + \frac{4 \sqrt{3}}{15} x^3 - \frac{9}{5} x^3 y + \frac{3}{25} x^4
\label{ffx91c}
\end{equation} 
vanishes on exactly two lines in ${\mathbb R}^2$. 
The $j$-invariant of this homogeneous part is $-96$.
\end{proposition}

\noindent
{\it Proof.} The first statement is equivalent to saying that the cubic polynomial $$t^3 + 9 t^2 - \frac{9}{5} t + \frac{3}{25} $$ in one variable has only one real root. Being equal to $t^3 + 9 (t-\frac{1}{10})^2 + \frac{3}{100},$ it has no positive roots. It is also positive at $t=0$, and its derivative is negative there, hence it has only one negative root. \hfill $\Box$
\medskip

The function (\ref{ffx91c})
has a singularity of type $A_5$ at the origin. Indeed, its lower quasihomogeneous part with the weights $(3,1)$ is  $\frac{1}{2} x^2 + x y^3 \equiv \frac{1}{2}(x+y^3)^2 -\frac{1}{2}y^6$. 

In the coordinates $\tilde x = x - x_0$ and $\tilde y = y- y_0$ centered   at the point $ (x_0, y_0) = \left( - \frac{5 \sqrt{3}}{24},  \frac{\sqrt{3}}{8}\right) $
 the polynomial (\ref{ffx91c}) is  
\begin{equation}
-\frac{5}{3^2 \cdot 2^8} + \frac{3}{80} ( \tilde x - 5 \tilde y)^2 - \frac{\sqrt{3}}{120} ( \tilde x - 5 \tilde y) (7 \tilde x^2 - 82 \tilde x \tilde y - 5 \tilde y^2) + \tilde x \tilde y ^3 + 9 \tilde x^2 \tilde y^2 - \frac{9}{5} \tilde x^3 \tilde y + \frac{3}{25} \tilde x^4.
\end{equation}
In the coordinates $X \equiv \tilde x - 5 \tilde y$ and $\tilde y $ it appears as 
$$ -\frac{5}{3^2 \cdot 2^8} + \frac{3}{80} X^2 - \frac{\sqrt{3}}{120} X (- 240 \tilde y^2 + X (-12 \tilde y + 7 X)) + 
$$
$$
+ 80 \tilde y^4 + X \left(16 \tilde y^3 +X^2 \left(\frac{3}{25} X + \frac{3}{25} \tilde y\right)\right).
$$
The lower quasihomogeneous part of this function (minus the constant term) with the weights $(2, 1)$ is 
$$\frac{1}{80}\left(\sqrt{3}X + 80 \tilde y^2\right)^2 , $$
in particular it is degenerate. 
So its singularity at this point is of type $A_k,$ $k \geq 4$.
The sum of its Milnor number and the Milnor number of the critical point at the origin cannot exceed 9, so it is exactly of type $A_4$ and can be written in the form $$ -\frac{5}{3^2 \cdot 2^8} + \xi^2 + \eta^5$$ in some local coordinates. 

We can slightly perturb  the polynomial (\ref{ffx91c}) so that its critical point of type $A_5$ splits into three saddlepoints with the critical value 0 and two local minima with slightly lower values, and the critical point of type $A_4$ splits into two saddlepoints with the critical value $-\frac{5}{3^2 \cdot 2^8}$ and two slightly lower local minima.

\begin{figure}
\unitlength 0.5mm
\begin{picture}(110,120)
\put(90,0){\line(0,1){120}}
\bezier{500}(50,120)(76,63)(40,60)
\bezier{200}(40,60)(1,57)(0,55)
\bezier{50}(0,55)(1,53)(20,51)
\bezier{300}(20,51)(88,49)(90,40)
\bezier{200}(90,40)(100,10)(110,0)
\bezier{500}(60,120)(82,63)(40,54)
\bezier{500}(40,54)(86,53)(88,120)
\bezier{300}(94,0)(92,40)(100,0)
\end{picture} \qquad \qquad 
\begin{picture}(110,120)
\put(90,0){\line(0,1){120}}
\bezier{500}(50,120)(76,63)(40,60)
\bezier{200}(40,60)(1,57)(0,55)
\bezier{50}(0,55)(1,53)(20,51)
\bezier{300}(20,51)(108,48)(90,42)
\bezier{30}(90,42)(85,40)(90,37)
\bezier{200}(90,37)(94,35)(110,0)
\bezier{500}(60,120)(82,63)(44,54)
\bezier{500}(44,54)(86,51)(88,120)
\bezier{300}(94,0)(92,40)(100,0)
\bezier{40}(44,54)(38,56)(34,54)
\bezier{40}(44,54)(38,52)(34,54)
\bezier{40}(34,54)(26,56)(26,54)
\bezier{40}(34,54)(26,52)(26,54)
\end{picture}
\caption{Critical level sets of the polynomial (\ref{ffx91c}) and of its perturbation  realizing Fig.~\ref{DGX921C}}
\label{??C}
\end{figure}

\begin{proposition}
The D-graph of the obtained Morse polynomial  is shown in Fig.~\ref{DGX921C}.
\end{proposition}

\noindent
{\it Proof} is analogous to that of Proposition \ref{proDg}: only one of the D-graphs shown in Figs.~\ref{DGX921A}, \ref{DGX921B}, \ref{DGX921C} and \ref{DGX921D} splits into the standard Coxeter-Dynkin graphs of types $A_5$ and $A_4$. 
 \hfill $\Box$ 

\subsection{ Realization of the D-graph of Fig.~\ref{DGX921D}} 
\label{reale7a2}

\begin{lemma}
\label{lemma101}
The set $\{E_7+A_2\}$ of the space of degree four polynomials with the principal parts of type $X_9^1$ consists of four connected strata. Any polynomial from this set can be reduced to the form 
\begin{equation}
\pm x^4 +  t x^3 \pm x y^3, \quad t >0,
\label{nff}
\end{equation}
$($with all possible combinations of the signs$)$ by an orientation-preserving  affine transformation of ${\mathbb R}^2$. The $j$-invariant of the principal homogeneous part of $($\ref{nff}$)$ is  0.
\end{lemma}

\begin{figure}
\unitlength=0.5mm

\begin{picture}(100,80)
\put(60,0){\line(0,1){80}}
\bezier{800}(100,75)(63,50)(60,30)
\bezier{400}(60,30)(45,70)(25,35)
\bezier{500}(25,35)(15,20)(0,10)
\bezier{600}(35,30)(25,15)(0,5)
\bezier{400}(35,30)(37,25)(40,23)
\bezier{400}(40,23)(52,15)(58,0)
\end{picture}
\qquad \qquad \qquad
\begin{picture}(100,90)
\put(65,0){\line(0,1){90}}
\bezier{800}(0,20)(40,80)(60,80)
\bezier{400}(60,80)(80,80)(52,45)
\bezier{200}(52,45)(44,35)(48,45)
\bezier{300}(48,45)(60,60)(100,80)
\bezier{400}(42,0)(30,11)(23,21)
\bezier{200}(23,21)(22,26)(27,26)
\bezier{200}(27,26)(31,24)(25,16)
\bezier{400}(25,16)(17.3,4)(6.8,0)
\put(25.8,17){\circle*{1.5}}
\put(65,59.7){\circle*{1.5}}
\put(65,79){\circle*{1.5}}
\put(65,63.5){\circle*{1.5}}
\put(61,57){\circle*{1.5}}
\end{picture}
\caption{Level sets of perturbations towards the D-graph of Fig.~\ref{DGX921D}}
\label{DGX921Dd}
\end{figure}
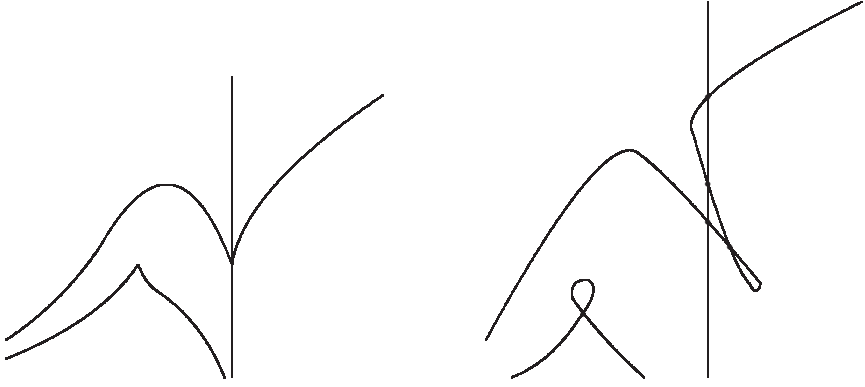

\noindent
{\it Proof.} The 3-jet of a degree four polynomial $\varphi$ at a point of type $E_7$ can be reduced to the form $x^3$ by such a transformation.  The coefficient at $y^4$ in the same coordinates will be equal to 0, and the coefficient at $x y^3$ not. By dilating the coordinate $y$,  the last coefficient can be made $\pm 1$. If $\varphi$ has another critical point of type $A_2$, then by a transformation of type $\tilde y = y+ \theta x$ (preserving the lower quasihomogeneous part of $\varphi$)  this point can be moved to the axis $\{y=0\}$. 
In the coordinates obtained, the function $\varphi$ has the form 
$$x^3 \pm x y^3 + a x^4 + b x^3 y + c x^2 y^2 .$$
Both partial derivatives and the Hessian determinant of $\Phi$ should vanish at the point $(\alpha, 0)$ for some $\alpha \neq 0$. This gives us a system of three equations on the parameters $a, b, c, \alpha$. It is easy to calculate that this system implies $b=c=0$, $a \alpha = - \frac{3}{4}$. So, we can choose $a$ arbitrarily; by dilating the coordinates the function obtained can be reduced to the form (\ref{nff}). 
The qualitative behavior of the level set $\{\varphi =0\}$ and of some component of the set $\{\varphi = -\frac{27}{256}\}$ of the polynomial $\varphi \equiv x^4 + x^3 - x y^3$ is shown in Fig.~\ref{DGX921Dd} left.

Different signs at the monomial $x^4$ in (\ref{nff}) provide different signs of the critical values at the point of type $A_2$, so the corresponding polynomials belong to different strata of the set $\{E_7+A_2\}$.  Also, the zero level set of any polynomial (\ref{nff}) defines an orientation of the plane, so the polynomials with equal signs at $x^4$ but different signs at $x y^3$ define opposite orientations. \hfill $\Box$ \medskip

Two critical points of the polynomial $\varphi$ can be perturbed independently in such a way that the level sets $\{\tilde \varphi =0 \}$ and $\{\tilde \varphi = -\frac{27}{256}\}$ of the new polynomial $\tilde \varphi$ behave topologically as shown in Fig.~\ref{DGX921Dd} right (with a different scaling than in Fig.~\ref{DGX921Dd} left), cf. page 16 in \cite{AC} or page 12 in \cite{Vcau}. In particular, this polynomial $\tilde \varphi$ has four local minima and five saddlepoints.

\begin{proposition}
The D-graph of the obtained Morse polynomial $\tilde \varphi$ is shown in Fig.~\ref{DGX921D}.
\end{proposition}

\noindent
{\it Proof.} By Theorem \ref{mthmdX91}, the D-graphs of Morse perturbations of $X_9^1$ singularities without local maxima can have only four shapes shown in Figs. \ref{DGX921A}, \ref{DGX921B}, \ref{DGX921C}, and \ref{DGX921D}. Only the last of these contains a subgraph equivalent to the standard Coxeter-Dynkin graph of type $E_7$, which necessarily occurs in the  D-graph of the polynomial $\tilde \varphi$. \hfill $\Box$

\subsection{Proof of Proposition \ref{refprop} for $X_9^1$}
\label{proofref1}

This proof is a variation of the proof for $X_9^+$ given in \S \ref{proofref0}. Again, we can assume that the polynomials $f$ and $\tilde f$ have the same critical values and systems of paths in ${\mathbb C}^1$ that define the vanishing cycles. We construct the paths in the parameter space $\Theta \cap {\mathbb R}^9$ of the versal deformation (\ref{vers0}) connecting them with some polynomials having a critical point of type $E_7$ with zero critical value and a point of type $A_2$ with the value $-\frac{27}{256}$; these paths have the same Lyashko--Looijenga projections to the space $\mbox{Sym}^9({\mathbb C}^1)$ of collections of critical values, and their almost final points should be generic Morse polynomial realizing the D-graph of Fig.~\ref{DGX921D}. 

Unlike \S~\ref{proofref0}, we cannot assume that the  polynomials corresponding to the endpoints of these paths coincide with each other: by Lemma \ref{lemma101} they can also be equal to two mutually symmetric polynomials $ x^4 + x^3 \pm x y^3$. On the other hand, if these endpoints are the same, then our paths just coincide, because the versal deformations of the $E_7$ and $A_2$ singularities have no real symmetries. For the same reason, if these endpoints are different, then our entire  paths are taken into each other by the same symmetry $f_\lambda(x, y) \leftrightarrow f_\lambda(x,-y)$ which takes these endpoints to each other. So, also the starting points of these paths are either coincident or symmetric. \hfill $\Box$

\subsection{Polynomials with two asymptotes and nine critical points, having both local minima and maxima}

\unitlength 1mm
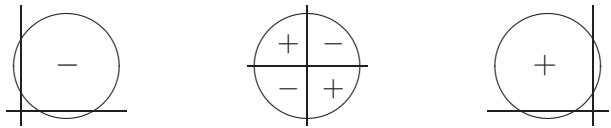
\begin{figure}
\unitlength 1.3mm
\begin{picture}(25, 20)
\bezier{200}(2,10)(2,18)(10,18)
\bezier{200}(10,18)(18,18)(18,10)
\bezier{200}(18,10)(18,2)(10,2)
\bezier{200}(10,2)(2,2)(2,10)
\put(0,3){\line(1,0){20}}
\put(3,0){\line(0,1){20}}
\put(10,10){\makebox(0,0)[cc]{\tiny $-$}}
\end{picture} \qquad \qquad
\begin{picture}(25,20)
\bezier{200}(2,10)(2,18)(10,18)
\bezier{200}(10,18)(18,18)(18,10)
\bezier{200}(18,10)(18,2)(10,2)
\bezier{200}(10,2)(2,2)(2,10)
\put(14,14){\makebox(0,0)[cc]{\tiny $-$}}
\put(6,6){\makebox(0,0)[cc]{\tiny $-$}}
\put(6,14){\makebox(0,0)[cc]{\tiny $+$}}
\put(14,6){\makebox(0,0)[cc]{\tiny $+$}}
\put(0,10){\line(1,0){20}}
\put(10,0){\line(0,1){20}}
\end{picture} \qquad \qquad
\begin{picture}(20,20)
\bezier{200}(2,10)(2,18)(10,18)
\bezier{200}(10,18)(18,18)(18,10)
\bezier{200}(18,10)(18,2)(10,2)
\bezier{200}(10,2)(2,2)(2,10)
\put(0,3){\line(1,0){20}}
\put(17,0){\line(0,1){20}}
\put(10,10){\makebox(0,0)[cc]{\tiny $+$}}
\end{picture}
\caption{Polynomials with two asymptotes, having both minima and maxima}
\label{x91nz}
\end{figure}

For any standard polynomial of class $X_9^1$ written in the normal form 
from Table \ref{t1},
three factors of this formula can be perturbed by the monomials of lower degrees in such a way that the zero set of the obtained polynomial will look as shown in one of the diagrams of Fig.~\ref{x91nz}. The numbers $(m_-, m_{\times}, m_+)$ of minima, saddlepoints and maxima of these polynomials are equal respectively to $(1, 5, 3), (2, 5, 2)$ and $(3, 5, 1)$.
The D-graph of the polynomial shown on the left (respectively, in the middle) side of this figure can be computed by the Gusein-Zade--A'Campo method and 
indeed is as shown in Fig.~\ref{X9011} left (respectively, right). 

Obviously, there are polynomials of each of these three types that are invariant under some reflections. Therefore, each of these types is represented by a single connected component of the space of Morse polynomials.

\subsection{Polynomials with two asymptotes and fewer than nine real critical points}
\label{x1rest}

By the proved part of Theorem \ref{mthmdX91}, the passports  $(m_-, m_{\times}, m_+)$ of these polynomials can only take values from the set $(0, 4, 3), $ $(3, 4, 0), $ $(1, 4, 2), $ $(2, 4, 1), $  $(0, 3, 2), (1, 3, 1), $ $(2, 3, 0), (0, 2, 1), $ $(1, 2, 0),$ and $(0, 1, 0)$. Moreover, there is exactly one virtual component with each of these values. Below we show that all these passports  (perhaps except the first two) can be realized by polynomials invariant under some reflections, so by Proposition \ref{refprop} each of these passports corresponds to  a single connected component of the space of Morse polynomials. The proof of the analogous fact for the passports $(3,4,0)$ and $(0,4,3)$ will be different.
\smallskip
 
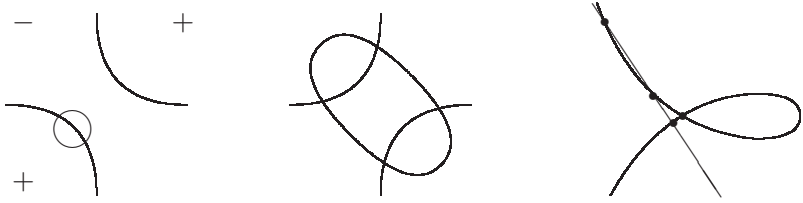
\begin{figure}
\unitlength 0.4mm
\begin{picture}(60,60)
\bezier{400}(30,0)(30,30)(0,30)
\bezier{400}(30,60)(30,30)(60,30)
\put(22,22){\circle{12}}
\put(55,55){\small $+$}
\put(2,2){\small $+$}
\put(2,55){\small $-$}
\end{picture} \qquad \quad
\begin{picture}(60,60)
\bezier{400}(30,60)(30,30)(0,30)
\bezier{400}(30,0)(30,30)(60,30)
\bezier{300}(40,40)(20,60)(10,50)
\bezier{300}(40,40)(60,20)(50,10)
\bezier{300}(20,20)(0,40)(10,50)
\bezier{300}(20,20)(40,0)(50,10)
\end{picture} \qquad \quad
\unitlength 0.3 mm
\begin{picture}(100,100)
\bezier{500}(16,0)(40,43)(80,45)
\bezier{300}(80,45)(100,45)(100,35)
\bezier{300}(100,35)(100,25)(80,25)
\bezier{800}(10,85)(30,30)(80,25)
\put(8,85){\line(2,-3){57}}
\put(13.5,77){\circle*{3}}
\put(34.9,44.1){\circle*{3}}
\put(43.8,32){\circle*{3}}
\put(48,35){\circle*{3}}
\end{picture}
\caption{Morse functions of $X_9^1$ type with few critical points}
\label{realx91l}
\end{figure}

The passport $(0,1,0)$ is realized by the polynomial 
$ x^4 + x^2 -y^4 -y^2$. The passports $(1,2,0)$ and $(0,2,1)$ are realized by the  polynomials $x^4 - y^4 \pm (x^2 + y^2)$.

\begin{proposition}
The passport $(1, 3, 1)$ is realized by the polynomial
\begin{equation}(x y-1)((x+1)^2 +(y+1)^2 -\varepsilon), \ \varepsilon > 0, 
\label{tq}
\end{equation}
 see Fig.~\ref{realx91l} left. 
\end{proposition}

\noindent
{\it Proof.} This polynomial has at most seven real critical points, because  two factors of (\ref{tq}) have two other intersection points in the complex domain. 
It is clear from the diagram that it has three critical points on the line $\{x=y\}$: a minimum, a maximum and a saddlepoint, and additionally two saddlepoints with zero critical value. If there are more critical points, then there are exactly two of them, they are symmetric with respect to the line $\{x=y\}$ and have the same Morse index. Depending on this index, the 
Euler number $m_- - m_\times + m_+$ of all real critical points will be either $-3$ or $+1$, which contradicts the index of the principal part of $X_9^1$ singularities at infinity.
\hfill $\Box$

\begin{proposition}
The passport $(2, 3, 0)$ is realized by an arbitrary polynomial 
\begin{equation}  f \equiv -y^4+y^2(x+2\varepsilon)+x^4 +\varepsilon x^3 -\varepsilon^2 x^2 - \varepsilon^3 x
\label{nea}
\end{equation}
$($which is even in the variable $y)$ with sufficiently small positive $\varepsilon$, and the passport $(0,3,2)$ is realized by the same polynomial multiplied by $-1$.
\end{proposition}

\noindent
{\it Proof.} The equation $f'_y=0$ provides two choices $y=0$ and $2y^2=x+2\varepsilon$. The equality $y=0$ substituted into the condition $f'_x=0$ gives us an equation on $x$ with three real roots $x = -\varepsilon, x= \frac{1\pm\sqrt{17}}{8} \varepsilon$. At all these roots the coefficient $x+2\varepsilon$ of the term $y^2$ in (\ref{nea}) is positive, so they are two local minima and one saddlepoint. 

Substituting the second condition gives the equation $$4 x^3 + 3 \varepsilon x^2 +\left(\frac{1}{2}-2\varepsilon \right) x +\left(\frac{\varepsilon}{2} - \varepsilon^3\right) = 0. $$
For $\varepsilon =0 $ it has two non-real roots, hence the same will hold for neighboring values of $\varepsilon$. Each of these roots corresponds to two critical points of the polynomial (\ref{nea}) which thus 
has at least four non-real critical points. The Euler number $m_- - m_\times + m_+$ of all its real critical points should be equal to $-1$, therefore the remaining two critical points are real and are saddlepoints. \hfill $\Box$

\begin{proposition}
The passport $(1,4,2)$ is realized by polynomials $$(x^2 + 2 A x y + y^2 -2)(x y + 1)$$ with $A \in (0,1)$, and the passport $(2,4,1)$ is realized by  polynomials $(x^2 - 2A x y + y^2 -2) (x y-1)$. 
\end{proposition}

\noindent
{\it Proof.} Each of these polynomials has seven  obvious critical points, see Fig.~\ref{realx91l} center. Two additional  real critical points would be symmetric with respect to the line $\{x=y\}$, in particular they would have equal Morse indices, which contradicts the balance of Morse indices of all critical points. \hfill $\Box$
\medskip

Finally, to realize the passport $(3,4,0)$ we first take the polynomial $x^3y + x y^3 + x^3 $, which has a singularity of class $E_7$ at the origin and no other real critical points, and then use 
a Morse perturbation of the latter singularity shown in Fig.~\ref{realx91l} right, see also p.~16 in \cite{AC}. Denote the Morse polynomial obtained by $f_\lambda$.

\begin{proposition}
The Morse function $f_\lambda(x,y)$ just constructed and its mirror image $f_\lambda(x, -y)$ belong to the same connected component of the space of Morse polynomials of degree four with non-discriminantal principal parts.
\end{proposition}

\noindent
\begin{table}
\caption{Disorienting 1-cycle in the formal graph}
\label{virtcycle}
\begin{tabular}{||ccccccc|cc|}
\hline
-2  & 0  & 0  & 0  & 0  & 1  & 0  & 0  & 0 \\
0 & -2   & 0 &  0  & 1  & 0 &  1  & 1  & -1 \\
 0  & 0 & -2  & 1  & 0 &  0 &  1 &  0 &  0 \\
 0  & 0  & 1 & -2  & 0  & 0  & 0  & 0 &  0 \\
 0 &  1 &  0 &  0 & -2 & 0 & 0 & -1 & 0 \\
 1 & 0 & 0  & 0 & 0 & -2 & 0 & -1 & -1 \\
 0  & 1  & 1 & 0 & 0 & 0 & -2 & -1 & 0 \\
 0  & 1 &  0  & 0 & -1 & -1 & -1&  -2 & -1 \\
 0  & -1 &  0 & 0 & 0 & -1 & 0 & -1 & -2 \\
\hline
0 & 0 & 0 & -1 & -1 & -1 & 0 & -1 & -1 \\
\hline
0 & 0 & 0 & 1 & 1 & 1 & 1 & & \\
\hline
\end{tabular} \qquad 
\begin{tabular}{||ccccccc|cc|}
\hline
-2  & 0  & 0  & 0  & 1  & 0  & 1  & 1  & -1 \\
0 & -2   & 0 &  0  & 0  & 1 &  0  & 0  & 0 \\
 0  & 0 & -2  & 1  & 0 &  0 &  1 &  0 &  0 \\
 0  & 0  & 1 & -2  & 0  & 0  & 0  & 0 &  0 \\
 1 &  0 &  0 &  0 & -2 & 0 & 0 & -1 & 0 \\
 0 & 1 & 0  & 0 & 0 & -2 & 0 & -1 & -1 \\
 1  & 0  & 1 & 0 & 0 & 0 & -2 & -1 & 0 \\
 1  & 0 &  0  & 0 & -1 & -1 & -1&  -2 & -1 \\
 -1  & 0 &  0 & 0 & 0 & -1 & 0 & -1 & -2 \\
\hline
0 & 0 & 0 & -1 & -1 & -1 & 0 & -1 & -1 \\
\hline
0 & 0 & 0 & 1 & 1 & 1 & 1 & & \\
\hline
\end{tabular}
\end{table}
\noindent
{\it Proof.} For any strictly Morse polynomial $f_\lambda$ of this component, consider the triangle in ${\mathbb R}^2$ formed by its three local minima ordered by their critical values, and 
define the sign of $f_\lambda$ as the orientation of this triangle.
Obviously, the function $f_\lambda$ and its mirror image have opposite signs. A chain of elementary surgeries not changing the isotopy class of the Morse polynomials (i.e., not containing surgeries $s1$ and $s3$) has different signs of its endpoints if and only if the surgeries of type $s2$, at which the critical values at the points of minima  collide, occur an odd number of times along it. 
Define the {\em disorienting cocycle} of a virtual component as the 1-cocycle  of this graph taking values in ${\mathbb Z}_2$ and equal to the parity of the number of such $s2$-surgeries along the cycles in this graph. Our component of the space of Morse polynomials is chiral if and only if this cocycle of the associated virtual component is trivial.

The list of virtual Morse functions with the passport $(m_-, m_\times, m_+) = (3, 4, 0)$ contains two elements shown in Table \ref{virtcycle}. They can be obtained from each other by a single surgery of type $s4$ with non-real critical values crossing the real line between the sixth and the seventh real critical values. Also, these two virtual Morse functions are connected by the chain of two surgeries of type $s2$ permuting the fifth and the sixth critical values  and the first and the second critical values, respectively. Only one of these three surgeries changes the parity of the order of the minima, so the disorienting cocycle is non-trivial on this cycle of length three, and our isotopy class is achiral.
\hfill $\Box$

\begin{remark} \rm
This cycle was also found by a computer search on the graph.
\end{remark}

\section{Construction of polynomials with four real asymptotes ($X_9^2$)} 
\label{reaax2}

\subsection{Realization of the  D-graph of Fig.~\ref{DGX92A}} 
\label{reala5d4m}

 The polynomial
\begin{equation}
F=  x^4 - 2 \sqrt{2} x^2 y^2 + y^4 -
 \frac{4}{3} x^3 + 2 \sqrt{2}  x y^2  
\label{realx92}
\end{equation}
 has a singularity of type $D^-_4$ at the origin.
It also has a critical point with the coordinates $(1, 0)$. In the local coordinates $(\hat x = x- 1 , y) $ it is equal to 
\begin{equation}
-\frac{1}{3} + 2 \hat x^2 + \frac{8}{3}\hat x^3 -2\sqrt{2} \hat x y^2 + \hat x^4 - 2\sqrt{2}\hat x^2 y^2 + y^4 ,
\label{redd}
\end{equation}
so the lowest quasihomogeneous part of the polynomial $F+ \frac{1}{3}$ (with weights $(2, 1)$) is
$$ 2 \hat x^2 - 2\sqrt{2} \hat x y^2 + y^4 \equiv (\sqrt{2} \hat x - y^2)^2$$
In the coordinates $z = \sqrt{2}x - y^2$ and $y$,  the lower quasihomogeneous part of $F+\frac{1}{3}$ with the weights $(3,1)$ is  $z^2 - \frac{\sqrt{2}}{3} y^6 , $ so the entire function $F$ has a singularity of class $A_5$ at this point.

\begin{figure}
\unitlength 0.4 mm
\begin{picture}(120,120)
\put(0,0){\line(1,0){120}}
\put(0,0){\line(0,1){120}}
\put(120,120){\line(-1,0){120}}
\put(120,120){\line(0,-1){120}}
\bezier{800}(20,0)(63,60)(20,120)
\bezier{700}(0,100)(70,60)(80,0)
\bezier{700}(0,20)(70,60)(80,120)
\bezier{800}(95,120)(80,65)(71,60)
\bezier{600}(71,60)(60,50)(120,20)
\bezier{600}(71,60)(60,70)(120,100)
\bezier{800}(95,0)(80,55)(71,60)
\put(4,4){\small $-$}
\put(4,112){\small $-$}
\put(4,58){\small $+$}
\put(110,4){\small $-$}
\put(110,112){\small $-$}
\put(57,4){\small $+$}
\put(57,112){\small $+$}
\put(41.3,50){\circle*{1.5}}
\put(41.3,70){\circle*{1.5}}
\put(51,60){\circle*{1.5}}
\put(71,60){\circle*{1.5}}
\put(81.8,78){\circle*{1.5}}
\put(81.8,42){\circle*{1.5}}
\end{picture}
\caption{Level sets for a function realizing Fig.~\ref{DGX92A}}
\label{miss92alevel}
\end{figure}
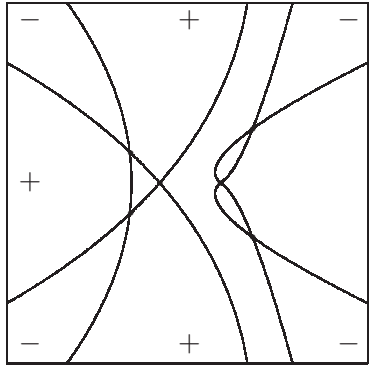

Perturbing independently these two critical points as indicated (up to a sign)  on pp. 13--14 of \cite{AC}, we obtain a polynomial that has two level sets 
looking topologically as shown in Fig.~\ref{miss92alevel}. Namely, three curves bounding the curvilinear triangle on the left belong to the set where the function is equal to 0;
the curves on the right belong to a similar set with the value $-\frac{1}{3} $.

By the construction, the D-graph of the obtained polynomial can be decomposed into the canonical Coxeter-Dynkin graphs of types $D_4$ and $A_5$ by removing some of its edges, so it can only be  the D-graph of Fig.~\ref{DGX92A}.

The $j$-invariant of the principal part of (\ref{realx92}) is $\frac{125}{27}$.

  \subsection{Realization of Fig.~\ref{DGX92B}} 
\label{reald6a3}

\begin{proposition}
\label{prol3}
A polynomial with the principal part of type $X_9^2$ has a critical point of class $D_6^-$ with zero critical value and a critical point of class $A_3$ if and only if  in appropriate affine coordinates it is equal to 
 \begin{equation}
f = x^2 y - x y^3 + \frac{a}{3} x^3 + \frac{3a^3}{4}  x^4 + 3a^2 x^3 y + \frac{3a}{2} x^2 y^2 , a \neq 0.
\label{nfff0}
\end{equation}
The $j$-invariant of the degree four homogeneous part of any such polynomial is  5.
\end{proposition}

\noindent
{\it Proof.} {\bf ``If''}. The principal (lower) quasihomogeneous part of (\ref{nfff0}) with the weights (2, 1) is $x^2 y- x y^3 \equiv y (x-\frac{y^2}{2})^2 -\frac{y^5}{4}$, so it has the normal form $D_6^-$ from Table \ref{t1} in the coordinates $\sqrt[5]{2}(x-\frac{y^2}{2})$ and $\frac{y}{\sqrt[5]{4}}$. According to \S 12 of \cite{AVGZ}, this entire function can then be reduced to this normal form by a local diffeomorphism with the unit Jacobian matrix at the origin.

The other critical point of (\ref{nfff0}) is $(x,y) = \left(-\frac{1}{3a^2},0\right)$. In the local coordinates $\tilde x \equiv a(x + \frac{1}{3a^2})+y$ and $y$, this function is  
$$ \frac{-1}{2^2 3^4 a^5} + \frac{1}{6a^3} \tilde x^2 +\tilde x\left(-\frac{2}{3a}\tilde x^2+\frac{1}{a^2}y^2\right) +\frac{1}{4a}y^4 + \tilde x \left(\frac{2}{a} y^3 + \tilde x \left(-\frac{3}{a}y^2 +\frac{3}{4a}\tilde x ^2\right)\right).$$
Its principal quasihomogeneous part  is $$\frac{1}{6a^3}\tilde x^2 + \frac{1}{a^2}\tilde x y^2 + \frac{1}{4a}y^4 \equiv \frac{1}{6a^3}\left((\tilde x+3a y^2)^2 - \frac{15}{2}y^4\right), $$
so this singularity is of type $A_3$.
The qualitative behavior of some components of the critical level sets of this polynomial is shown in Fig.~\ref{nfnf} left.
 \smallskip

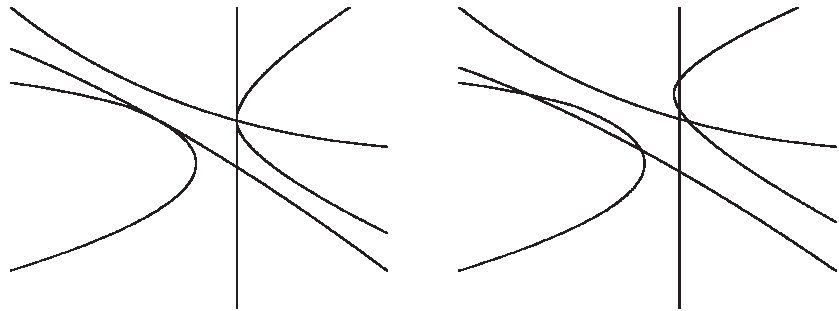
\begin{figure}
\unitlength 0.5mm
\begin{picture}(100,80)
\put(60,0){\line(0,1){80}}
\bezier{400}(90,80)(60,60)(60,50)
\bezier{400}(60,50)(60,40)(100,20)
\bezier{800}(0,80)(40,48)(100,43)
\bezier{800}(0,69)(60,54)(57,0)
\bezier{400}(0,60)(15,58)(25,54)
\bezier{400}(25,54)(80,35)(0,10)
\end{picture} \qquad
\begin{picture}(100,80)
\put(58.5,0){\line(0,1){80}}
\bezier{400}(90,80)(57,65)(57,57)
\bezier{400}(57,57)(57,47)(100,23)
\bezier{800}(0,80)(40,47)(100,43)
\bezier{800}(0,67)(55,53)(56,0)
\bezier{400}(0,60)(15,58)(25,54)
\bezier{400}(25,54)(80,35)(0,10)
\end{picture}
\caption{Critical level sets for a polynomial of class $D_6^- + A_3$ and its morsification}
\label{nfnf}
\end{figure}

{\bf ``Only if''.} In some affine coordinates centered at the $D_6^-$ point, the 3-jet of a polynomial of degree four is equal to $x^2y$, the coefficient at $y^4$ is  0 and the coefficient at $x y^3$ is not. By dilating the coordinates (maybe with negative coefficients) the sum of these two monomials can be made equal to $x^2y-x y^3$. Suppose that our polynomial also has a point of type $A_3$. By  substituting $\tilde y = y + \theta \cdot x$ we can move this point to the line $\{y=0\}$. 
In the coordinates $x$ and $\tilde y$, the obtained polynomial has the form 
 \begin{equation}
x^2 \tilde y - x \tilde y^3 + \frac{a}{3} x^3 + \frac{b}{4} x^4 + c x^3 \tilde y + \frac{d}{2} x^2 \tilde y^2. 
\label{nfff}
\end{equation}

The rest  of the proof follows the scheme of the proof  of Lemma \ref{lem77}: in its notation, the solutions of the resulting system of equations are parameterized by the coefficient $a$, which takes any non-zero value, and the remaining coefficients are expressed via it as $b=3a^3,$ $c=3a^2$, $d=3a$, $\alpha = - \frac{1}{3 a^2}$, $\beta = - \frac{1}{a}$, $\gamma = -3$.  \hfill $\Box$

\begin{corollary}
The set of degree four polynomials with the principal homogeneous parts of type $X_9^2$, having critical points of types $D_6^-$ and $A_3$, consists of four connected strata.
\end{corollary}

\noindent
{\it Proof.} There are at most four such strata, because the group of affine transformations by which we can reduce such polynomials to the form (\ref{nfff0}) consists of two connected components, and the set of polynomials of this form also consists of two components depending on the sign of the parameter $a$. The sign of the critical value of any polynomial of the form (\ref{nfff0}) at its point of type $A_3$ is opposite to the sign of $a$, therefore such polynomials with different signs of $a$ belong to different components. 
Also, for any value of $a$ any polynomial of the form (\ref{nfff0}) in some coordinates uniquely defines those coordinates. Therefore, the orientation of the plane defined by these coordinates is constant along any connected component of the set of such polynomials, while the reflection of the plane changes this orientation and the component. \hfill $\Box$
\medskip 

Perturbing any  polynomial of the form (\ref{nfff0}) with a positive value of the parameter  $a$ at  two critical points as shown (up to a sign) in pp. 13--14 of \cite{AC}, we obtain a polynomial with six saddlepoints and three local minima, see Fig.~\ref{nfnf} right. 

\begin{proposition}
The D-graph of the obtained polynomial
is as shown in Fig.~\ref{DGX92B}. 
\end{proposition}

\noindent
{\it Proof.} By Theorem \ref{mthmdX92}, there are only three possible  D-graphs of degree four polynomials with nine real critical points, no local maxima, and the principal part vanishing on four different real lines. We do not have the case shown in Fig.~\ref{DGX92A} because it does not contain a subgraph of type $D_6^-$. We also do not have the case shown in Fig.~\ref{DGX92C} which cannot be split into two subgraphs of types $D_6^-$ and $A_3$. \hfill $\Box$ \medskip

The graph of Fig.~\ref{DGX92B} has no symmetries, so analogous to Proposition \ref{achir} this construction realizes two different (mutually symmetric)  isotopy classes of Morse polynomials.
 
\subsection{Realization of Fig.~\ref{DGX92C}} 
\label{realx2c}

The homogeneous degree four polynomial 
$$x^4 - 6 x^2y^2 + y^4 \equiv \mbox{Re} (x+ i y)^4 $$
obviously vanishes on four lines and splits into the product of two real quadratic forms
\begin{equation}
\label{four}
\left(x^2 - (3+\sqrt{8})y^2\right)\left(x^2 - (3-\sqrt{8})y^2\right).
\end{equation}

\unitlength 0.4mm
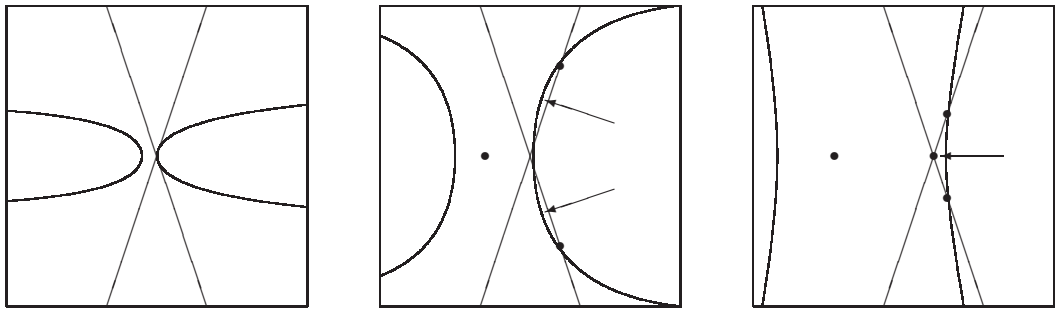
\begin{figure}
\begin{picture}(100,100)
\put(0,0){\line(1,0){100}}
\put(0,0){\line(0,1){100}}
\put(100,100){\line(-1,0){100}}
\put(100,100){\line(0,-1){100}}
\put(70,50){\line(1,3){16.7}}\put(70,50){\line(1,-3){16.7}}\put(70,50){\line(-1,3){16.7}}\put(70,50){\line(-1,-3){16.7}}
\bezier{700}(100,68)(70.1,62)(70.1,50)
\bezier{700}(100,32)(70.1,38)(70.1,50)
\bezier{700}(0,65)(25,60)(25,50)
\bezier{700}(0,35)(25,40)(25,50)
\put(45,50){\circle*{2.5}}
\end{picture} \qquad
\begin{picture}(100,100)
\put(0,0){\line(1,0){100}}
\put(0,0){\line(0,1){100}}
\put(100,100){\line(-1,0){100}}
\put(100,100){\line(0,-1){100}}
\put(60,50){\line(1,3){16.7}}
\put(60,50){\line(1,-3){16.7}}
\put(60,50){\line(-1,3){16.7}}
\put(60,50){\line(-1,-3){16.7}}
\bezier{700}(100,100)(61,85)(61,50)
\bezier{700}(100,0)(61,15)(61,50)
\put(5,50){\circle*{2.5}}
\put(69.5,78){\circle*{2.5}}
\put(69.5,22){\circle*{2.5}}
\put(88,60){\vector(-3,1){23}}
\put(88,40){\vector(-3,-1){23}}
\end{picture} \qquad
\begin{picture}(100,100)
\put(0,0){\line(1,0){100}}
\put(0,0){\line(0,1){100}}
\put(100,100){\line(-1,0){100}}
\put(100,100){\line(0,-1){100}}
\put(60,50){\line(1,3){16.7}}
\put(60,50){\line(1,-3){16.7}}
\put(60,50){\line(-1,3){16.7}}
\put(60,50){\line(-1,-3){16.7}}
\bezier{700}(70,100)(64,65)(64,50)
\bezier{700}(70,0)(64,35)(64,50)
\put(64.5,64){\circle*{2.5}}
\put(64.5,36){\circle*{2.5}}
\put(60,50){\circle*{2.5}}
\put(83.3,50){\vector(-1,0){21}}
\end{picture}
\caption{Realization of Fig.~\ref{DGX92C}}
\label{realX92C}
\end{figure}
\noindent
Let us keep the second factor of (\ref{four}) and perturb the first one as 
$(x-\lambda)^2 - (3+\sqrt{8})y^2 + \varepsilon$, $\varepsilon>0$.
Then for some choices of $\varepsilon $ and $\lambda$ the zero level set of the obtained function $F$ will look as shown (in different scalings) in three diagrams of Fig.~\ref{realX92C}. This zero set has five self-intersection points and bounds three domains of negative values of $F$, each of which contains a local minimum of $F$. In addition, there is a saddlepoint in an unbounded domain of the set where $F<0$ whose value is below all other eight critical values. 

\begin{proposition}
The  D-graph of the resulting function $F$ is as shown in Fig.~\ref{DGX92C}.
\end{proposition}

\noindent
{\it Proof.} The intersection indices of eight cycles which vanish either at the local minima of the polynomial $F$ or at the self-intersection points of its zero level set can be computed by the Gusein-Zade--A'Campo method \cite{AC}, \cite{GZ} and are described by the extended Coxeter-Dynkin graph  \ 
\begin{picture}(60,10)
\put(0,0){\line(1,0){60}}
\put(0,0){\circle*{1.5}}\put(10,0){\circle*{1.5}}\put(20,0){\circle*{1.5}}\put(30,0){\circle*{1.5}}\put(40,0){\circle*{1.5}}\put(50,0){\circle*{1.5}}\put(60,0){\circle*{1.5}}\put(30,10){\circle*{1.5}}
\put(30,0){\line(0,1){10}}
\end{picture} \ \ of type $\tilde E_7$.
Neither of the two competing  D-graphs shown in Figs. \ref{DGX92A} and \ref{DGX92B} contains such a subgraph. \hfill $\Box$

\subsection{Proof of Proposition \ref{refprop} for $X_9^2$} 
\label{proofref2}

This proof almost repeats the one from \S~\ref{proofref1} and is based on Proposition \ref{prol3}. Namely, we again draw a path $I$ connecting the polynomial $f$ and the polynomial (\ref{nfff0}) with $a=1$.
In this case, the real versal deformations of the two components $D_6^-$ and $A_3$ of the bisingularity $\{D_6^- + A_3\}$ have real symmetries of order 2, but neither of the corresponding symmetries of the Dynkin subgraphs of the D-graph of Fig.~\ref{DGX92B} can be continued to a symmetry of the entire graph. Therefore, the path $\tilde I$ defined as before  approaches either the same polynomial (\ref{nfff0}) or  the mirror image of (\ref{nfff0}). In the first case it coincides with $I$ in a neighborhood of its endpoint and hence everywhere, so that $\tilde f = f$; in the second case it is the mirror image of $I$, and its starting point $\tilde f$ is the mirror image of $f$. \hfill $\Box$

\subsection{Realization of Morse polynomials with four asymptotes that have both local minima and maxima or fewer than nine real critical points}

It is easy to see that the passport  (2, 6, 1) is realized by the polynomial $$(x-1)(y-1)(x^2-3 x y + y^2), $$ and the passport  (1, 6, 2) by minus this polynomial. 

The polynomial $$x^4 -2Ax^2y^2 + y^4 + 2 (A^2-1)x^2, \quad A>1, $$
has four saddlepoints with the coordinates $(x,y) = (\pm 1, \pm \sqrt{A})$, one critical point of type $A_3$ at the origin and no other critical points. By adding the monomial $\varepsilon y^2$ (respectively, $-\varepsilon y^2$) with sufficiently small $\varepsilon >0$, we do not change the type of the saddlepoints, but make the $A_3$ point to a single real minimum point (respectively, to two minima and one saddlepoint), thus realizing the passport (1, 4, 0) (respectively, (2, 5, 0)). The same polynomials multiplied by $-1$ realize the passports (0, 4, 1) and (0, 5, 2).   

The polynomial $$x^4 - 2 A x^2y^2 + y^4 + 2  x^2 -2 y^2, \quad A>1, $$
has exactly three saddlepoints with the coordinates $x=0, $ $y = 0, 1$ or $-1$,
thus realizing the passport (0, 3, 0).

\begin{proposition}
The passport $(1, 5, 1)$ is realized by the polynomials 
\begin{equation}
   x^4 - 2A x^2y^2 + y^4 + \frac{4}{3}x^3 + 4  \varepsilon x y^2 - 4\varepsilon^4 x, \quad  A>1, 
\label{nif}
\end{equation}
 for sufficiently small $\varepsilon>0$.
\end{proposition}

\noindent
{\it Proof.} 
 The main part $  x^4 - 2A x^2y^2 + y^4 + \frac{4}{3} x^3  $ of (\ref{nif}) has three saddlepoints 
with the coordinates $(-1,0)$ and $\left(\frac{1}{A^2-1}, \frac{\pm \sqrt{A}}{A^2-1} \right), $ and   a singularity of type $E_6$ at the origin. 
The perturbing small terms \ $4\varepsilon x y^2 -4\varepsilon^4 x$ \ 
do not change the type of saddlepoints, but the resulting function (\ref{nif}) has two critical points on the line $\{y=0\},$ arising from the perturbation of the $E_6$ point. Their $x$ coordinates are  the roots of the equation $x^3 + x^2-\varepsilon^4$ equal to $\pm \varepsilon^2 + O(\varepsilon^3)$ (while the third root approximately equal to $-1$ corresponds to one of the three saddlepoints fixed previously). The Hessian matrix of (\ref{nif}) is  $\begin{pmatrix}12 x^2 + 8 x & 0\\
0 & -4Ax^2 + 8 \varepsilon x
\end{pmatrix}$
on the line $\{y=0\}$, in particular it
is positive definite at one of these two critical points and is negative definite at the other. The remaining real critical points of (\ref{nif}) (in addition to these two and to three saddlepoints inherited from the unperturbed polynomial) appear in at most two pairs symmetric with respect to the line $\{y=0\}$, in particular having equal Morse indices. By the index theorem, the number of saddlepoints among them is by 2 greater than the number of local extrema. This is possible only if there is exactly one pair of such critical points, and they are saddlepoints. \hfill $\Box$
\medskip

All  polynomials constructed in this subsection are invariant under certain reflections, so by Proposition \ref{refprop} the corresponding virtual components are presented by single components of the set of Morse polynomials.

\section{Proof of Theorem \ref{tabadjp}}
\label{adj}

\subsection{Adjacency $\{A_5 + D_4^+\} \rightsquigarrow \{X_9^+\}$}
\label{adj1}

The approximation of an $X_9^+$ singularity by the stratum $\{A_5+D_4^+\}$ is realized by the family of polynomials 
\begin{equation}
\label{bbb}
\frac{x^4}{2} + x^2y^2 + \frac{y^4}{4} + t \left(x^2 y + \frac{y^3}{3}\right), \ t \neq 0.
\end{equation}
The critical point of every polynomial in this family at the origin is of type $D_4^+$. Also, the Taylor decomposition of (\ref{bbb}) at the point $(x,y) = (0, -t)$ (in the coordinates $x$ and $\tilde y = y+ t$) is  
$$-\frac{t^4}{12} + \frac{1}{2}(t \tilde y - x^2)^2 - \frac{2}{3}t \tilde y^3 + x^2 \tilde y^2 + \frac{1}{4} \tilde y^4. $$ In the coordinates $x$ and $Y \equiv t \tilde y - x^2$, the lower quasihomogeneous part with the weights $(1,3)$ of its nonconstant part is $$\frac{1}{2} Y^2 + \frac{1}{3t^2} x^6, $$ so this  is a critical point of type $A_5$.

The $j$-invariant of the principal homogeneous part of (\ref{bbb}) is $\left(\frac{5}{3}\right)^3$.

\subsection{Adjacency  $\{A_5 + D_4^+\} \rightsquigarrow \{X_9^1\}$}
\label{adj3}
It is realized by the family of polynomials
\begin{equation}\frac{t^2}{2} x^2 + x y^3 + \frac{2t}{3} x^3 + \frac{1}{4} x^4 .
\label{x1add}
\end{equation}
Indeed, any such polynomial has a singularity of type $A_5$ at the origin. In addition, its Taylor expansion at the point $(-t  , 0)$ in the coordinates $\tilde x =x+t$ and $y$ is  
$$ \frac{t^4}{12}  - \frac{t}{3}\tilde x^3 - t y^3 + \tilde x y^3 + \frac{1}{4} \tilde x^4  $$
with the lower homogeneous part $-\frac{t}{3} \tilde x^3 - t y^3$ of type $D_4^+$.

The $j$-invariant of the principal homogeneous part of (\ref{x1add}) is  0.

\subsection{Adjacency $\{E_6 + A_3\} \rightsquigarrow \{X_9^1\}$}
\label{adj38}
It is provided by the family of polynomials $$-x^4 + y^4 + t x^3 $$
with the $j$-invariant equal to 1.

\subsection{Obstructions to adjacencies}
\label{adj2}

Many approximations of the strata \{$X_9^{\ast}$\} by bisingularities are forbidden by the Euler characteristic consideration: the sum of the local indices of the gradient vector field at all real critical points of any function of class $X_9^{\ast}$ is equal to the similar index of the principal (degree four) homogeneous part of this function. This index is equal to $1$ for $X_9^{\pm}$, to $-1$ for $X_9^1$, to $-3$ for $X_9^2$, to 0 for $A_{2k}$, $D_{2k}^+$ and $E_6$, to $-2$ for $D_{2k}^-$, to $-1$ for $D_{2k-1}$ and $E_7$, and to $1$ or $-1$ for $A_{2k-1}$.  The cases where the adjacency is forbidden for this reason are indicated by the word ``(index)'' in Table \ref{tabadj}. 

There are no approximations of the $X_9$ singularities by the $D_7$ stratum (even in the complex domain), see \cite{siersma}; so we get ``No''s in the line for $D_7 + A_2$. Moreover, degree four real polynomials with non-discriminantal principal parts cannot have singularities of type $D_6^+$. Indeed, the degree three homogeneous part of any critical point of class $D_k,$ $k>4$, can be reduced to the form $x^2 y $ by a choice of affine coordinates. If in these coordinates the coefficient at the monomial $y^4$ of the degree four polynomial is not zero, then the singularity is of type $D_5$. If this coefficient is 0 but the coefficient at $x y^3$ is not, then the singularity is of type $D_6^-$, and if both these coefficients are equal to zero then it is non-isolated. For this reason, also the row of Table~\ref{tabadj} for $D_6^+ + A_3$ contains only negations.

The functions with two critical points of types $D_5$ and $D_4$ are forbidden by Bezout's theorem: all partial derivatives of a fourth degree polynomial should vanish on the line connecting  such two points. 

Finally, the negations marked in Table~\ref{tabadj} by ``(D-graph)'' are related to the following remark (cf. \cite{lyashko}, \cite{Jaw2}).

\begin{proposition}
\label{pro18}
If a fourth degree  polynomial with a principal part of type $X_9^{\ast}$ has two real critical points of types $\Xi$ and $\tilde \Xi$ with $\mu(\Xi)+\mu(\tilde \Xi)=9$, then the set of vertices of one of the D-graphs of this type  $X_9^{\ast}$  can be divided into two subsets of cardinality $\mu(\Xi)$ and $\mu(\tilde \Xi)$ in such a way that 

a$)$ all edges connecting vertices of different subsets are directed from  vertices of one of  these subsets to vertices of the other, and

b$)$ the edges of our D-graph, both vertices of which belong to one of these subsets, form a Dynkin graph describing the intersection form of vanishing cycles of the corresponding simple singularity $\Xi$ or $\tilde \Xi$ $($for some choice of bases of these vanishing cycles$)$.  
\end{proposition}

\noindent
{\it Proof.} Consider the polynomial $f_\tau$ with very small $\tau \neq 0$ from a family that realizes adjacency, and 
a Morse perturbation of $f_\tau$ that splits its two critical points as in the Gusein-Zade--A'Campo method (keeping all critical values at the saddlepoints  equal to the critical values of the corresponding initial critical point). 
The D-graph of the obtained perturbation satisfies the conditions of the proposition. \hfill $\Box$
\medskip

We know from Theorems \ref{mthmdX90}--\ref{mthmdX92} the complete lists of D-graphs of all Morse functions with $X_9^\ast$  principal parts and can check whether they admit such decompositions. The indication ``No (D-graph)'' in a cell of Table \ref{tabadj} means that the condition of Proposition \ref{pro18} is not satisfied for the corresponding classes $\Xi$, $\tilde \Xi$ and $X_9^*$. 


}

\end{document}